\documentclass[10pt]{article} 
\usepackage{
graphicx,subfigure,epstopdf,animate,latexsym,amsmath,amssymb,amsthm,amsfonts,bm,tikz,booktabs,ltxtable,color,tabularx,multirow,float,cases,setspace,times,contour,cite,geometry,extarrows,bbm,indentfirst,mathtools,lineno,soul}
\usepackage[colorlinks,linkcolor=blue,anchorcolor=blue,citecolor=blue]{hyperref}
\usepackage[normalem]{ulem}
\usepackage[ruled]{algorithm2e}
\graphicspath{{figure/}}
\numberwithin{equation}{section}
\numberwithin{table}{section}
\numberwithin{figure}{section}
\numberwithin{footnote}{section}

\newtheorem{theorem}{Theorem}[section]

\newtheorem{remark}{Remark}[section]

\newtheorem{example}{Example}[section]

\newtheorem{alg}{Algorithm}[section]
\begin{document}
\title{Thermodynamically Consistent Algorithms for Models of Diblock Copolymer Solutions Interacting with Electric and Magnetic Fields}
\author{
Xiaowen Shen\footnote{sxw@csrc.ac.cn, Beijing Computational Science Research Center, Beijing 100193, China.}
,
Qi Wang\footnote{qwang@math.sc.edu, Department of Mathematics, University of South Carolina, Columbia, SC 29208, USA.}
}
\date{\today}
\maketitle

\begin{abstract}
We derive thermodynamically consistent models for diblock copolymer solutions coupled with the electric and magnetic field, respectively. These models satisfy the second law of thermodynamics and therefore are therefore thermodynamically consistent. We then design a set of 2nd order, linear, semi-discrete schemes for the models using the energy quadratization method and the supplementary variable method, respectively, which preserve energy dissipation rates of the models. The spatial discretization is carried out subsequently using 2nd order finite difference methods, leading to fully discrete algorithms that preserve discrete energy-dissipation-rates of the models so that the resulting fully discrete models are thermodynamically consistent. Convergence rates are numerically confirmed through mesh refinement tests and several numerical examples are given to demonstrate the role of the mobility in pattern formation, defect removing effect of both electric and magnetic fields as well as the hysteresis effect for applied external fields in copolymer solutions.
\end{abstract}


\section{Introduction}
A polymer is a class of natural or synthetic substances composed of macromolecules. The block copolymer comprises of two or more homopolymer subunits linked by covalent bonds. The block copolymer solution is a mixture of copolymers and solvent. In the diblock copolymer solution, whose copolymer consists of polymer $A$ and polymer $B$, the solvent can be water or another homopolymer that serves as a solvent for $A/B$ polymers. Diblock copolymers can self-assemble into various ordered mesoscopic structures as the result of polymer-polymer and polymer-solvent interaction, which are useful in materials science and engineering, in particular, nanotechnologies \cite{Hamley_Copolymer_physics, Hamley_Copolymer_science, Hamley_Copolymer_solution}. In experiments, copolymer solutions often undergo macroscopic phase separation between the copolymer and solvent or microscopic phase separation between polymer $A$ and polymer $B$ components. A number of meso-structures of copolymer melts and solutions have been experimentally observed and theoretically predicted via the self-consistent mean field theory as well as the phase field theory. In 2D, there are spot and lamellar morphologies; whereas in 3D cases, the class of possible stable structures is rather abundant including lamellae, spheroids, cylinders and gyrations \cite{Choksi_2009_3D_copolymer_melt, Yamada_2014_Copolymer_solvent, LiBaohui_2015_Double_hydrophilic_copolymer, Makatsoris_2019_3D_diblock_copolymer, ChenChuanjun_2020_Copolymer_melt, YangXiaofeng_2020_Hydrodynamic_copolymer_melt}.

The phase field approach and associated numerical simulations have proven to be effective means to investigate dynamics in copolymer systems. The seminal work of Ohta and Kawasaki \cite{Ohta_1986_Copolymer_melts} laid the foundation of the phase field models for diblock copolymer system. Ohta \textit{et al.} \cite{Ohta_1998_Copolymer_homopolymer} also considered models for mixtures of diblock copolymers and a third immiscible phase. Choksi and Ren \cite{Choksi_2003_Derive_copolymer_melt,Choksi_2005_Derive_copolymer_homopolymer} rederived the models with rigorous mathematical details. Recently, some rigorous analyses on minimizers of the free energy functional or stable steady states of the  Ohta–Kawasaki model in both 2D and 3D were reported \cite{Williams_2017_2D_Ohta-Kawasaki, Williams_2019_3D_Ohta-Kawasaki}, providing insightful information about the dynamics and their attractors in the model of copolymer systems. 

Applying an external electric field or magnetic field to copolymer systems is one of the efficient ways to manipulate and produce various nanostructured materials in self-assembly involving copolymers \cite{Seymour_2015_Phase_field_magnet_electric, Boker_2017_Electric_copolymer, Faghihi_2018_Magnetic_islands}. The electric or magnetic field assisted self-assembly is a cutting-edge research topic in materials science and engineering today.

When the electric field is applied to a copolymer system, dipoles are induced so that the electric field is fully coupled to dynamics of the copolymer system. The existing phase field models for this problem, however, have adopted a over-simplified approach which results in a model that does not respect the thermodynamical consistency in the coupled system of equations \cite{Russell_2005_Electric_copolymer, WuXiangfa_2008_Copolymer_with_electric, Zvelindovsky_2009_Electric_copolymer_lamellae, Glasner_2016_Electric_copolymer_melt, ShenJie_2017_EQ_copolymer, YangXiaofeng_2021_Triblock_copolymer, ChenChuanjun_2020_AC_copolymer}. In fact, all existing models coupling the electric field are based on two original papers \cite{Krausch_2003_Electric_copolymer_solution, Russel_2004_Electric_copolymer}, where significant simplified approximations were applied. As the result, the models do not satisfy the second law of thermodynamics. This thermodynamical inconsistency prevents one from devising structure-preserving numerical approximations, which are thermodynamically consistent at the discrete level, to the important class of models for copolymer solutions with a wide range of applications in fabrications of novel materials. Thermodynamical consistency here refers to the positive entropy production principle in nonequilibrium thermodynamics expressed either as the generalized Onsager principle or as the second law of thermodynamics. It is an important principle to follow when devising mathematical models for any material systems. It is also an important criterion to respect when designing numerical algorithms.

In this paper, we will first derive a thermodynamically consistent phase field model for the diblock copolymer solution or melt system interacting with electric and magnetic fields, following fundamental principles of electrodynamics \cite{Landau_Electrodynamics} and the generalized Onsager principle \cite{WangQi_2016_GOP, WangQi_2021_GOP_application}. In particular, the model for copolymer systems with coupled electric field respects the energy dissipation law, overcoming the shortcomings in the previous models for copolymer systems coupled with the electric field due to over-simplification. These thermodynamically consistent models allow us to design efficient and accurate structure-preserving numerical algorithms to simulate field assisted processing of diblock copolymer solutions based on the dynamical models while respecting the model specific energy dissipation rate.

Around 2011, Badia and Tierra \textit{et al.} explored an idea of transforming the free energy into a quadratic functional to derive energy stable schemes \cite{Badia_2011_EQ, Gonzalez_2013_EQ, Gonzalez_2014_EQ_CH}. Later the idea was amplified and systematically applied to many specific thermodynamic and hydrodynamic models by Yang, Zhao, Shen and Wang \textit{et al.} \cite{WangQi_2018_EQ, ShenJie_2018_SAV, ShenJie_2019_SAV, ShenJie_2020_Lagrange_SAV}, which is termed energy quadratization. The energy quadratization (EQ) strategy bypasses the traditional complicated ones to arrive at semi-discrete, 2nd order or higher order in time linear schemes readily. This strategy is so general that it has little restriction on the specific expression of the free energy density so long as it is bounded below required by physics. Shen, Xu and Yang proposed the scalar auxiliary variable (SAV)  method based upon the EQ strategy for gradient flow systems so long as their free energies are bounded below \cite{ShenJie_2018_SAV, ShenJie_2019_SAV}. In \cite{WangQi_2018_EQ}, Wang \textit{et al.} summarized previous works they had done using the EQ strategy and presented a general framework to discretize the thermodynamically consistent models in forms of partial differential equations to arrive at linear, 2nd order, energy stable numerical schemes. Yang and Shen \textit{et al.} applied the EQ and SAV strategies to copolymer melts and copolymer-homopolymer mixtures and proposed a set of energy stable schemes \cite{ShenJie_2017_EQ_copolymer, YangXiaofeng_2019_CH_log_energy, ChenChuanjun_2020_Copolymer_melt}. Zhang \textit{et al.} used a stabilized SAV method to discretize a magnetic field coupled copolymer melt model and a hydrodynamic model to develop energy stable schemes recently \cite{YangXiaofeng_2020_Magnetic_copolymer, YangXiaofeng_2020_Hydrodynamic_copolymer_melt}.

For a dissipative system, a numerical scheme is called energy stable if it respects energy decay. Among energy stable schemes, the ones that preserve the energy-dissipation-rate are better than the ones that simply preserve energy dissipation inequalities since the energy dissipation rate in the energy stable discrete scheme may not be calculated accurately or even correctly despite of the energy decay property in general. Usually, excessive dissipation is added to the discrete equation system to ensure energy stability. Therefore, developing numerical schemes that preserve the energy-dissipation-rate is desirable and very important. There are no numerical works that preserve energy-dissipation-rates of diblock copolymer systems on the original free energy and the corresponding energy-dissipation-rate yet. 

In \cite{GongYuezheng_2020_SVM, HongQi_2020_SVM}, Wang \textit{et al.} proposed a new temporal discretization paradigm for developing energy-production-rate preserving numerical approximations to thermodynamically consistent partial differential equations with deduced equations, called the supplementary variable method (SVM). The central idea behind this method is to introduce a set of supplementary variables to the thermodynamically consistent model to remove the over-determinedness of the system of equations consisting of the partial differential equations and the deduced equations (the energy definition and the energy-dissipation-rate equation). The introduction of supplementary variables allows one to enforce the deduced equations (i.e., the energy dissipation rate equation in the diblock copolymer model) at the discrete level as intended. This numerical paradigm can be applied to any partial differential equation systems with deduced equations including thermodynamically consistent gradient flow models as well as  hydrodynamic models. Upon discretization, some implementations of this method resemble that of the Lagrange multiplier method proposed by Cheng, Liu and Shen recently \cite{ShenJie_2020_Lagrange_SAV}, and they also share some similarity to the extended SAV approach discussed by Yang and Dong in \cite{DongSuchuan_2020_gPAV} and the projection method in \cite{Hairer_GNI}. This numerical strategy allows a great degree of freedom in the implementation of the supplementary variables in the governing system of equations compared to the other methods mentioned above.

As the second component of this paper, we design a series of energy dissipation rate preserving numerical algorithms for the thermodynamically consistent, electric and magnetic fields coupled phase field model for diblock copolymer solutions, respectively. These schemes are linear, 2nd order and energy-dissipation-rate preserving based on the energy quadratization and the supplementary variable strategy, respectively. The SVM is presented as a prediction-correction scheme, in which the solution and the energy dissipation rate are predicted using 2nd order schemes in the predictive step and the solution and the energy dissipation rate are corrected within the same order to satisfy the discrete equation that respects the energy dissipation rate. The solution of the supplementary variable consumes a negligibly small computational cost, making this method quite efficient compared to the EQ method. We then use a developed scheme to numerically simulate transient dynamics of the diblock copolymer solution coupled with electric and magnetic fields.

The rest of the paper is organized as follows. In Section 2, using the Onsager principle along with the theory on electrodynamics, we derive thermodynamically consistent phase field models for the diblock copolymer solution coupled with an electric field and a magnetic field, respectively. In Section 3, we recall numerical strategies in the EQ method and SVM for gradient flow models. In Section 4, we develop a set of fully discrete, linear, 2nd order, structure-preserving schemes. In Section 5, we numerically demonstrate the accuracy and energy stability of proposed schemes. Then we use one of the SVM schemes to numerically simulate transient dynamics in the copolymer solution system driven by electric and magnetic fields. In Section 6, we give some concluding remarks.

\section{Thermodynamically consistent models}
We derive the thermodynamically consistent models for the diblock copolymer solution that couples the electric and magnetic fields properly. We begin the derivation with a brief review of the generalized Onsager principle which is a constructive method for developing thermodynamically consistent models.

\subsection{Onsager principle}
Nonequilibrium phenomena require models derived from principles of nonequilibrium thermodynamics to describe. The generalized Onsager principle is a linear response theory for developing dynamical models for nonequilibrium phenomena near the equilibrium of any matter systems \cite{Onsager_1931_RRIP1, Onsager_1931_RRIP2, OnsagerMachlup_1953_Irreversible, WangQi_2016_GOP, WangQi_2021_GOP_application}. Assuming the thermodynamical variables describing nonequilibrium phenomena of the matter system are denoted as $\bm\phi(\bm x,t)\in\mathbb R^N$, where $N$ is the number of the variables, and the free energy of the matter system is prescribed by a functional $F[\bm\phi]$ in isothermal cases (for nonisothermal cases, we would have to use an entropy functional instead). The Onsager linear response theory states
\begin{equation}
\bm R \dot{\bm\phi} = - \bm M\bm\mu ,
\quad
\bm\mu = \frac{\delta F}{\delta\bm\phi} ,
\end{equation}
where $\bm R$ and $\bm M$ are operators, $\bm M^{-1}\bm R$ is known as the friction and $\bm R^{-1}\bm M$ as the mobility operator (or matrix) when they exist. This equation system provides relaxation dynamics for the nonequilibrium state to return to equilibrium in dissipative systems or to oscillate in nondissipative systems. Without loss of generality, we assume $\bm R = I$ in this work so that $\bm M$ is the mobility operator.

In general, the mobility operator $\bm M(\bm\phi)$ can be decomposed into a symmetric part $\bm M_s(\bm\phi)$, satisfying the reciprocal relation \cite{Onsager_1931_RRIP1, Onsager_1931_RRIP2}, and an antisymmetric part $\bm M_a(\bm\phi)$ as follows
\begin{equation}
\bm M = \bm M_a + \bm M_s .
\end{equation}
The time rate of change of the free energy is given by
\begin{equation}
\frac{dF}{dt}
= - \int_\Omega \frac{\delta F}{\delta\bm\phi} \bm M \frac{\delta F}{\delta\bm\phi} d\bm x
= - \int_\Omega \frac{\delta F}{\delta\bm\phi} \bm M_s \frac{\delta F}{\delta\bm\phi} d\bm x ,
\end{equation}
under proper adiabatic boundary conditions on $\bm\phi$. Here, $\Omega$ is the domain that the system occupies. If $\bm M_s \geq 0$ is positive semi-definite, the system is a dissipative system; if $\bm M_s=0$, it is a conservative one.

\subsection{Free energy}
We adopt a free energy functional firstly derived by Ohta \textit{et al.} \cite{Ohta_1998_Copolymer_homopolymer} and then reformulated by Choksi and Ren \cite{Choksi_2005_Derive_copolymer_homopolymer}. Consider the diblock copolymer solution in a fixed domain $\Omega\subset\mathbb R^d(d=2,3)$. We utilize ${N_i}$ to denote the polymerization degree of species $i=A,B,S$ ($N_S=1$ when the solvent is water). We denote the local volume fraction of each component as $\phi_i(\bm x,t), i=A,B,S,$ respectively. In the following, we use either $i=A,B,S$ or $i=1,2,3$, or even omit them in some cases for convenience. The free energy functional consists of 2 parts
\begin{equation}
F = F_S + F_L,
\end{equation}
where $F_S, F_L$ represent the short and long range energy, respectively. The short range part $F_S$ is given by
\begin{equation}
F_S = \int_\Omega \Big[ \sum_i \frac{\gamma_i}{2}|\nabla\phi_i|^2
+ \sum_{i,j}\frac{\chi_{ij}}{2}\phi_i\phi_j + f(\bm\phi) \Big] d\bm x ,
\quad
f(\bm\phi) = \sum_i \frac{\phi_i}{N_i}\ln\phi_i ,
\end{equation}
where $\chi_{ij} ( = \chi_{ji} )$ is the interaction constant between polymer $i$ and polymer $j$ ($\chi_{ii}=0$). The long range part $F_L$ is given by
\begin{equation}
F_L = \sum_{i,j} \frac{\alpha_{ij}}{2} \int_\Omega\int_\Omega G(\bm x,\bm y)
[ \phi_i(\bm x)-\overline\phi_i ] [ \phi_j(\bm y)-\overline\phi_j ] d\bm y d\bm x ,
\end{equation}
where
\begin{equation}
\overline\phi_i = \frac{1}{|\Omega|}\int_\Omega \phi_i(\bm x) ~d\bm x , \quad i = A,B,S ,
\end{equation}
$G$ is the Green's function with homogeneous Neumann boundary conditions such that
\begin{equation}
- \Delta_{\bm y}G(\bm x,\bm y)=\delta(\bm x-\bm y) \quad{\rm in}~\Omega ,
\end{equation}
$\delta(\bm x-\bm y)$ denotes the Dirac measure giving unit mass to the point $\bm x$, and $\alpha_{ij}$ are the coupling constants, $\gamma_i$ measure the strength of the conformational entropy. The Kuhn length of the copolymer is set to unity. In the reformulation of Choksi and Ren \cite{Choksi_2005_Derive_copolymer_homopolymer}, parameters of the model admit the following relations
\begin{equation}
\gamma_i = \frac{\varepsilon^2}{\overline\phi_i} ,
\quad
(\alpha_{ij}) = \frac{3}{2}\varepsilon\gamma
\left[ \begin{matrix}
1/(\overline\phi_A\overline\phi_A) & - 1/(\overline\phi_A\overline\phi_B) & 0
\\
- 1/(\overline\phi_A\overline\phi_B) & 1/(\overline\phi_B\overline\phi_B) & 0
\\
0 & 0 & 0
\end{matrix} \right] .
\end{equation}
We define
\begin{equation}
\psi_j(\bm x) := - \int_\Omega G(\bm x,\bm y)[\phi_j(\bm y)-\overline\phi_j] d\bm y
\quad \Leftrightarrow \quad
\Delta\psi_j = \phi_j - \overline\phi_j ,
\end{equation}
with Neumann boundary conditions $\nabla\psi_j\cdot\bm n|_{\partial\Omega}=0, ~j=A,B,S$. Then, we obtain an equivalent free energy
\begin{equation}
F[\bm\phi] = \int_\Omega \Big[
\sum_i \frac{\gamma_i}{2}|\nabla\phi_i|^2
+ \sum_{i,j}\frac{\chi_{ij}}{2}\phi_i\phi_j + f(\bm\phi) \Big] d\bm x
- \sum_{i,j} \frac{\alpha_{ij}}{2} \int_\Omega
[ \phi_i(\bm x)-\overline\phi_i ] \psi_j(\bm x) d\bm x .
\end{equation}

\subsection{Thermodynamically consistent diblock copolymer solution model}
The time rate of change of the free energy given above is calculated as follows
\begin{equation}
\frac{dF}{dt} = \int_\Omega
\sum_i\Big(-\gamma_i\Delta\phi_i+\sum_{j}\chi_{ij}\phi_j+\frac{\partial f}{\partial\phi_i}\Big)
\frac{\partial\phi_i}{\partial t} d\bm x
+ \int_{\partial\Omega} \sum_i\gamma_i\nabla\phi_i\cdot\bm n\frac{\partial\phi_i}{\partial t} ds
- \sum_{i,j}\alpha_{ij}\int_\Omega\psi_j\frac{\partial\phi_i}{\partial t} d\bm x .
\end{equation}
In the adiabatic case, the surface energy dissipation vanishes at the zero Neumann boundary conditions:
\begin{equation}\label{boundary condition of phi}
\nabla\phi_i \cdot \bm n |_{\partial\Omega} = 0 , \quad i = A, B, S .
\end{equation}
The chemical potentials are given by
\begin{equation}
\begin{aligned}
\mu_A &= \frac{\delta F}{\delta\phi_A}
= - \gamma_A\Delta\phi_A + \sum_{j}\chi_{Aj}\phi_j + \frac{\partial f}{\partial\phi_A}
- \alpha_{AA}\psi_A - \alpha_{AB}\psi_B ,
\\
\mu_B &= \frac{\delta F}{\delta\phi_B}
= - \gamma_B\Delta\phi_B + \sum_{j}\chi_{Bj}\phi_j + \frac{\partial f}{\partial\phi_B}
- \alpha_{AB}\psi_A - \alpha_{BB}\psi_B ,
\\
\mu_S &= \frac{\delta F}{\delta\phi_S}
= - \gamma_S\Delta\phi_S + \sum_{j}\chi_{Sj}\phi_j + \frac{\partial f}{\partial\phi_S} .
\end{aligned}
\end{equation}
Now we derive the transport equations for the copolymer solution system using the Onsager principle. The transport equations are proposed as
\begin{equation}
\frac{\partial\phi_i}{\partial t} = j_i , \quad i=A,B,S,
\end{equation}
where $j_i$ is the  rate of production for $\phi_i,i=A,B,S$, respectively. We add the incompressibility condition to the free energy as a constraint through Lagrange multiplier $L$:
\begin{equation}
E[\bm\phi,L] = F[\bm\phi] + \int_\Omega L(\bm x,t)(\phi_A+\phi_B+\phi_S-1) d\bm x .
\end{equation}
The time derivative of the free energy becomes
\begin{equation}
\frac{dE}{dt} = \int_\Omega \Big[
(\mu_A+L)\frac{\partial\phi_A}{\partial t}
+ (\mu_B+L)\frac{\partial\phi_B}{\partial t}
+ (\mu_S+L)\frac{\partial\phi_S}{\partial t} \Big] d\bm x .
\end{equation}
We apply the Onsager principle to obtain
\begin{equation}
[j_A, j_B, j_S]^\top = - \bm M [\mu_A+L, \mu_B+L, \mu_S+L]^\top ,
\end{equation}
where
\begin{equation}
\bm M = - \nabla\cdot(M\nabla) =
- \nabla\cdot \Big( \left[ \begin{matrix}
M_{11}(\bm\phi) & M_{12}(\bm\phi) & M_{13}(\bm\phi)
\\
M_{12}(\bm\phi) & M_{22}(\bm\phi) & M_{23}(\bm\phi)
\\
M_{13}(\bm\phi) & M_{23}(\bm\phi) & M_{33}(\bm\phi)
\end{matrix} \right] \nabla \Big)
\end{equation}
is the mobility operator and $M_{ij}$ are the coefficients. So the dynamical equations are given as follows
\begin{equation}\label{CH}
\left\{
\begin{aligned}
&\Delta\psi_A=\phi_A-\overline\phi_A, \quad \Delta\psi_B=\phi_B-\overline\phi_B,
\\
&\frac{\partial\phi_A}{\partial t} = \nabla\cdot\big[M_{11}\nabla(\mu_A+L)\big]
+ \nabla\cdot\big[M_{12}\nabla(\mu_B+L)\big] + \nabla\cdot\big[M_{13}\nabla(\mu_S+L)\big] ,
\\
&\frac{\partial\phi_B}{\partial t} =\nabla\cdot\big[M_{12}\nabla(\mu_A+L)\big]
+ \nabla\cdot\big[M_{22}\nabla(\mu_B+L)\big] + \nabla\cdot\big[M_{23}\nabla(\mu_S+L)\big] ,
\\
&\frac{\partial\phi_S}{\partial t} =\nabla\cdot\big[M_{13}\nabla(\mu_A+L)\big]
+ \nabla\cdot\big[M_{23}\nabla(\mu_B+L)\big] + \nabla\cdot\big[M_{33}\nabla(\mu_S+L)\big] ,
\\
&\nabla\cdot\Big(\sum_{i,j}M_{ij}\nabla L\Big) = -\nabla\cdot\Big(\sum_{i,j}M_{ij}\nabla\mu_i\Big).
\end{aligned}
\right.
\end{equation}
The adiabatic boundary conditions for the copolymer model are
\begin{equation}
\bm n\cdot M\cdot\nabla[\mu_A+L,\mu_B+L,\mu_S+L]^\top\big|_{\partial\Omega}=\bm 0,
~
\bm n\cdot\nabla\phi_i|_{\partial\Omega}=0, i=A,B,S,
~
\nabla\psi_j\cdot\bm n|_{\partial\Omega}=0, j=A,B.
\end{equation}
Under these boundary conditions, we obtain the energy dissipation rate
\begin{equation}
\frac{dE}{dt} = -\int_\Omega
\nabla[\mu_A+L,\mu_B+L,\mu_S+L]\cdot M\cdot\nabla[\mu_A+L,\mu_B+L,\mu_S+L] d\bm x .
\end{equation}
The model is energy dissipative provided $(M_{ij})\geq 0$. We note that the Onsager principle also allows one to propose nonadiabatic boundary conditions in which energy can exchange with the surrounding. For simplicity, we focus on the adiabatic case in this study.

\subsection{Electric field coupled model}
Copolymers are responsive to external field such as electric and magnetic fields. Field assisted processing is one of the promising means to control the morphology of the copolymer system. The electric and magnetic fields impact on the copolymer system by modifying the free energy landscape. Developing a thermodynamically consistent relaxation model for the copolymer system subject to external fields is an important first step in understanding dynamics of the system. Finding stable steady state (or low energy state) along with the energy decaying path is of great importance \cite{Choksi_2011_2D_copolymer_melt, Choksi_2015_Metastable_states, ZhangLei_2020_Energy_landscape, ZhangLei_2020_Find_saddle_points}. In this paper, we utilize the PDE paradigm to describe the dynamical path, whose governing equations are derived from the Onsager principle. This is one legitimate description of relaxation dynamics from the nonequilibrium state to return to the equilibrium state through the mobility operator.

We consider the additional free energy due to the existence of electric field $\bm E$
\begin{equation}
F_{\rm e} = \int_\Omega\frac{1}{2}\bm D\cdot\bm E ~ d\bm x
= \int_\Omega\frac{1}{2\epsilon}\bm D\cdot\bm D ~ d\bm x ,
\end{equation}
where $\bm D=\epsilon\bm E$ is the electric displacement, $\epsilon$ is the dielectric constant, which can be approximated as \cite{Krausch_2003_Electric_copolymer_solution, Russel_2004_Electric_copolymer}
\begin{equation}
\epsilon(\bm x, t) = \epsilon_0 + \epsilon_1\varphi ,
\quad
\varphi = (\phi_A-\phi_B) - (\overline\phi_A-\overline\phi_B) ,
\end{equation}
where $\epsilon_0$ is the vacuum permittivity, $\epsilon_1$ is a constant related to the material. Because of the inhomogeneity of $\epsilon(\bm x,t)$, the electric field $\bm E$ inside the material deviates from the applied electric field $\bm E_0$ and is given by
\begin{equation}
\bm E = \bm E_0 - \nabla\Phi ,
\end{equation}
where $\Phi$ is the electric potential for the induced electric field. $\Phi$ can be expressed in terms of $\varphi$ or $\bm\phi$, by solving the following Poisson equation (Gauss theorem)
\begin{equation}\label{Maxwell_equation}
\nabla\cdot\bm D = \nabla\cdot\big[(\epsilon_0+\epsilon_1\varphi)(\bm E_0-\nabla\Phi)\big]
= - \epsilon_0\Delta\Phi + \epsilon_1\bm E_0\cdot\nabla\varphi
- \epsilon_1\nabla\varphi\cdot\nabla\Phi - \epsilon_1\varphi\Delta\Phi
= 0 ,
\end{equation}
where we use the Dirichlet boundary condition $\Phi|_{\partial\Omega}=\Phi_0$ for the electric potential and $\Phi_0$ is a prescribed value (we use $\Phi_0=0$ in this study).

Following the work of \cite{Landau_Electrodynamics}, we have
\begin{equation}
\frac{dF_{\rm e}}{dt} = \int_\Omega - \frac{1}{2\epsilon^2}|\bm D|^2
\frac{\partial\epsilon}{\partial\varphi}\frac{\partial\varphi}{\partial t} d\bm x
= \int_\Omega - \frac{\epsilon_1}{2}|\bm E|^2
\Big(\frac{\partial\phi_A}{\partial t}-\frac{\partial\phi_B}{\partial t}\Big) d\bm x .
\end{equation}
Since the total free energy is given by
\begin{equation}
E[\bm\phi,L] = F[\bm\phi] + F_{\rm e}[\bm\phi]
+ \int_\Omega L(\bm x,t)(\phi_A+\phi_B+\phi_S-1) d\bm x .
\end{equation}
the time derivative of the total free energy becomes
\begin{equation}
\frac{dE}{dt} = \int_\Omega \Big[
(\mu_A+\mu_{\rm e}+L)\frac{\partial\phi_A}{\partial t}
+ (\mu_B-\mu_{\rm e}+L)\frac{\partial\phi_B}{\partial t}
+ (\mu_S+L)\frac{\partial\phi_S}{\partial t} \Big] d\bm x ,
~
\mu_{\rm e} = -\frac{\epsilon_1}{2}|\bm E_0-\nabla\Phi|^2 .
\end{equation}
We denote $\widehat\mu_A=\mu_A+\mu_{\rm e}, ~\widehat\mu_B=\mu_B-\mu_{\rm e}, ~\widehat\mu_S=\mu_S$. Using the Onsager principle, we obtain the model coupled with the electric field
\begin{equation}
\left\{
\begin{aligned}
&\Delta\psi_A=\phi_A-\overline\phi_A, \quad \Delta\psi_B=\phi_B-\overline\phi_B,
\\
&\nabla\cdot\big[(\epsilon_0+\epsilon_1\varphi)(\bm E_0-\nabla\Phi)\big] = 0 ,
\\
& \frac{\partial\phi_A}{\partial t} =
\nabla\cdot\Big[M_{11}\nabla(\widehat\mu_A+L)\Big]
+ \nabla\cdot\Big[M_{12}\nabla(\widehat\mu_B+L)\Big]
+ \nabla\cdot\big[M_{13}\nabla(\widehat\mu_S+L)\big] ,
\\
& \frac{\partial\phi_B}{\partial t} =
\nabla\cdot\Big[M_{12}\nabla(\widehat\mu_A+L)\Big]
+ \nabla\cdot\Big[M_{22}\nabla(\widehat\mu_B+L)\Big]
+ \nabla\cdot\big[M_{23}\nabla(\widehat\mu_S+L)\big] ,
\\
& \frac{\partial\phi_S}{\partial t} =
\nabla\cdot\Big[M_{13}\nabla(\widehat\mu_A+L)\Big]
+ \nabla\cdot\Big[M_{23}\nabla(\widehat\mu_B+L)\Big]
+ \nabla\cdot\big[M_{33}\nabla(\widehat\mu_S+L)\big] ,
\\
& \nabla\cdot\Big(\sum_{i,j}M_{ij}\nabla L\Big) =
- \nabla\cdot\Big(\sum_{i,j}M_{ij}\nabla\widehat\mu_i\Big) .
\end{aligned}
\right.
\end{equation}
The adiabatic boundary conditions for the electric field coupled copolymer model are
\begin{equation}
\bm n\cdot M\cdot\nabla\big[\widehat\mu_A+L,\widehat\mu_B+L,\widehat\mu_S+L\big]^\top
\big|_{\partial\Omega} \!=\!\bm 0,
~
\bm n\cdot\nabla\phi_i\big|_{\partial\Omega} \!=\! 0, i\!=\!A,\!B,\!S,
~
\nabla\psi_j\cdot\bm n|_{\partial\Omega} \!=\! 0, j\!=\!A,\!B,
~
\Phi\big|_{\partial\Omega} \!=\! \Phi_0 .
\end{equation}
The energy dissipation rate is given by
\begin{equation}
\frac{dE}{dt} = - \int_\Omega \nabla
\big[ \widehat\mu_A+L, \widehat\mu_B+L, \widehat\mu_S+L \big]
\cdot M\cdot \nabla
\big[ \widehat\mu_A+L, \widehat\mu_B+L, \widehat\mu_S+L \big] d\bm x .
\end{equation}
The model is energy dissipative provided $(M_{ij})\geq 0$.

\begin{remark}
Consider a diblock copolymer melt system. If we drop nonlinear terms in \eqref{Maxwell_equation} and approximate it by
\begin{equation}
\Delta\Phi = \frac{\epsilon_1}{\epsilon_0}\bm E_0\cdot\nabla\varphi ,
\end{equation}
we obtain a simplified electric field coupled model for diblock copolymer melts
\begin{equation}
\frac{\partial\varphi}{\partial t} =
M\Delta\Big(\frac{\delta F}{\delta\varphi}-\frac{\epsilon_1}{2}|\bm E_0-\nabla\Phi|^2\Big)
= M\Delta\frac{\delta F}{\delta\phi}
+ M\frac{\epsilon_1^2}{\epsilon_0}(\bm E_0-\nabla\Phi)\cdot\nabla\nabla\varphi\cdot\bm E_0
- M\epsilon_1\nabla\nabla\Phi:\nabla\nabla\Phi^\top .
\end{equation}
Furthermore, we drop the last term and use approximation $\bm E_0-\nabla\Phi\approx\bm E_0$ in the above, we recover the approximate model used in \cite{Krausch_2003_Electric_copolymer_solution, Russel_2004_Electric_copolymer, Russell_2005_Electric_copolymer, WuXiangfa_2008_Copolymer_with_electric, Zvelindovsky_2009_Electric_copolymer_lamellae, Glasner_2016_Electric_copolymer_melt, ShenJie_2017_EQ_copolymer, YangXiaofeng_2021_Triblock_copolymer, ChenChuanjun_2020_AC_copolymer}
\begin{equation}
\frac{\partial\varphi}{\partial t}
= M\Delta\frac{\delta F}{\delta\phi}
+ M\frac{\epsilon_1^2}{\epsilon_0}\bm E_0\cdot\nabla\nabla\varphi\cdot\bm E_0
~\Rightarrow~
\frac{\partial\varphi}{\partial t}
= M\Delta\frac{\delta F}{\delta\phi}
+ M\frac{\epsilon_1^2}{\epsilon_0}E_1^2\partial^2_x\varphi ,
~\text{when}~\bm E_0=[E_1,0]^\top\in\mathbb R^2 .
\end{equation}
This model is not thermodynamically consistent!
\end{remark}

\subsection{Magnetic field coupled model}
When the magnetic field is present, we consider the corresponding free energy due to the imposed magnetic field $\bm B_0$. The free energy corresponding to the polymer-magnetic field interaction must be added with an additional term \cite{Landau_Electrodynamics}
\begin{equation}
F_{\rm m} = \int_\Omega\frac{\gamma_m}{2}|\bm B_0\cdot\nabla\varphi|^2 d\bm x ,
\quad
\varphi = (\phi_A-\phi_B) - (\overline\phi_A-\overline\phi_B) ,
\end{equation}
where $\gamma_m$ is the constant parameter related to magnetic permeability. Taking the time derivative of this free energy term, we obtain
\begin{equation}
\frac{dF_{\rm m}}{dt} = \int_\Omega
- \gamma_m\bm B_0\cdot\nabla\nabla(\phi_A-\phi_B)\cdot\bm B_0
\Big( \frac{\partial\phi_A}{\partial t}-\frac{\partial\phi_B}{\partial t} \Big) d\bm x
+ \int_{\partial\Omega}\gamma_m\bm B_0\cdot\nabla\varphi\bm B_0\cdot\bm n
\Big( \frac{\partial\phi_A}{\partial t}-\frac{\partial\phi_B}{\partial t} \Big) ds ,
\end{equation}
where the boundary integral vanishes due to zero Neumann boundary condition \eqref{boundary condition of phi}. Since the total free energy is given by
\begin{equation}
E[\bm\phi,L] = F[\bm\phi] + F_{\rm m}[\bm\phi]
+ \int_\Omega L(\bm x,t)(\phi_A+\phi_B+\phi_S-1) d\bm x ,
\end{equation}
the time derivative of the total free energy becomes
\begin{equation}
\frac{dE}{dt} = \int_\Omega \Big[
(\mu_A+\mu_{\rm m}+L)\frac{\partial\phi_A}{\partial t}
+ (\mu_B-\mu_{\rm m}+L)\frac{\partial\phi_B}{\partial t}
+ (\mu_S+L)\frac{\partial\phi_S}{\partial t} \Big] d\bm x ,
~
\mu_{\rm m} = - \gamma_m\bm B_0\cdot\nabla\nabla(\phi_A-\phi_B)\cdot\bm B_0 .
\end{equation}
We denote $\widehat\mu_A=\mu_A+\mu_{\rm m}, ~\widehat\mu_B=\mu_B-\mu_{\rm m}, ~\widehat\mu_S=\mu_S$. Using the Onsager principle, we obtain the model coupled with the magnetic field
\begin{equation}
\left\{
\begin{aligned}
&\Delta\psi_A=\phi_A-\overline\phi_A, \quad \Delta\psi_B=\phi_B-\overline\phi_B,
\\
& \frac{\partial\phi_A}{\partial t} =
\nabla\cdot\Big[M_{11}\nabla(\widehat\mu_A+L)\Big]
+ \nabla\cdot\Big[M_{12}\nabla(\widehat\mu_B+L)\Big]
+ \nabla\cdot\big[M_{13}\nabla(\widehat\mu_S+L)\big] ,
\\
& \frac{\partial\phi_B}{\partial t} =
\nabla\cdot\Big[M_{12}\nabla(\widehat\mu_A+L)\Big]
+ \nabla\cdot\Big[M_{22}\nabla(\widehat\mu_B+L)\Big]
+ \nabla\cdot\big[M_{23}\nabla(\widehat\mu_S+L)\big] ,
\\
& \frac{\partial\phi_S}{\partial t} =
\nabla\cdot\Big[M_{13}\nabla(\widehat\mu_A+L)\Big]
+ \nabla\cdot\Big[M_{23}\nabla(\widehat\mu_B+L)\Big]
+ \nabla\cdot\big[M_{33}\nabla(\widehat\mu_S+L)\big] ,
\\
& \nabla\cdot\Big(\sum_{i,j}M_{ij}\nabla L\Big) =
- \nabla\cdot\Big(\sum_{i,j}M_{ij}\nabla\widehat\mu_i\Big) .
\end{aligned}
\right.
\end{equation}
The adiabatic boundary conditions for the magnetic field coupled copolymer model are
\begin{equation}
\bm n\cdot M\cdot\nabla\big[\widehat\mu_A+L,\widehat\mu_B+L,\widehat\mu_S+L\big]^\top
\big|_{\partial\Omega} = \bm 0,
~
\bm n\cdot\nabla\phi_i\big|_{\partial\Omega} = 0, i=A,B,S,
~
\nabla\psi_j\cdot\bm n|_{\partial\Omega}=0, j=A,B.
\end{equation}
The energy dissipation rate is given by
\begin{equation}
\frac{dE}{dt} = - \int_\Omega \nabla
\big[ \widehat\mu_A+L, \widehat\mu_B+L, \widehat\mu_S+L \big]
\cdot M\cdot \nabla
\big[ \widehat\mu_A+L, \widehat\mu_B+L, \widehat\mu_S+L \big] d\bm x .
\end{equation}
The model is energy dissipative provided $(M_{ij})\geq 0$.

Next, we discuss how to design numerical algorithms that respect the thermodynamical consistency in the discrete models. We call them the thermodynamically consistent algorithms.

\section{Numerical strategies in time}
We discuss a couple of numerical strategies for discretizing the field-coupled copolymer models in time in this section and then discuss the spatial discretization in the next section. We focus on two methods: the energy quadratization (EQ) method and the supplementary variable method (SVM).

\subsection{Energy quadratization (EQ) method}
We recall the EQ procedure applied to a gradient flow system
\begin{equation}
\frac{\partial\bm\phi}{\partial t} = - \bm M \bm\mu ,
\quad
\bm\mu = \frac{\delta F}{\delta\bm\phi} ,
\quad
F[\bm\phi] = \int_\Omega
\Big[ \frac{1}{2}\bm\phi\mathcal L\bm\phi + f(\bm\phi,\nabla\bm\phi) \Big] d\bm x ,
\end{equation}
where $\mathcal L$ is a linear, self-adjoint, positive definite operator. The representation of the EQ process is not unique. For example, we introduce an intermediate variable
\begin{equation}
q = \sqrt{ f(\bm\phi, \nabla\bm\phi) + C } ,
\end{equation}
where constant $C>0$ is large enough to make $q$ real-valued for all $\bm\phi$, to quadratize the free energy density. The free energy then is rewritten into a quadratic form
\begin{equation}
F[\bm\phi, q] = \int_\Omega
\Big( \frac{1}{2}\bm\phi\mathcal L\bm\phi + q^2 - C \Big) d\bm x .
\end{equation}
The reformulated chemical potential is
\begin{equation}
\bm\mu = \mathcal L\bm\phi + 2q\frac{\partial q}{\partial\bm\phi}
- \nabla\cdot\Big( 2q\frac{\partial q}{\partial\nabla\bm\phi} \Big) .
\end{equation}
The time evolutionary equation of the intermediate variable is
\begin{equation}
\frac{\partial q}{\partial t} =
\frac{\partial q}{\partial\bm\phi}\cdot\frac{\partial\bm\phi}{\partial t}
+ \frac{\partial q}{\partial\nabla\bm\phi}\cdot\nabla\frac{\partial\bm\phi}{\partial t}
= \sum_i \Big(
\frac{\partial q}{\partial\phi_i}\frac{\partial\phi_i}{\partial t}
+ \frac{\partial q}{\partial\nabla\phi_i}\cdot\nabla\frac{\partial\phi_i}{\partial t}
\Big) .
\end{equation}
We obtain the equivalent system
\begin{equation}
\left\{
\begin{aligned}
\frac{\partial\bm\phi}{\partial t} &= - \bm M
\Big[ \mathcal L\bm\phi + 2q\frac{\partial q}{\partial\bm\phi}
- \nabla\cdot\Big( 2q\frac{\partial q}{\partial\nabla\bm\phi} \Big) \Big] ,
\\
\frac{\partial q}{\partial t} &=
\frac{\partial q}{\partial\bm\phi}\cdot\frac{\partial\bm\phi}{\partial t}
+ \frac{\partial q}{\partial\nabla\bm\phi}\cdot\nabla\frac{\partial\bm\phi}{\partial t} ,
\end{aligned}
\right.
\end{equation}
with consistent initial condition for $q$ and the boundary conditions for the original variables. The reformulated model shares the same energy dissipation rate as the original one
\begin{equation}
\frac{dF}{dt} = - \int_\Omega
\Big[ \mathcal L\bm\phi + 2q\frac{\partial q}{\partial\bm\phi}
- \nabla\cdot\Big( 2q\frac{\partial q}{\partial\nabla\bm\phi} \Big) \Big]
\bm M
\Big[ \mathcal L\bm\phi + 2q\frac{\partial q}{\partial\bm\phi}
- \nabla\cdot\Big( 2q\frac{\partial q}{\partial\nabla\bm\phi} \Big) \Big]
d\bm x
\end{equation}
with the consistent adiabatic boundary conditions on $\bm\phi$.

\subsection{Supplementary variable method}
EQ methods preserve the energy dissipation rate for an EQ reformulated quadratic free energy. Here, we present the supplementary variable method for the gradient flow model that preserves the energy dissipation rate of the original free energy functional. Consider the following gradient flow system
\begin{equation}\label{SVM gradient flow}
\frac{\partial\bm\phi}{\partial t} = - \bm M\bm\mu ,
\quad
\bm\mu = \frac{\delta F}{\delta\bm\phi} ,
\end{equation}
where $\bm\phi\in\mathbb R^N$ is the thermodynamic variable describing nonequilibrium phenomena of the system, $\bm M(\bm\phi)$ is the mobility operator, the free energy is expressed as
\begin{equation}\label{SVM free energy}
F = \frac{1}{2}(\bm\phi, \mathcal L\bm\phi) + \big( f(\bm\phi), 1 \big) ,
\end{equation}
$(\cdot,\cdot)$ is the $L^2$ inner product, $\mathcal L$ is a linear self-adjoint positive definite operator, $f$ is bounded below for physically accessible state of $\bm\phi$. In SVM paradigm, we view $F$ as one of the thermodynamical variables. Under proper boundary conditions, we obtain an over-determined system with $N+2$ equations and $N+1$ unknowns $[\bm\phi,F]$:
\begin{subequations}\label{SVM equations}
\begin{align}
& \frac{\partial\bm\phi}{\partial t} = - \bm M\big(\mathcal L\bm\phi + f'(\bm\phi)\big) ,
\label{SVM phi_t}
\\
& \frac{dF}{dt} = - (\bm\mu, \bm M\bm\mu) , \label{SVM F_t}
\\
& F = \frac{1}{2}(\bm\phi, \mathcal L\bm\phi) + \big( f(\bm\phi), 1 \big) .
\end{align}
\end{subequations}
This system is consistent and solvable should gradient flow model \eqref{SVM gradient flow} is solvable with free energy \eqref{SVM free energy} and proper initial and boundary conditions.

However, this system is structurally unstable in that if the equations are discretized in \eqref{SVM equations}, there is a high likelihood that the discretized system would be inconsistent so that there would be no solutions of the discretized system. To address this catastrophic problem, we augment \eqref{SVM equations} with sufficient number of supplementary variables properly so that the augmented system is well-determined and structural stable. For this gradient flow model, we need 1 supplementary variable. Here we present an implementation of SVM in a prediction-correction formulation. There are other ways SVM can be implemented which we will not elaborate here.

\subsubsection{Predictive step}
We rewrite \eqref{SVM phi_t} and \eqref{SVM F_t} as
\begin{subequations}
\begin{align}
& \frac{\partial\bm\phi}{\partial t} = - \bm M\big( \mathcal L\bm\phi + f'(\bm\phi) \big) ,
\label{SVM equation phi}
\\
& \bm\mu = \mathcal L\bm\phi + f'(\bm\phi) ,
\\
& \frac{dF}{dt} = - (\bm\mu, \bm M\bm\mu) .
\end{align}
\end{subequations}
Assuming $\bm\phi^{n-1},\bm\phi^n$ are known, we define
\begin{equation}
\bm\mu^n = \mathcal L\bm\phi^n + f'(\bm\phi^n) ,
\quad
F^n = \frac{1}{2}(\bm\phi^n, \mathcal L\bm\phi^n) + \big( f(\bm\phi^n), 1 \big) ,
\quad
\overline{\bm\phi}^{n+\frac{1}{2}} = \frac{3\bm\phi^n-\bm\phi^{n-1}}{2} .
\end{equation}
Firstly, we approximate the energy-dissipation-rate equation of the system as follows
\begin{subequations}
\begin{align}
& \frac{\widetilde{\bm\phi}^{n+\frac{1}{2}} - \bm\phi^n}{\triangle t/2} =
- \bm M(\overline{\bm\phi}^{n+\frac{1}{2}})
\left[\mathcal L\widetilde{\bm\phi}^{n+\frac{1}{2}}+f'(\overline{\bm\phi}^{n+\frac{1}{2}})\right] ,
\\
& \widetilde{\bm\mu}^{n+\frac{1}{2}} = \mathcal L\widetilde{\bm\phi}^{n+\frac{1}{2}}
+ f'(\widetilde{\bm\phi}^{n+\frac{1}{2}}) ,
\label{SVM predict mu}
\\
& \frac{\widetilde F^{n+1}-F^n}{\triangle t} = - \left( \widetilde{\bm\mu}^{n+\frac{1}{2}},
~ \bm M(\widetilde{\bm\phi}^{n+\frac{1}{2}})\widetilde{\bm\mu}^{n+\frac{1}{2}} \right) ,
\label{SVM predict F}
\end{align}
\end{subequations}
where $\triangle t$ is the time step. It follows that
\begin{equation}\label{SVM predict}
\widetilde{\bm\phi}^{n+\frac{1}{2}} \!=\!
\Big[\! 1 \!+\! \frac{\triangle t}{2}\bm M(\overline{\bm\phi}^{n+\frac{1}{2}})\mathcal L \Big]^{-1}\!
\Big[\! \bm\phi^n \!-\! \frac{\triangle t}{2}\bm M(\overline{\bm\phi}^{n+\frac{1}{2}})
f'(\overline{\bm\phi}^{n+\frac{1}{2}}) \Big] ,
\widetilde F^{n+1} \!=\! F^n \!-\! \triangle t\left( \widetilde{\bm\mu}^{n+\frac{1}{2}},
\bm M(\widetilde{\bm\phi}^{n+\frac{1}{2}})\widetilde{\bm\mu}^{n+\frac{1}{2}} \right) ,
\end{equation}
which is a 2nd order approximation to the model.

\subsubsection{Corrective step}
Secondly, we impose \eqref{SVM predict F} as a constraint for equation \eqref{SVM equation phi}, i.e.
\begin{equation}\label{SVM constraint}
F[\bm\phi(t_{n+1})] = \widetilde F^{n+1} .
\end{equation}
To solve equation \eqref{SVM equation phi} with the constraint \eqref{SVM constraint}, we modify \eqref{SVM equation phi} by a time-dependent supplementary variable $\alpha(t)$ together with a user supplied function $\bm g(\bm\phi)$
\begin{equation}\label{SVM system}
\left\{
\begin{aligned}
& \frac{\partial\bm\phi}{\partial t} =
- \bm M\big(\mathcal L\bm\phi+f'(\bm\phi)\big)+\alpha(t)\bm g(\bm\phi) , \quad t_n<t\leq t_{n+1},
\\
& F[\bm\phi(t_{n+1})] = \widetilde F^{n+1} .
\end{aligned}
\right.
\end{equation}
There is a great deal of flexibility in determining how to modify the gradient flow model with a supplementary variable. Here we present 4 implementations. In theory, there are infinitely many choices to choose from.
\begin{itemize}
\item \textbf{SVM1}:
We take
\begin{equation}\label{SVM g1}
\bm g(\bm\phi) = - \bm M f'(\bm\phi) ,
\end{equation}
which leads to
\begin{equation}\label{SVM1 system}
\left\{
\begin{aligned}
& \frac{\partial\bm\phi}{\partial t} = - \bm M\big(\mathcal L\bm\phi+(1+\alpha)f'(\bm\phi)\big) ,
\\
& \bm\mu = \mathcal L\bm\phi + (1+\alpha)f'(\bm\phi) ,
\\
& \frac{dF}{dt} = - (\bm\mu, \bm M\bm\mu) .
\end{aligned}
\right.
\end{equation}
This is a new system with a modified chemical potential.

\item \textbf{SVM2}:
We take
\begin{equation}\label{SVM g2}
\bm g(\bm\phi) = - \bm M\big( \mathcal L\bm\phi + f'(\bm\phi) \big) ,
\end{equation}
which leads to
\begin{equation}
\left\{
\begin{aligned}
& \frac{\partial\bm\phi}{\partial t} = - (1+\alpha)\bm M\big(\mathcal L\bm\phi+f'(\bm\phi)\big) ,
\\
& \bm\mu = \mathcal L\bm\phi + f'(\bm\phi) ,
\\
& \frac{dF}{dt} = - (\bm\mu, \bm M\bm\mu) .
\end{aligned}
\right.
\end{equation}
This is a new system with a modified mobility.
\end{itemize}
Note that \eqref{SVM1 system} is equivalent to the reformulated model proposed in \cite{ShenJie_2020_Lagrange_SAV} using what they call the Lagrange multiplier method. The supplementary variable is a perturbation variable since when $\alpha(t)=0$, the modified/perturbed system reduces to the original one. By supplementing the over-determined system with a new variable without an equation, we effectively remove the over-determinedness by embedding the system into a higher dimensional phase space! We reiterate that there is a great deal of flexibility to place the supplementary variable in the model, which differentiates it from the Lagrange multiplier method and the projection method in structural preserving approximations. We introduce 2 other implementations in following numerical procedures.

\subsubsection{Semi-discrete numerical algorithms in time}
We apply the implicit-explicit Crank-Nicolson scheme in time to \eqref{SVM system} to arrive at 
\begin{equation}\label{SVM semi scheme}
\left\{
\begin{aligned}
& \delta^+_t\bm\phi^n = - \bm M(\widetilde{\bm\phi}^{n+\frac{1}{2}})
\Big[\mathcal L\bm\phi^{n+\frac{1}{2}}+f'(\widetilde{\bm\phi}^{n+\frac{1}{2}})\Big]
+ \alpha^{n+\frac{1}{2}}\bm g(\widetilde{\bm\phi}^{n+\frac{1}{2}}) ,
\\
& F[\bm\phi^{n+1}] = \widetilde F^{n+1} ,
\end{aligned}
\right.
\end{equation}
where
\begin{equation}
\delta^+_t\bm\phi^n = \frac{\bm\phi^{n+1}-\bm\phi^n}{\triangle t} ,
\quad
\bm\phi^{n+\frac{1}{2}} = \frac{\bm\phi^{n+1}+\bm\phi^n}{2} ,
\end{equation}
and $\widetilde{\bm\phi}^{n+\frac{1}{2}}, \widetilde F^{n+1}$ are obtained through \eqref{SVM predict}. Next we discuss how to solve system \eqref{SVM semi scheme} efficiently. Let
\begin{equation}\label{SVM procedure}
\left\{
\begin{aligned}
& \widehat{\bm\phi}^{n+1} =
\Big[1+\frac{\triangle t}{2}\bm M(\widetilde{\bm\phi}^{n+\frac{1}{2}})\mathcal L\Big]^{-1}
\Big[ \Big(1-\frac{\triangle t}{2}\bm M(\widetilde{\bm\phi}^{n+\frac{1}{2}})\mathcal L\Big)
\bm\phi^n - \triangle t\bm M(\widetilde{\bm\phi}^{n+\frac{1}{2}})
f'(\widetilde{\bm\phi}^{n+\frac{1}{2}}) \Big] ,
\\
& \check{\bm\phi}^{n+\frac{1}{2}} = \Big[1+\frac{\triangle t}{2}\bm M(\widetilde{\bm\phi}^{n+\frac{1}{2}})\mathcal L\Big]^{-1}
\bm g(\widetilde{\bm\phi}^{n+\frac{1}{2}}) ,
\\
& \beta = \alpha^{n+\frac{1}{2}} \triangle t .
\end{aligned}
\right.
\end{equation}
It follows from the first equation in \eqref{SVM semi scheme} that
\begin{equation}\label{SVM update phi}
\bm\phi^{n+1} = \widehat{\bm\phi}^{n+1} + \beta \check{\bm\phi}^{n+\frac{1}{2}} .
\end{equation}
We give 2 more implementations hinted by \eqref{SVM update phi}.
\begin{itemize}
\item \textbf{SVM3}:
We take
\begin{equation}\label{SVM g3}
\bm g(\bm\phi) = \Big( 1+\frac{\triangle t}{2}\bm M\mathcal L \Big)\bm\phi ,
\end{equation}
which leads to $\check{\bm\phi}^{n+\frac{1}{2}} = \widetilde{\bm\phi}^{n+\frac{1}{2}}$.

\item \textbf{SVM4}:
We take
\begin{equation}\label{SVM g4}
\check{\bm\phi}^{n+\frac{1}{2}} = \widehat{\bm\phi}^{n+1} .
\end{equation}
\end{itemize}
We plug \eqref{SVM update phi} into the second equation in \eqref{SVM semi scheme} to arrive at
\begin{equation}\label{SVM algebraic equation}
F[\widehat{\bm\phi}^{n+1} + \beta \check{\bm\phi}^{n+\frac{1}{2}}] = \widetilde F^{n+1} ,
\end{equation}
which is an algebraic equation for $\beta$. In general, it can have multiple solutions, but one of them must be close to $0$ and it approaches $0$ as $\triangle t\rightarrow 0$. So we solve for this solution using an iterative method such as the Newton iteration or an optimization method with $0$ as the initial condition, it converges to a solution close to $0$ when $\triangle t$ is not too large. After obtaining $\beta$, we update $\bm\phi^{n+1}$ through \eqref{SVM update phi}.

\subsubsection{Existence of the supplementary variable}
The solution existence of the algebraic equation is guaranteed under the conditions given in the following theorems.
\begin{theorem}\label{SVM alpha existence}
If $\big(\bm\mu^n, g(\bm\phi^n)\big) \neq 0$, there exists $\triangle t^*$ such that \eqref{SVM algebraic equation} defines a unique function $\beta(\triangle t)$ for all $\triangle t\in[0,\triangle t^*]$ and scheme \eqref{SVM semi scheme} is 2nd order.
\begin{proof}
For $\triangle t,\beta$ in a neighborhood of $(0,0)$, we define the real valued function
\begin{equation}
u(\triangle	t,\beta) =
F[\widehat{\bm\phi}^{n+1} + \beta \check{\bm\phi}^{n+\frac{1}{2}}] - \widetilde F^{n+1}
= F[\widehat{\bm\phi}^{n+1} + \beta \check{\bm\phi}^{n+\frac{1}{2}}]
- F[\bm\phi^n] + \triangle t(\widetilde{\bm\mu}, \widetilde{\bm M}\widetilde{\bm\mu}) .
\end{equation}
It follows from \eqref{SVM predict mu}\eqref{SVM predict}\eqref{SVM procedure} that
\begin{equation}
u(0,0) = 0 ,
\quad
\frac{\partial u}{\partial\beta}(0,0)
= \left( \frac{\delta F}{\delta\bm\phi}(\bm\phi^n) , \bm g(\bm\phi^n) \right) \neq 0 .
\end{equation}
According to the implicit function theorem, there exists $\triangle t^* > 0$ such that equation $u(\triangle t,\beta)=0$ defines a unique smooth function $\beta=\beta(\triangle t)$ satisfying $\beta(0)=0$ and $u\big(\triangle t,\beta(\triangle t)\big)=0$ for all $\triangle t\in[0,\triangle t^*]$. Since $\widehat{\bm\phi}^{n+1}$ satisfies the following scheme
\begin{equation}
\frac{\widehat{\bm\phi}^{n+1}-\bm\phi^n}{\triangle t}
= - \bm M(\widetilde{\bm\phi}^{n+\frac{1}{2}})
\Big[ \mathcal L \frac{\widehat{\bm\phi}^{n+1}+\bm\phi^n}{2}
+ f'(\widetilde{\bm\phi}^{n+\frac{1}{2}}) \Big] ,
\end{equation}
through local truncation error analysis we have
\begin{equation}
\widehat{\bm\phi}^{n+1} = \bm\phi(t_{n+1}) + \mathcal O(\triangle t^3),
\quad
F[\widehat{\bm\phi}^{n+1}] = F[\bm\phi(t_{n+1})] + \mathcal O(\triangle t^3).
\end{equation}
In addition, we expand
\begin{equation}
u(\triangle t, \beta) = u(\triangle t, 0)
+ \beta\frac{\partial u}{\partial\beta}(\triangle t,0) + \mathcal O(\beta^2) ,
\end{equation}
with
\begin{equation}
u(\triangle t,0)=F[\widehat{\bm\phi}^{n+1}]-\widetilde F^{n+1}=\mathcal O(\triangle t^3) ,
\quad
\frac{\partial u}{\partial\beta}(\triangle t,0)
= \frac{\partial u}{\partial\beta}(0,0) + \mathcal O(\triangle t) .
\end{equation}
Then $\beta=\beta(\triangle t)=\mathcal O(\triangle t^3)$ and the proposed scheme is 2nd order.
\end{proof}
\end{theorem}

\begin{theorem}
If \eqref{SVM g2} is chosen i.e. $\bm g(\bm\phi)=-\bm M\bm\mu$, the sufficient condition in Theorem \ref{SVM alpha existence} reduces to $(\bm\mu^n, \bm M^n\bm\mu^n) \neq 0$, which is usually satisfied when the steady state is not reached.
\end{theorem}

\section{Fully discrete numerical algorithms}
We present fully discrete numerical schemes in this section. In this paper, we consider cases where $M_{ij}$ are all constants, where it is readily seen that the Lagrange multiplier for the volume fraction constraint is given by
\begin{equation}\label{explicit L}
L = - \frac{\sum_{i,j}M_{ij}\mu_i}{\sum_{i,j}M_{ij}} .
\end{equation}
We deduce from \eqref{CH} that
\begin{equation}
\left\{
\begin{aligned}
\frac{\partial\phi_A}{\partial t} &= \Delta( m_{11}\mu_A + m_{12}\mu_B + m_{13}\mu_S ) ,
\\
\frac{\partial\phi_B}{\partial t} &= \Delta( m_{21}\mu_A + m_{22}\mu_B + m_{23}\mu_S ) ,
\\
\frac{\partial\phi_S}{\partial t} &= \Delta( m_{31}\mu_A + m_{32}\mu_B + m_{33}\mu_S ) ,
\end{aligned}
\right.
\end{equation}
where
\begin{equation}
m_{kl} = M_{kl} - \big(\sum_jM_{kj}\big)\big(\sum_jM_{lj}\big)/\big(\sum_{i,j}M_{ij}\big) .
\end{equation}

\subsection{Numerical treatment of the logarithmic potential}
Following the work in \cite{Copetti_1992_CH_log_energy, YangXiaofeng_2019_CH_log_energy}, we regularize the logarithmic bulk potential by a $C^2$ piecewise function. More precisely, for any $0 < \sigma \ll 1$, the regularized function of
\begin{equation}
h(\phi) = \frac{\phi}{N}\ln\phi , \quad 0 < \phi < 1
\end{equation}
is defined by
\begin{equation}
\widehat h(\phi) =
\left\{
\begin{aligned}
&\frac{1}{N}\Big(\frac{\phi^2}{2\sigma}+\phi\ln\sigma-\frac{\sigma}{2}\Big), &&\phi\leq\sigma,
\\
&\frac{\phi}{N}\ln\phi, &&\phi\geq\sigma,
\end{aligned}
\right.
\end{equation}
In \cite{Copetti_1992_CH_log_energy}, the authors proved the error bound between the regularized PDE and the original PDE is controlled by $\sigma$ up to a constant. We use the modified logarithmic function with $\sigma = 0.01$ in all numerical simulations in this paper. However for simplicity, we retain the original notation in our presentation of schemes.

\subsection{Notations used in discretization}
We consider a time domain $[0,T]$ for $T\in\mathbb R^+$. For $N_t\in\mathbb Z^+$ we define the time step as $\triangle t=T/N_t$ and the $n$-th step as $t_n=n\triangle t~(n=0,\cdots,N_t)$. Denote
\begin{equation}
f^{n+\frac{1}{2}} = (f^n+f^{n+1})/2 .
\end{equation}
We define the temporal difference operator and extrapolation operator respectively as
\begin{equation}
\delta^+_t f^n = (f^{n+1}-f^n)/{\triangle t} ,
\quad
\overline f^{n+\frac{1}{2}} = (3f^n-f^{n-1})/2 .
\end{equation}

For simplicity, we present a spatial discretization in 2D rectangular domain $\Omega=[0,L_x]\times [0, L_y]$. Let $N_x, N_y$ be two positive integers. The domain $\Omega$ is uniformly partitioned with mesh sizes $h_x=L_x/N_x, ~h_y=L_y/N_y$. The cell centered grid points are denoted by
\begin{equation}
\Omega_h = \{(x_i,y_j) ~|~ x_i=(i-1/2)h_x, ~y_j=(j-1/2)h_y, ~1\leq i\leq N_x, ~1\leq j\leq N_y \}.
\end{equation}
We define the following average and finite difference operators
\begin{equation}
\begin{gathered}
A_x f_{ij} = \frac{f_{i-\frac{1}{2},j} + f_{i+\frac{1}{2},j}}{2} , ~
A_y f_{ij} = \frac{f_{i,j-\frac{1}{2}} + f_{i,j+\frac{1}{2}}}{2} , ~
D_x f_{ij} = \frac{f_{i+\frac{1}{2},j} - f_{i-\frac{1}{2},j}}{h_x} , ~
D_y f_{ij} = \frac{f_{i,j+\frac{1}{2}} - f_{i,j-\frac{1}{2}}}{h_y} ,
\\
d_x f_{i-\frac{1}{2},j} = \frac{f_{ij}-f_{i-1,j}}{h_x} , ~
d_y f_{i,j-\frac{1}{2}} = \frac{f_{ij}-f_{i,j-1}}{h_y} , ~
\Delta_h f_{ij} = (d_xD_x + d_yD_y)f_{ij} .
\end{gathered}
\end{equation}
The discrete inner product and norm are defined respectively by
\begin{equation}
(f,g)_h = \sum^{N_x}_{i=1} \sum^{N_y}_{j=1} f_{ij}g_{ij}h_xh_y ,
\quad
\|f\|_h = \sqrt{(f,f)_h} .
\end{equation}
The following summation-by-parts formulas are valid analogous to the integration-by-parts formulas
\begin{equation}
(f, D_\alpha g)_h + (D_\alpha f, g)_h = 0 ,
\quad
(f, d_\alpha g)_h + (d_\alpha f, g)_h = 0 ,
\quad \alpha = x, y ,
\end{equation}
when $f,g$ satisfy proper boundary conditions making the corresponding boundary integrals vanish.

On the cell-centered grid, the zero Neumann boundary conditions we use in this study can be easily implemented using the scheme and the discrete summation-by-parts formulae apply.

\subsection{EQ algorithm for diblock copolymer solutions}
We present a 2nd order, linear algorithm using the EQ method. We introduce an intermediate variable
\begin{equation}
q = \sqrt{f(\bm\phi) + C}
= \sqrt{\frac{\phi_A}{N_A}\ln\phi_A+\frac{\phi_B}{N_B}\ln\phi_B+\frac{\phi_S}{N_S}\ln\phi_S + C} ,
\end{equation}
where $C>0$ is large enough to make $q>0$. We have
\begin{equation}
q'(\bm\phi)=\Big[\frac{\partial q}{\partial\phi_A}, \frac{\partial q}{\partial\phi_B},
\frac{\partial q}{\partial\phi_S}\Big]^\top ,
\quad
\frac{\partial q}{\partial\phi_i}=\frac{1}{2q}\frac{1+\ln\phi_i}{N_i}, \quad i=A,B,S.
\end{equation}
For simplicity, we use notations
\begin{equation}\label{scheme notations}
\begin{gathered}
\bm\phi = [\phi_A, \phi_B, \phi_S]^\top ,
\quad
(m_{kl}) = \left[\begin{matrix}
m_{11} & m_{12} & m_{13} \\ m_{12} & m_{22} & m_{33} \\ m_{13} & m_{23} & m_{33}
\end{matrix}\right] ,
\quad
F_h[\bm\phi] = \frac{1}{2}(\bm\phi,\mathcal L\bm\phi)_h + \big(f(\bm\phi),1\big)_h ,
\\
\mathcal L\bm\phi =
\left[ \begin{matrix}
-\gamma_A\Delta\phi_A + \chi_{AA}\phi_A + \chi_{AB}\phi_B + \chi_{AS}\phi_S
- \alpha_{AA}\psi_A - \alpha_{AB}\psi_B
\\
-\gamma_B\Delta\phi_B + \chi_{AB}\phi_A + \chi_{BB}\phi_B + \chi_{BS}\phi_S
- \alpha_{AB}\psi_A - \alpha_{BB}\psi_B
\\
-\gamma_S\Delta\phi_S + \chi_{AS}\phi_A + \chi_{BS}\phi_B + \chi_{SS}\phi_S
\end{matrix} \right] .
\end{gathered}
\end{equation}

\begin{alg}[Fully discrete EQ algorithm]
Assuming $\phi^{n-1}_A, \phi^{n-1}_B, \phi^{n-1}_S, \phi^n_A, \phi^n_B, \phi^n_S ~(\phi^{n-1}_A+\phi^{n-1}_B+\phi^{n-1}_S=\phi^n_A+\phi^n_B+\phi^n_S=1)$ are known, we solve $\phi^{n+1}_A, \phi^{n+1}_B, \phi^{n+1}_S$ through
\begin{equation}
\left\{
\begin{aligned}
&\Delta_h(\psi_A)^{n+\frac{1}{2}}_{ij} = (\phi_A)^{n+\frac{1}{2}}_{ij} - \overline\phi_A ,
\quad
\Delta_h(\psi_B)^{n+\frac{1}{2}}_{ij} = (\phi_B)^{n+\frac{1}{2}}_{ij} - \overline\phi_B ,
\\
&\delta_t^+\bm\phi^n_{ij} = (m_{kl})\Delta_h \Big[ \mathcal L_h\bm\phi^{n+\frac{1}{2}}_{ij}
+ 2q^{n+\frac{1}{2}}_{ij}q'(\overline{\bm\phi}^{n+\frac{1}{2}}_{ij}) \Big] ,
\\
&\delta_t^+q^n_{ij} = q'(\overline{\bm\phi}^{n+\frac{1}{2}}_{ij})\cdot\delta_t^+\bm\phi^n_{ij} ,
\end{aligned}
\right.
\end{equation}
where $1\leq i\leq N_x, ~1\leq j\leq N_y$.
\end{alg}

\subsection{SVM algorithms for diblock copolymer solutions}
Here we present \textbf{SVM2} algorithm only. The other 3 SVM algorithms are developed analogously using the other 3 SVM strategies. With notations in \eqref{scheme notations}, we denote
\begin{equation}\label{constant_CH f_phi}
f'(\bm\phi)=\Big[\frac{\partial f}{\partial\phi_A}, \frac{\partial f}{\partial\phi_B},
\frac{\partial f}{\partial\phi_S}\Big]^\top ,
\quad
\frac{\partial f}{\partial\phi_i}=\frac{1+\ln\phi_i}{N_i}, ~ i=A,B,S.
\end{equation}
In \textbf{SVM2} strategy, we also use function
\begin{equation}
\bm g(\bm\phi) = (m_{kl})\Delta\big[\mathcal L\bm\phi+f'(\bm\phi)\big] .
\end{equation}
\begin{alg}[Fully discrete SVM algorithm]~
\\
\noindent\textbf{Prediction}:
Assuming $\phi^{n-1}_A, \phi^{n-1}_B, \phi^{n-1}_S, \phi^n_A, \phi^n_B, \phi^n_S ~(\phi^{n-1}_A+\phi^{n-1}_B+\phi^{n-1}_S=\phi^n_A+\phi^n_B+\phi^n_S=1), F^n$ are known, we compute $\widetilde\phi^{n+\frac{1}{2}}_A, \widetilde\phi^{n+\frac{1}{2}}_B, \widetilde\phi^{n+\frac{1}{2}}_S$ and $F^{n+1}$ by
\begin{equation}
\left\{
\begin{aligned}
& \Delta_h(\widetilde\psi_A)^{n+\frac{1}{2}}_{ij}
= (\widetilde\phi_A)^{n+\frac{1}{2}}_{ij} - \overline\phi_A ,
\quad
\Delta_h(\widetilde\psi_B)^{n+\frac{1}{2}}_{ij}
= (\widetilde\phi_B)^{n+\frac{1}{2}}_{ij} - \overline\phi_B ,
\\
& \frac{\widetilde{\bm\phi}^{n+\frac{1}{2}}_{ij}-\bm\phi^n_{ij}}{\triangle t/2} =
(m_{kl})\Delta_h\Big[ \mathcal L_h\widetilde{\bm\phi}^{n+\frac{1}{2}}_{ij}
+ f'(\overline{\bm\phi}^{n+\frac{1}{2}}_{ij}) \Big] ,
\\
& \frac{F^{n+1}-F^n}{\triangle t} =
\Big( \mathcal L\widetilde{\bm\phi}^{n+\frac{1}{2}}+f'(\widetilde{\bm\phi}^{n+\frac{1}{2}}),
(m_{kl})\Delta
\big[ \mathcal L\widetilde{\bm\phi}^{n+\frac{1}{2}}+f'(\widetilde{\bm\phi}^{n+\frac{1}{2}}) \big]
\Big)_h ;
\end{aligned}
\right.
\end{equation}
\noindent\textbf{Correction}: We compute $\phi^{n+1}_A, \phi^{n+1}_B, \phi^{n+1}_S$ and $\alpha^{n+\frac{1}{2}}$ by
\begin{equation}
\left\{
\begin{aligned}
&\Delta_h(\psi_A)^{n+\frac{1}{2}}_{ij} = (\phi_A)^{n+\frac{1}{2}}_{ij} - \overline\phi_A ,
\quad
\Delta_h(\psi_B)^{n+\frac{1}{2}}_{ij} = (\phi_B)^{n+\frac{1}{2}}_{ij} - \overline\phi_B ,
\\
&\delta_t^+\bm\phi^n_{ij} = (m_{kl})\Delta_h
\Big[ \mathcal L_h\bm\phi^{n+\frac{1}{2}}_{ij} + f'(\widetilde{\bm\phi}^{n+\frac{1}{2}}_{ij}) \Big]
+ \alpha^{n+\frac{1}{2}}\bm g(\widetilde{\bm\phi}^{n+\frac{1}{2}}_{ij}) ,
\\
&F_h[\bm\phi^{n+1}] = F^{n+1} ,
\end{aligned}
\right.
\end{equation}
where $1\leq i\leq N_x, ~1\leq j\leq N_y$.
\end{alg}

\subsection{SVM algorithm for the electric field coupled model}
In the electric field coupled model, there is an additional Poisson equation for the electric potential $\Phi$. With \eqref{scheme notations}, we add the following notations in the description of the scheme.
\begin{equation}
\begin{gathered}
(\overline\mu_{\rm e})^{n+\frac{1}{2}}_{ij} =
- \frac{\epsilon_1}{2}\Big|\bm E_0-[d_xA_x,d_yA_y]^\top\overline\Phi^{n+\frac{1}{2}}_{ij}\Big|^2 ,
\quad
(\widetilde\mu_{\rm e})^{n+\frac{1}{2}}_{ij} =
- \frac{\epsilon_1}{2}\Big|\bm E_0-[d_xA_x,d_yA_y]^\top\widetilde\Phi^{n+\frac{1}{2}}_{ij}\Big|^2 ,
\\
\overline{\bm h}^{n+\frac{1}{2}}_{ij} = f'(\overline{\bm\phi}^{n+\frac{1}{2}}_{ij})
+ \Big[ (\overline\mu_{\rm e})^{n+\frac{1}{2}}_{ij},
- (\overline\mu_{\rm e})^{n+\frac{1}{2}}_{ij}, 0 \Big]^\top ,
\quad
\widetilde{\bm h}^{n+\frac{1}{2}}_{ij} = f'(\widetilde{\bm\phi}^{n+\frac{1}{2}}_{ij})
+ \Big[ (\widetilde\mu_{\rm e})^{n+\frac{1}{2}}_{ij},
- (\widetilde\mu_{\rm e})^{n+\frac{1}{2}}_{ij}, 0 \Big]^\top ,
\\
\bm g(\widetilde{\bm\phi}) = (m_{kl})\Delta(\mathcal L\widetilde{\bm\phi}+\widetilde{\bm h}) ,
\quad
E_h[\bm\phi,\Phi] = F_h[\bm\phi]
+ \frac{1}{2}\big( \epsilon_0+\epsilon_1(\phi_A-\phi_B-\overline\phi_A+\overline\phi_B) ,
|\bm E_0-\nabla\Phi|^2 \big)_h .
\end{gathered}
\end{equation}
\begin{alg}[Fully discrete SVM algorithm for the electric field coupled model]~
\\
\noindent\textbf{Prediction}:
Assuming $\phi^{n-1}_A, \phi^{n-1}_B, \phi^{n-1}_S, \phi^n_A, \phi^n_B, \phi^n_S ~(\phi^{n-1}_A+\phi^{n-1}_B+\phi^{n-1}_S=\phi^n_A+\phi^n_B+\phi^n_S=1), \Phi^n, E^n$ are known, we compute $\overline\Phi^{n+\frac{1}{2}}, \widetilde\phi^{n+\frac{1}{2}}_A, \widetilde\phi^{n+\frac{1}{2}}_B, \widetilde\phi^{n+\frac{1}{2}}_S$ and $E^{n+1}$ by
\begin{equation}
\left\{
\begin{aligned}
& \Big[\! \frac{\epsilon_0}{\epsilon_1} \!-\! \overline\phi_A \!+\! \overline\phi_B \!+\!
(\overline\phi_A)^{n\!+\!\frac{1}{2}}_{ij} \!\!-\!\! (\overline\phi_B)^{n\!+\!\frac{1}{2}}_{ij}
\!\Big] \! \Delta_h \! \overline\Phi^{n\!+\!\frac{1}{2}}_{ij}
\!\!=\!\! \big(\bm E_0\!-\![d_xA_x,\!d_yA_y]^\top\! \overline\Phi^{n\!+\!\frac{1}{2}}_{ij}\big)
\!\cdot\! [d_xA_x,\!d_yA_y]^\top\! \Big[\! (\overline\phi_A)^{n\!+\!\frac{1}{2}}_{ij}
\!\!-\!\! (\overline\phi_B)^{n\!+\!\frac{1}{2}}_{ij} \!\Big] ,
\\
& \Delta_h(\widetilde\psi_A)^{n+\frac{1}{2}}_{ij}
= (\widetilde\phi_A)^{n+\frac{1}{2}}_{ij} - \overline\phi_A ,
\quad
\Delta_h(\widetilde\psi_B)^{n+\frac{1}{2}}_{ij}
= (\widetilde\phi_B)^{n+\frac{1}{2}}_{ij} - \overline\phi_B ,
\\
& \frac{\widetilde{\bm\phi}^{n+\frac{1}{2}}_{ij}-\bm\phi^n_{ij}}{\triangle t/2} =
(m_{kl})\Delta_h\Big[ \mathcal L_h\widetilde{\bm\phi}^{n+\frac{1}{2}}_{ij}
+ \overline{\bm h}^{n+\frac{1}{2}}_{ij}) \Big] ,
\\
& \frac{E^{n+1}-E^n}{\triangle t} =
\Big( \mathcal L\widetilde{\bm\phi}^{n+\frac{1}{2}}+\widetilde{\bm h}^{n+\frac{1}{2}},
(m_{kl})\Delta
\big[\mathcal L\widetilde{\bm\phi}^{n+\frac{1}{2}}+\widetilde{\bm h}^{n+\frac{1}{2}}\big] \Big)_h ;
\end{aligned}
\right.
\end{equation}
\textbf{Correction}: We compute $\widetilde\Phi^{n+\frac{1}{2}}, \phi^{n+1}_A, \phi^{n+1}_B, \phi^{n+1}_S$ and $\alpha^{n+\frac{1}{2}}$ by
\begin{equation}
\left\{
\begin{aligned}
& \Big[\! \frac{\epsilon_0}{\epsilon_1} \!-\! \overline\phi_A \!+\! \overline\phi_B \!+\!
(\widetilde\phi_A)^{n\!+\!\frac{1}{2}}_{ij} \!\!-\!\! (\widetilde\phi_B)^{n\!+\!\frac{1}{2}}_{ij}
\!\Big] \! \Delta_h \! \widetilde\Phi^{n\!+\!\frac{1}{2}}_{ij}
\!\!=\!\! \big(\bm E_0\!-\![d_xA_x,\!d_yA_y]^\top\! \widetilde\Phi^{n\!+\!\frac{1}{2}}_{ij}\big)
\!\cdot\! [d_xA_x,\!d_yA_y]^\top\! \Big[\! (\widetilde\phi_A)^{n\!+\!\frac{1}{2}}_{ij}
\!\!-\!\! (\widetilde\phi_B)^{n\!+\!\frac{1}{2}}_{ij} \!\Big] ,
\\
&\Delta_h(\psi_A)^{n+\frac{1}{2}}_{ij} = (\phi_A)^{n+\frac{1}{2}}_{ij} - \overline\phi_A ,
\quad
\Delta_h(\psi_B)^{n+\frac{1}{2}}_{ij} = (\phi_B)^{n+\frac{1}{2}}_{ij} - \overline\phi_B ,
\\
&\delta_t^+\bm\phi^n_{ij} = (m_{kl})\Delta_h
\Big[ \mathcal L_h\bm\phi^{n+\frac{1}{2}}_{ij}+\widetilde{\bm h}^{n+\frac{1}{2}}_{ij} \Big]
+ \alpha^{n+\frac{1}{2}}\bm g(\widetilde{\bm\phi}^{n+\frac{1}{2}}_{ij}) ,
\\
&E_h[\bm\phi^{n+1}] = E^{n+1} ,
\end{aligned}
\right.
\end{equation}
where $1\leq i\leq N_x, ~1\leq j\leq N_y$.
\end{alg}

\subsection{SVM algorithm for the magnetic field coupled model}
In the magnetic field coupled model, the magnetic field is assumed a prescribed field. With \eqref{scheme notations}, we add the following notations in the description of the scheme.
\begin{equation}
\begin{gathered}
(\overline\mu_{\rm m})^{n+\frac{1}{2}}_{ij} = - \gamma_m\Big[
(B^2_1d_xD_x+B^2_2d_yD_y)(\overline\phi_A^{n+\frac{1}{2}}-\overline\phi_B^{n+\frac{1}{2}})_{ij}
+ 2B_1B_2d_xd_yA_xA_y(\overline\phi_A^{n+\frac{1}{2}}-\overline\phi_B^{n+\frac{1}{2}})_{ij} \Big] ,
\\
(\widetilde\mu_{\rm m})^{n+\frac{1}{2}}_{ij} = - \gamma_m\Big[
(B^2_1d_xD_x+B^2_2d_yD_y)(\widetilde\phi_A^{n+\frac{1}{2}}-\widetilde\phi_B^{n+\frac{1}{2}})_{ij}
+ 2B_1B_2d_xd_yA_xA_y(\widetilde\phi_A^{n+\frac{1}{2}}-\widetilde\phi_B^{n+\frac{1}{2}})_{ij} \Big] ,
\\
\overline{\bm h}^{n+\frac{1}{2}}_{ij} = f'(\overline{\bm\phi}^{n+\frac{1}{2}}_{ij})
+ \Big[ (\overline\mu_{\rm m})^{n+\frac{1}{2}}_{ij},
- (\overline\mu_{\rm m})^{n+\frac{1}{2}}_{ij}, 0 \Big]^\top ,
\quad
\widetilde{\bm h}^{n+\frac{1}{2}}_{ij} = f'(\widetilde{\bm\phi}^{n+\frac{1}{2}}_{ij})
+ \Big[ (\widetilde\mu_{\rm m})^{n+\frac{1}{2}}_{ij},
- (\widetilde\mu_{\rm m})^{n+\frac{1}{2}}_{ij}, 0 \Big]^\top ,
\\
\bm g(\widetilde{\bm\phi}) = (m_{kl})\Delta(\mathcal L\widetilde{\bm\phi}+\widetilde{\bm h}) ,
\quad
E_h[\bm\phi] = F_h[\bm\phi]
+ \frac{1}{2}\big(\gamma_m,|\bm B_0\cdot\nabla(\phi_A-\phi_B)|^2\big)_h .
\end{gathered}
\end{equation}

\begin{alg}[Fully discrete SVM algorithm for the magnetic field coupled model]~
\\
\noindent\textbf{Prediction}:
Assuming $\phi^{n-1}_A, \phi^{n-1}_B, \phi^{n-1}_S, \phi^n_A, \phi^n_B, \phi^n_S ~(\phi^{n-1}_A+\phi^{n-1}_B+\phi^{n-1}_S=\phi^n_A+\phi^n_B+\phi^n_S=1), E^n$ are known, we compute $\widetilde\phi^{n+\frac{1}{2}}_A, \widetilde\phi^{n+\frac{1}{2}}_B, \widetilde\phi^{n+\frac{1}{2}}_S$ and $E^{n+1}$ by
\begin{equation}
\left\{
\begin{aligned}
& \Delta_h(\widetilde\psi_A)^{n+\frac{1}{2}}_{ij}
= (\widetilde\phi_A)^{n+\frac{1}{2}}_{ij} - \overline\phi_A ,
\quad
\Delta_h(\widetilde\psi_B)^{n+\frac{1}{2}}_{ij}
= (\widetilde\phi_B)^{n+\frac{1}{2}}_{ij} - \overline\phi_B ,
\\
& \frac{\widetilde{\bm\phi}^{n+\frac{1}{2}}_{ij}-\bm\phi^n_{ij}}{\triangle t/2} =
(m_{kl})\Delta_h\Big[ \mathcal L_h\widetilde{\bm\phi}^{n+\frac{1}{2}}_{ij}
+ \overline{\bm h}^{n+\frac{1}{2}}_{ij}) \Big] ,
\\
& \frac{E^{n+1}-E^n}{\triangle t} =
\Big( \mathcal L\widetilde{\bm\phi}^{n+\frac{1}{2}}+\widetilde{\bm h}^{n+\frac{1}{2}},
(m_{kl})\Delta
\big[\mathcal L\widetilde{\bm\phi}^{n+\frac{1}{2}}+\widetilde{\bm h}^{n+\frac{1}{2}}\big] \Big)_h ;
\end{aligned}
\right.
\end{equation}
\textbf{Correction}: We compute $\phi^{n+1}_A, \phi^{n+1}_B, \phi^{n+1}_S$ and $\alpha^{n+\frac{1}{2}}$ by
\begin{equation}
\left\{
\begin{aligned}
&\Delta_h(\psi_A)^{n+\frac{1}{2}}_{ij} = (\phi_A)^{n+\frac{1}{2}}_{ij} - \overline\phi_A ,
\quad
\Delta_h(\psi_B)^{n+\frac{1}{2}}_{ij} = (\phi_B)^{n+\frac{1}{2}}_{ij} - \overline\phi_B ,
\\
&\delta_t^+\bm\phi^n_{ij} = (m_{kl})\Delta_h
\Big[ \mathcal L_h\bm\phi^{n+\frac{1}{2}}_{ij}+\widetilde{\bm h}^{n+\frac{1}{2}}_{ij} \Big]
+ \alpha^{n+\frac{1}{2}}\bm g(\widetilde{\bm\phi}^{n+\frac{1}{2}}_{ij}) ,
\\
&E_h[\bm\phi^{n+1}] = E^{n+1} ,
\end{aligned}
\right.
\end{equation}
where $1\leq i\leq N_x, ~1\leq j\leq N_y$.
\end{alg}

\subsection{Energy-dissipation-rate preservation and remarks}
Both the EQ and SVM algorithms are linear, 2nd order accurate and energy-dissipation-rate preserving. The proof for the EQ algorithm follows the standard routine \cite{WangQi_2018_EQ}, while the SVM algorithms preserve the energy-dissipation-rate by design. In fact, there are various formulations of the equations for the SVM strategy. In the presented algorithm, the dissipation rate is evaluated at the predicted value. We might as well use the final updated value in the energy-dissipation-rate without altering the order of the scheme. This demonstrates that a non-structure-preservation scheme can be turned into a structure-preserving scheme in the context of SVM through a correction in the order of the truncation error for any numerical schemes.

For variable mobility i.e. $M_{ij}(\bm\phi)$ are functions of phase variables, the Lagrange multiplier can no longer be given by an explicit formula like \eqref{explicit L}. So, it must be solved from the elliptic equation about $L$ in \eqref{CH}. We can propose EQ and SVM algorithms for that model as well although they are more complicated. Hence, that model along with the numerical treatment warrants a separate discussion.

If we replace all $\overline\square^{n+\frac{1}{2}}$ with $\square^{n+\frac{1}{2}}$ in Crank-Nicolson schemes, we obtain a 2nd order nonlinear energy stable scheme. If we replace all $\overline\square^{n+\frac{1}{2}}$ with $\square^n$, we obtain a two-level energy stable scheme, which is still linear but is of first order in time. In numerical experiments, we use the two-level scheme to compute the initial data for the second level values of the three-level schemes. This does not affect the overall accuracy of our 2nd order schemes. In addition, 2nd order BDF schemes can be used for the temporal discretization as well to achieve structure preservation in the SVM and EQ algorithms.

\section{Numerical results}
We first verify the order of the schemes in time and space using mesh refinement tests. Then we use \textbf{SVM2} algorithm to investigate (1). how the cross-coupling effect in the mobility affects transient dynamics in phase separation; (2). how the electric field accelerates phase separation in the copolymer system and removes transient defects from the material's domain; (3). how the magnetic field affects phase separation and defect dynamics. We use dimensionless domain $\Omega=[0,1]\times[0,1]$ and take $h_x=h_y=h=1/N$ in all numerical experiments.

\subsection{Mesh refinement test and computational efficiency}

\begin{figure}[H]
\centering
\subfigure[\textbf{EQ}]{\includegraphics[width=0.3\textwidth]{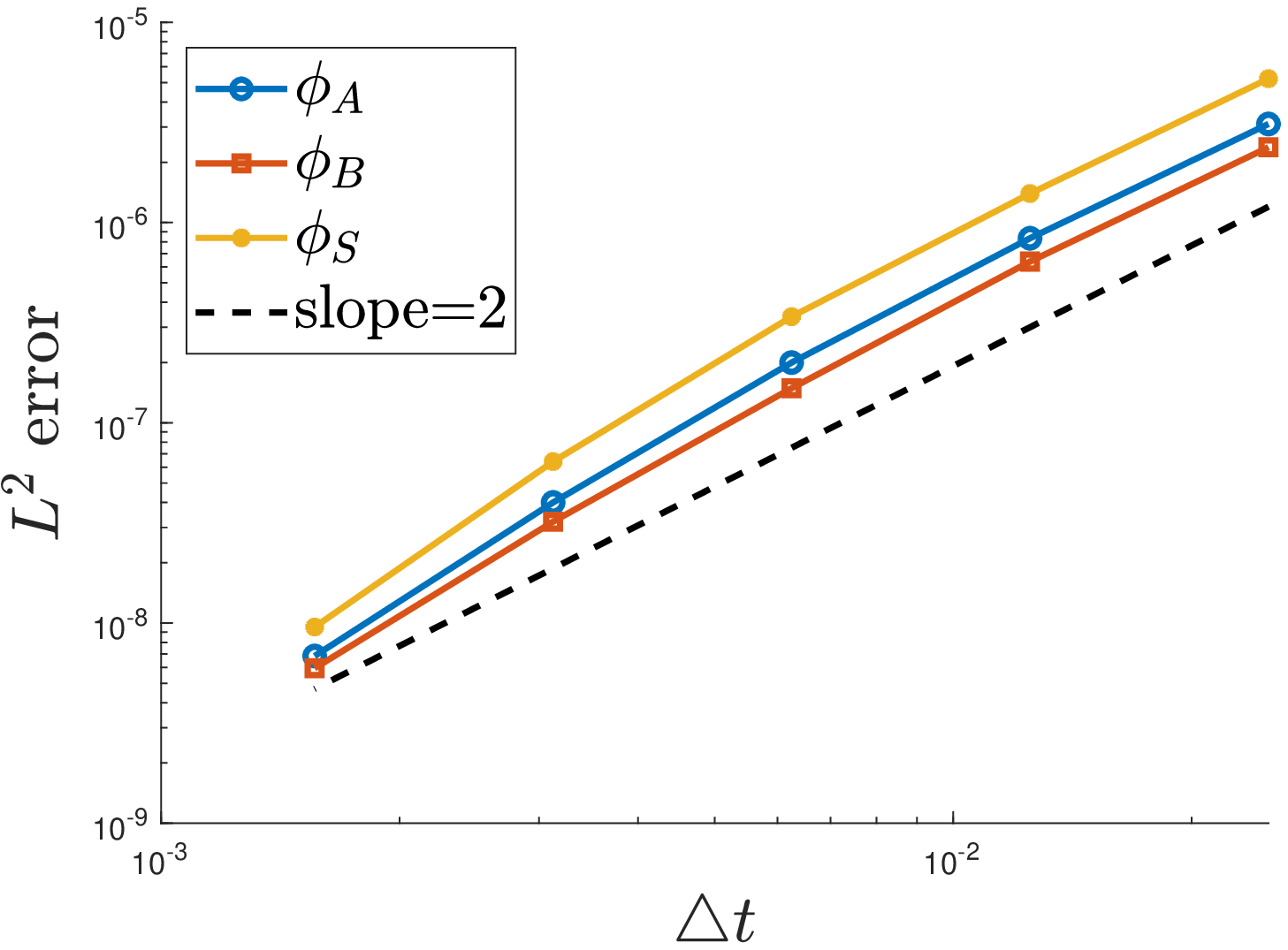}}
\subfigure[\textbf{SVM1}]{\includegraphics[width=0.3\textwidth]{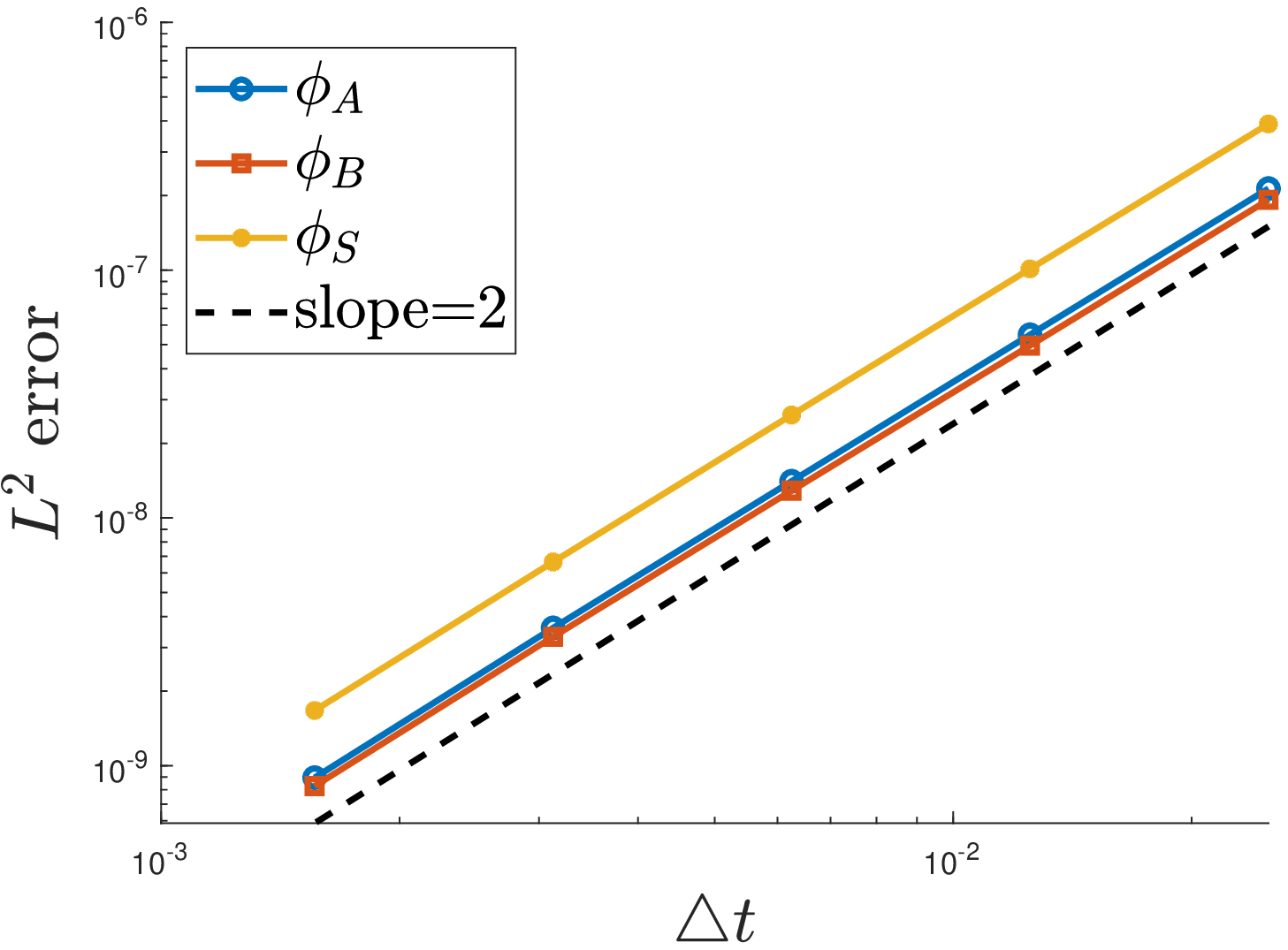}}
\subfigure[\textbf{SVM2}]{\includegraphics[width=0.3\textwidth]{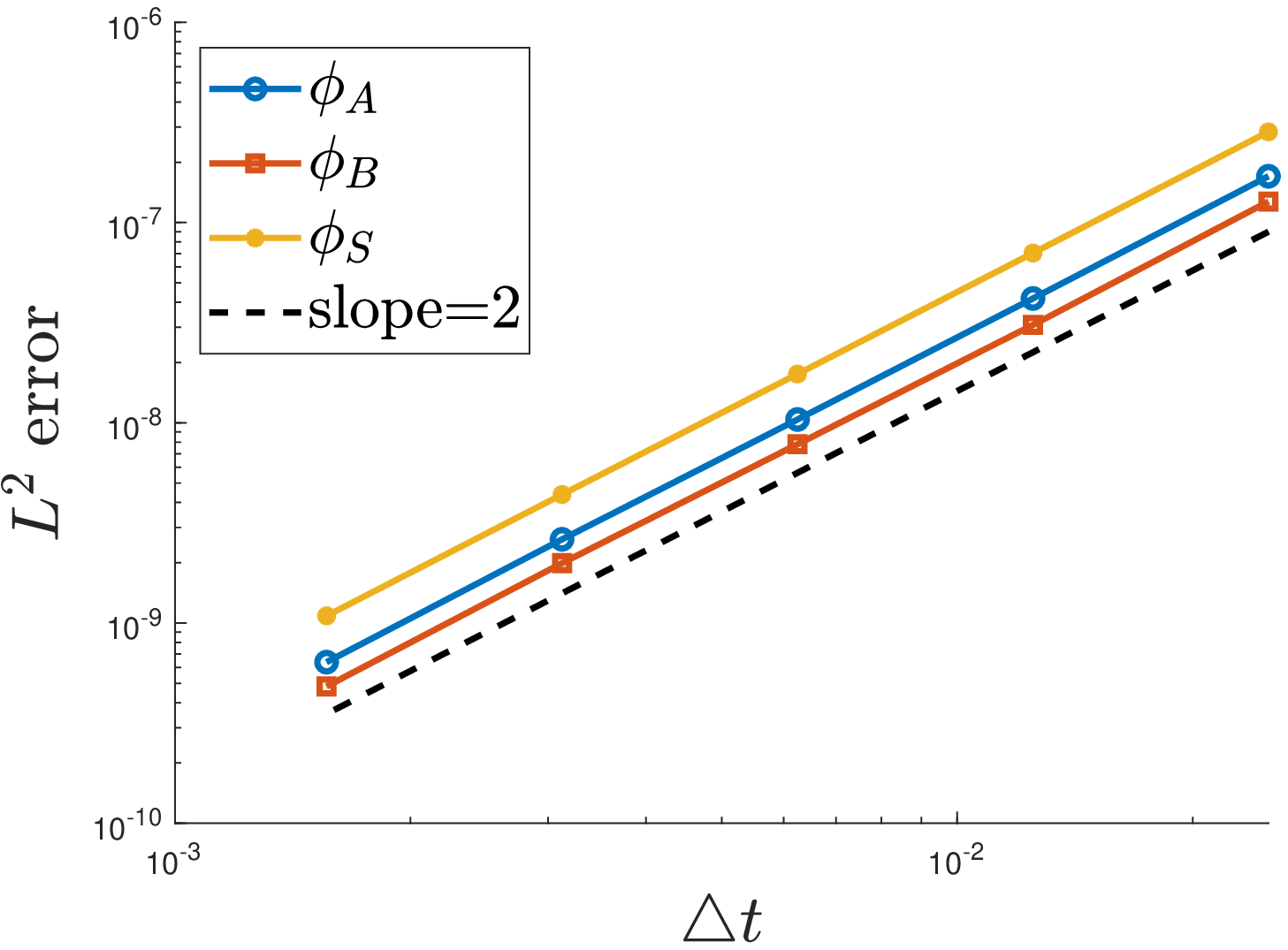}}
\\ \vspace{-0.3cm}
\subfigure[\textbf{SVM3}]{\includegraphics[width=0.3\textwidth]{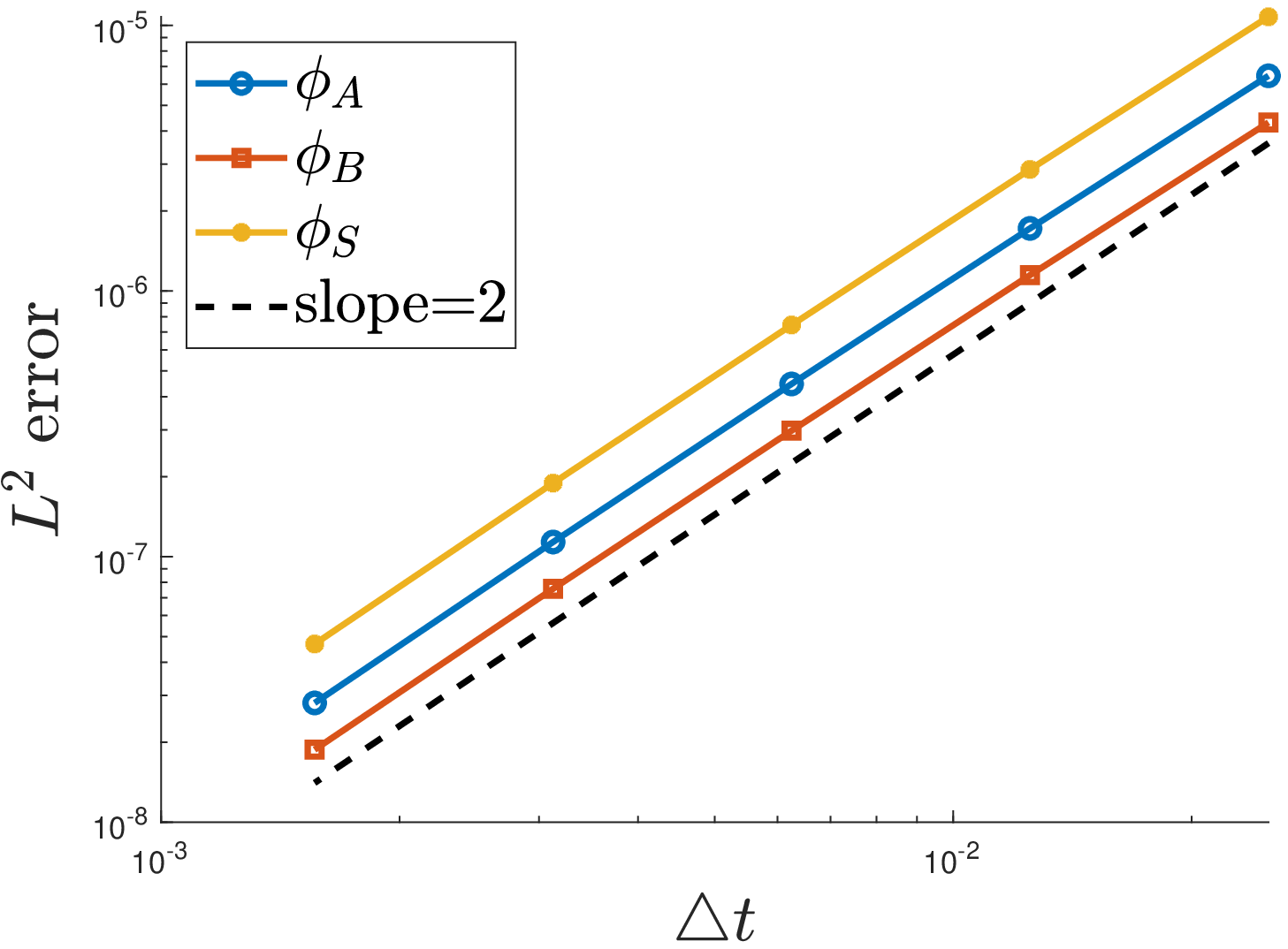}}
\subfigure[\textbf{SVM4}]{\includegraphics[width=0.3\textwidth]{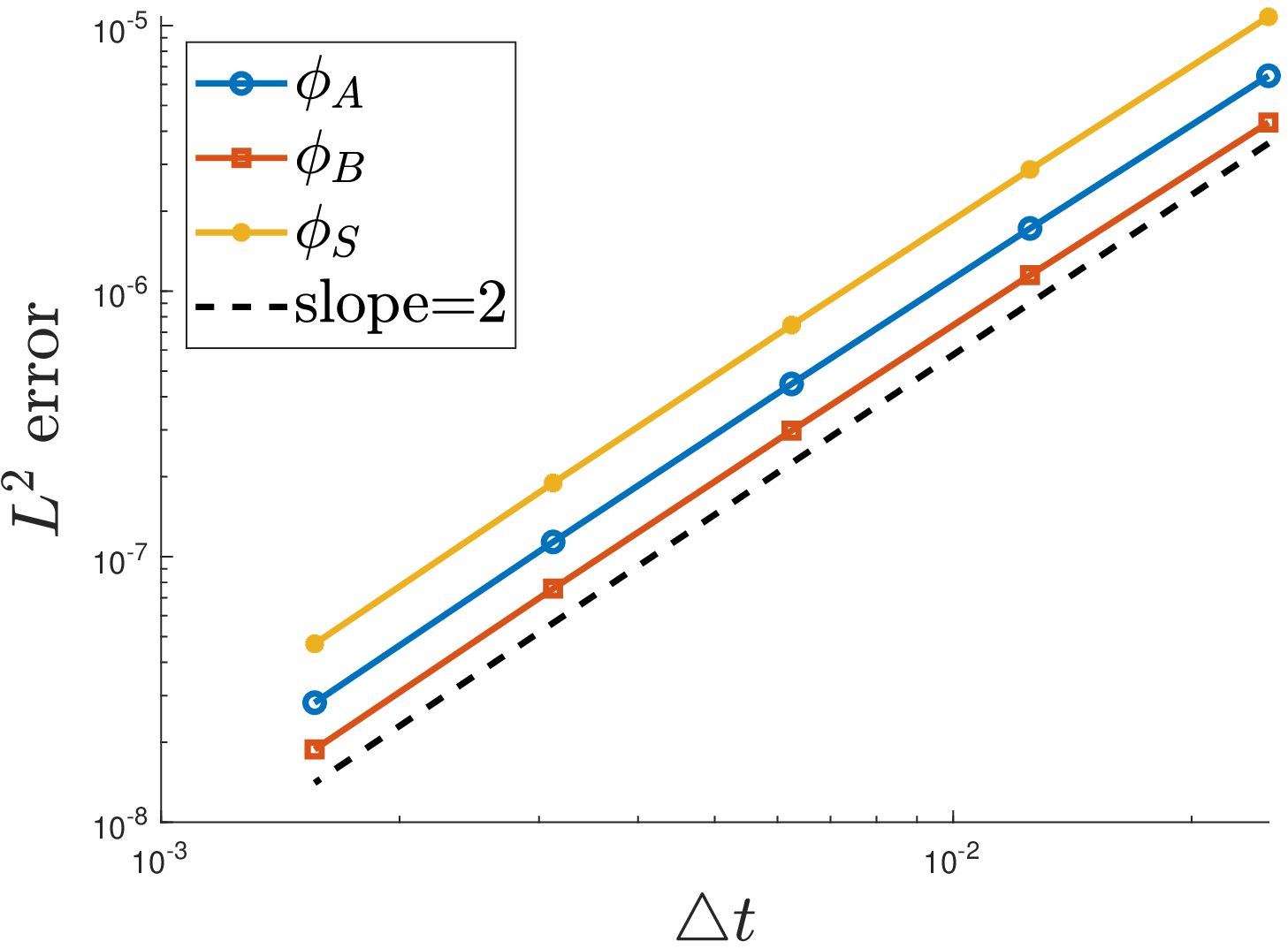}}
\\ \vspace{-0.3cm}
\caption{Temporal mesh refinement tests. The solution is computed up to $t=1$ with spatial mesh size $h=1/256$ and variable time steps. Subplots show the mesh refinement test in time for the convergence rate of the \textbf{EQ, SVM1, SVM2, SVM3, SVM4} scheme, respectively. \label{constant_CH_accuracy_time} }
\end{figure}

\begin{figure}[H]
\centering
\subfigure[\textbf{EQ}]{\includegraphics[width=0.31\textwidth]{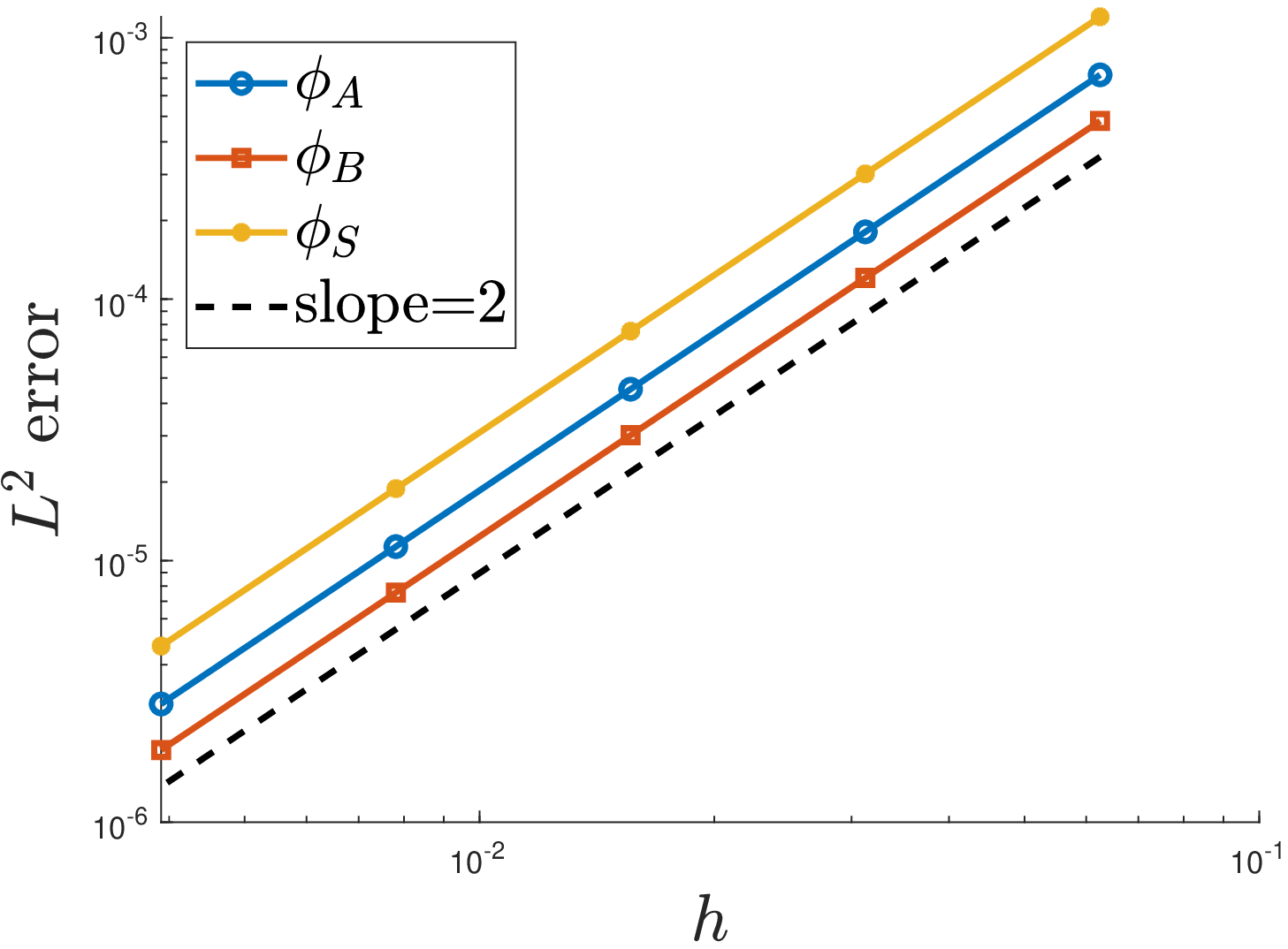}}
\subfigure[\textbf{SVM1}]
{\includegraphics[width=0.3\textwidth]{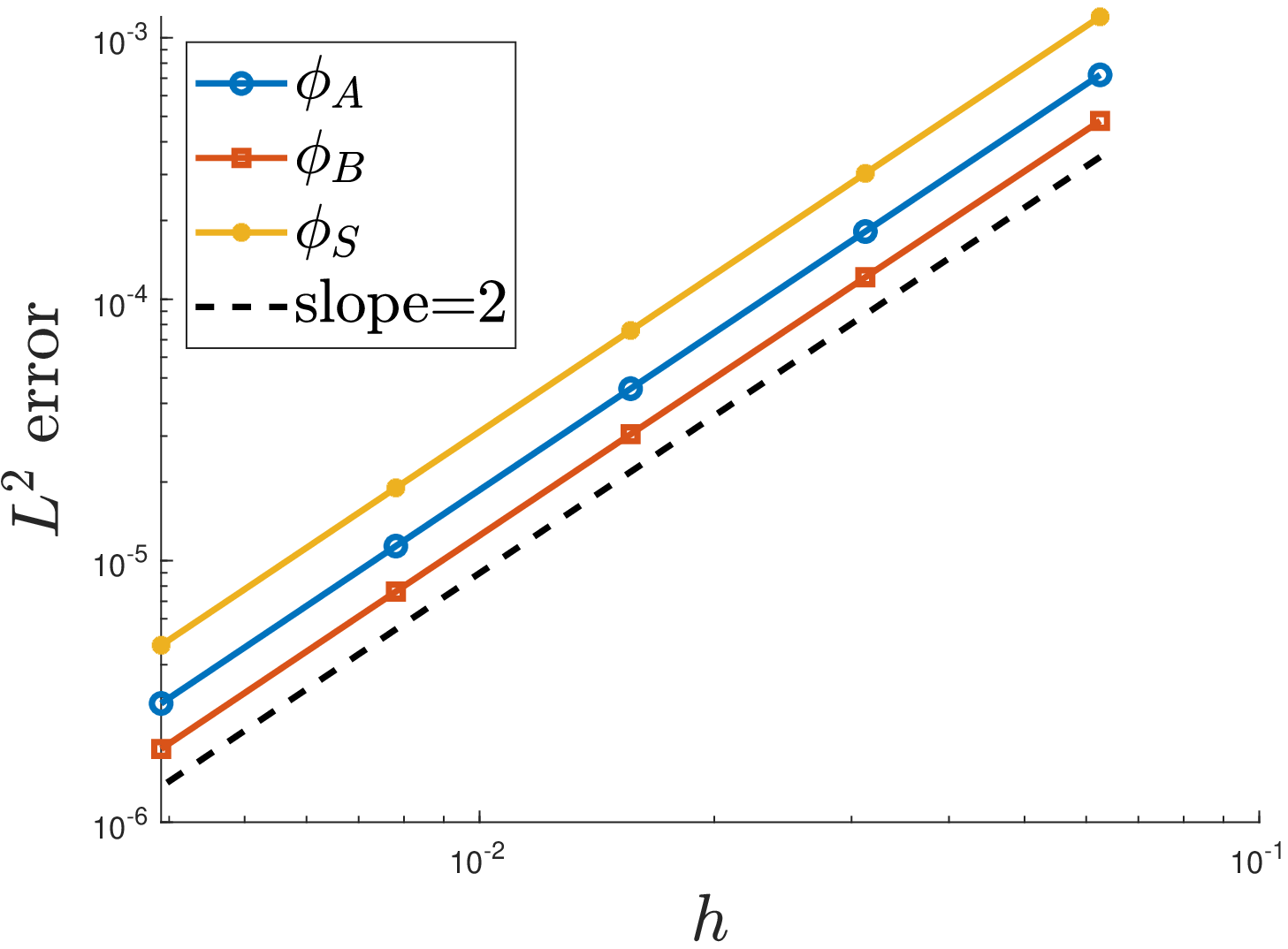}}
\subfigure[\textbf{SVM2}]
{\includegraphics[width=0.3\textwidth]{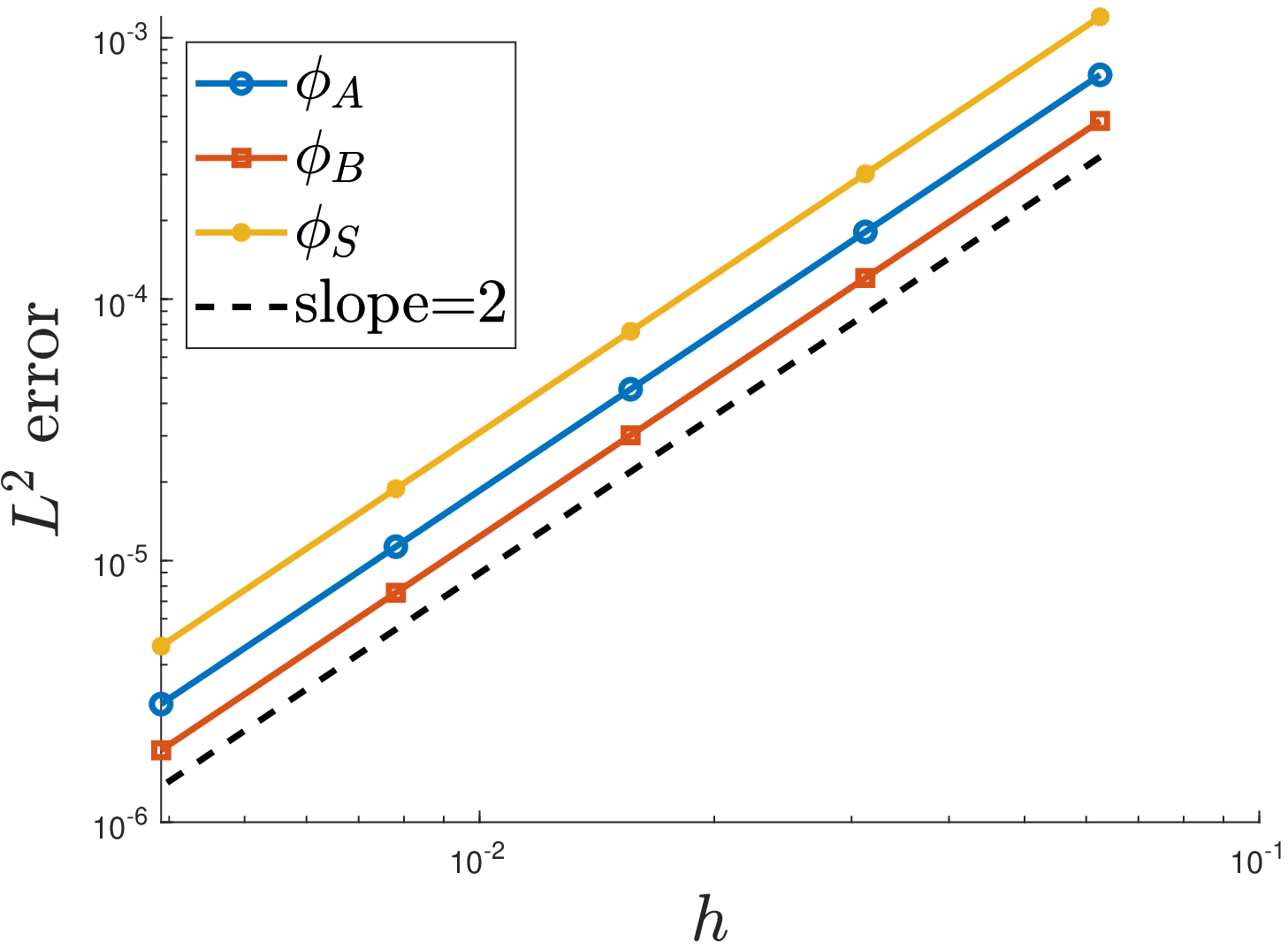}}
\\ \vspace{-0.3cm}
\subfigure[\textbf{SVM3}]
{\includegraphics[width=0.3\textwidth]{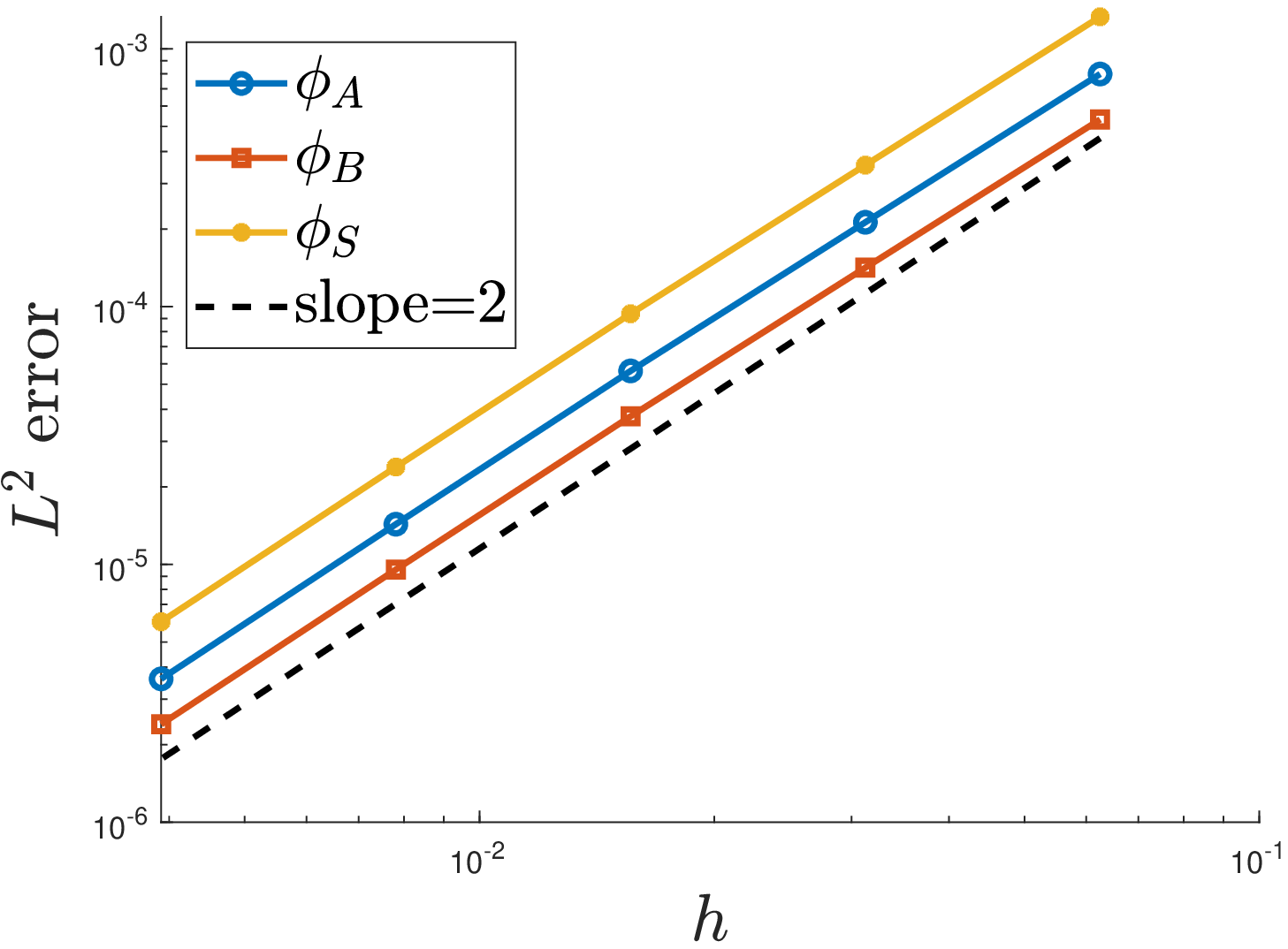}}
\subfigure[\textbf{SVM4}]
{\includegraphics[width=0.3\textwidth]{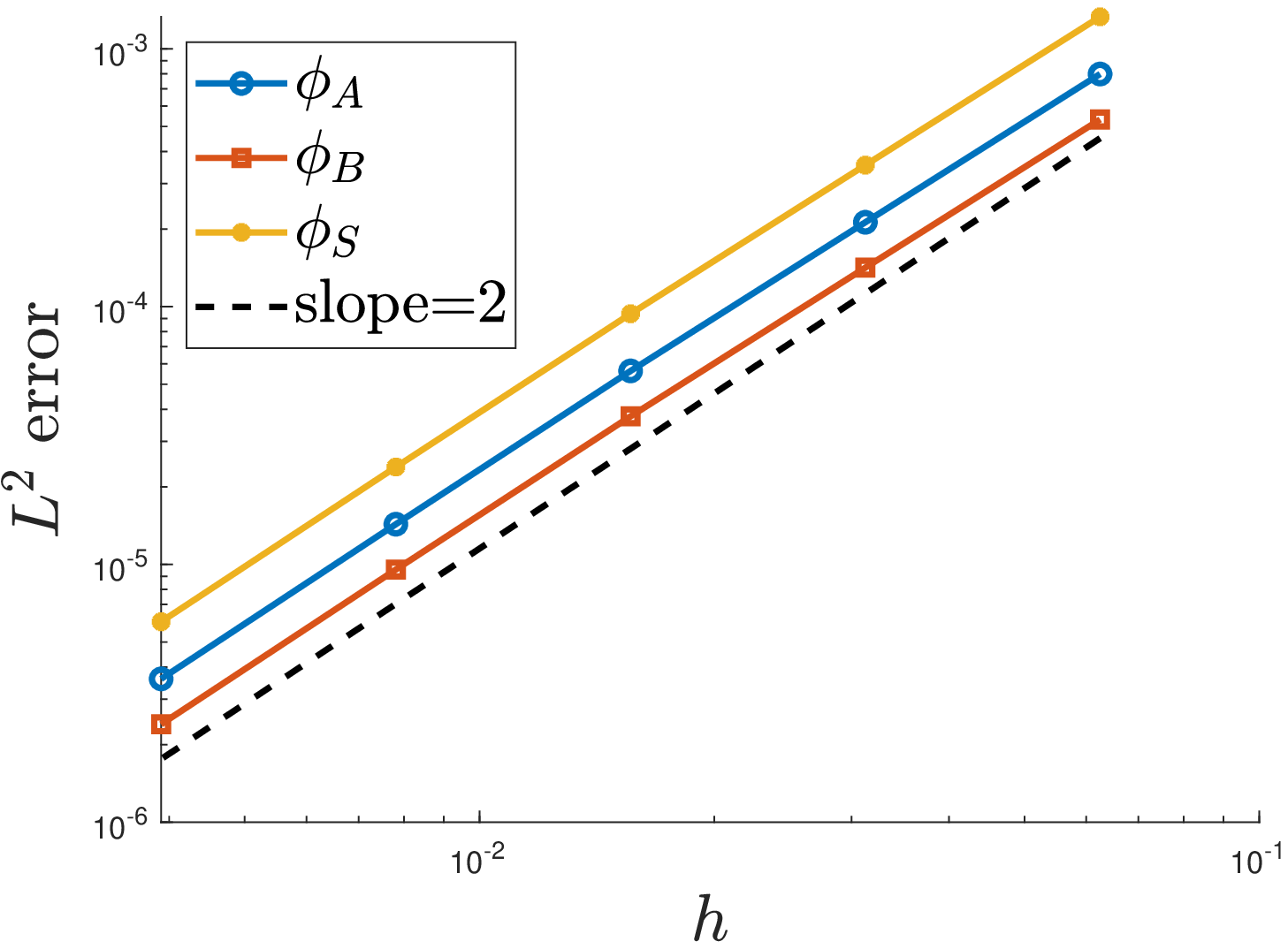}}
\\ \vspace{-0.3cm}
\caption{Spatial mesh refinement tests. The solution is computed up to $t=1$ at fixed time step $\triangle t=10^{-4}$. Subplots show the mesh refinement test in space for the convergence rate of the \textbf{EQ, SVM1, SVM2, SVM3, SVM4} scheme, respectively. \label{constant_CH_accuracy_space} }
\end{figure}

\begin{example}[\textbf{Mesh refinement}]\label{constant_CH mesh refinement}
We set $N_A=3, ~N_B=2, ~N_S=1, ~\chi_{AA}=\chi_{BB}=\chi_{SS}=0, ~\chi_{AB}=2, ~\chi_{AS}=3, ~\chi_{BS}=4, ~\varepsilon=0.1, ~\gamma=1$. The mobility matrix and initial conditions are given respectively by
\begin{equation}
M=10^{-5}\!\times\!\left[\begin{matrix} 4 & 1 & 2 \\ 1 & 5 & 3 \\ 2 & 3 & 6 \end{matrix}\right]\!,
\quad
\phi^0_A = 0.3(1+\cos\pi x\cos\pi y) ,
~
\phi^0_B = 0.2(1+\cos\pi x\cos\pi y) ,
~
\phi^0_S = 1 - \phi^0_A - \phi^0_B .
\end{equation}
\end{example}
First, we conduct the temporal accuracy test with spatial mesh size $h=1/256$ and time step $\triangle t=0.05\times 2^{-k},~k=0,\cdots,5$ respectively. We then conduct the spatial accuracy test with time step $\triangle t=10^{-4}$ and spatial mesh size $h=1/8\times 2^{-k},~k=0,\cdots,5$ respectively. The errors are calculated as the differences between the coarser step and the adjacent finer step. The discrete $L^2$ errors at $t=1$ are shown in Figures \ref{constant_CH_accuracy_time}, \ref{constant_CH_accuracy_space} respectively. The results demonstrate that all proposed schemes yield 2nd order convergence rates in both time and space.

\begin{figure}[H]
\centering
\includegraphics[width=0.35\textwidth]{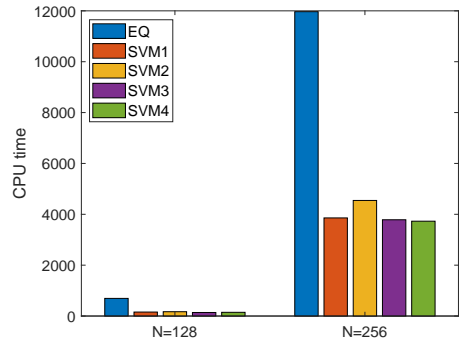}
\\ \vspace{-0.3cm}
\caption{For this model, \textbf{SVM} algorithms are more efficient than the \textbf{EQ} algorithm. In \textbf{SVM} algorithms, one solves a couple of linear system with all constant coefficients and a scalar nonlinear equation of the supplementary variable $\alpha$. Whereas, in \textbf{EQ} algorithm, one solves a linear system with variable coefficients. \label{constant_CH_accuracy_CPU} }
\end{figure}

\begin{figure}[H]
\centering
\subfigure[$h=1/256$]{\includegraphics[width=0.3\textwidth]{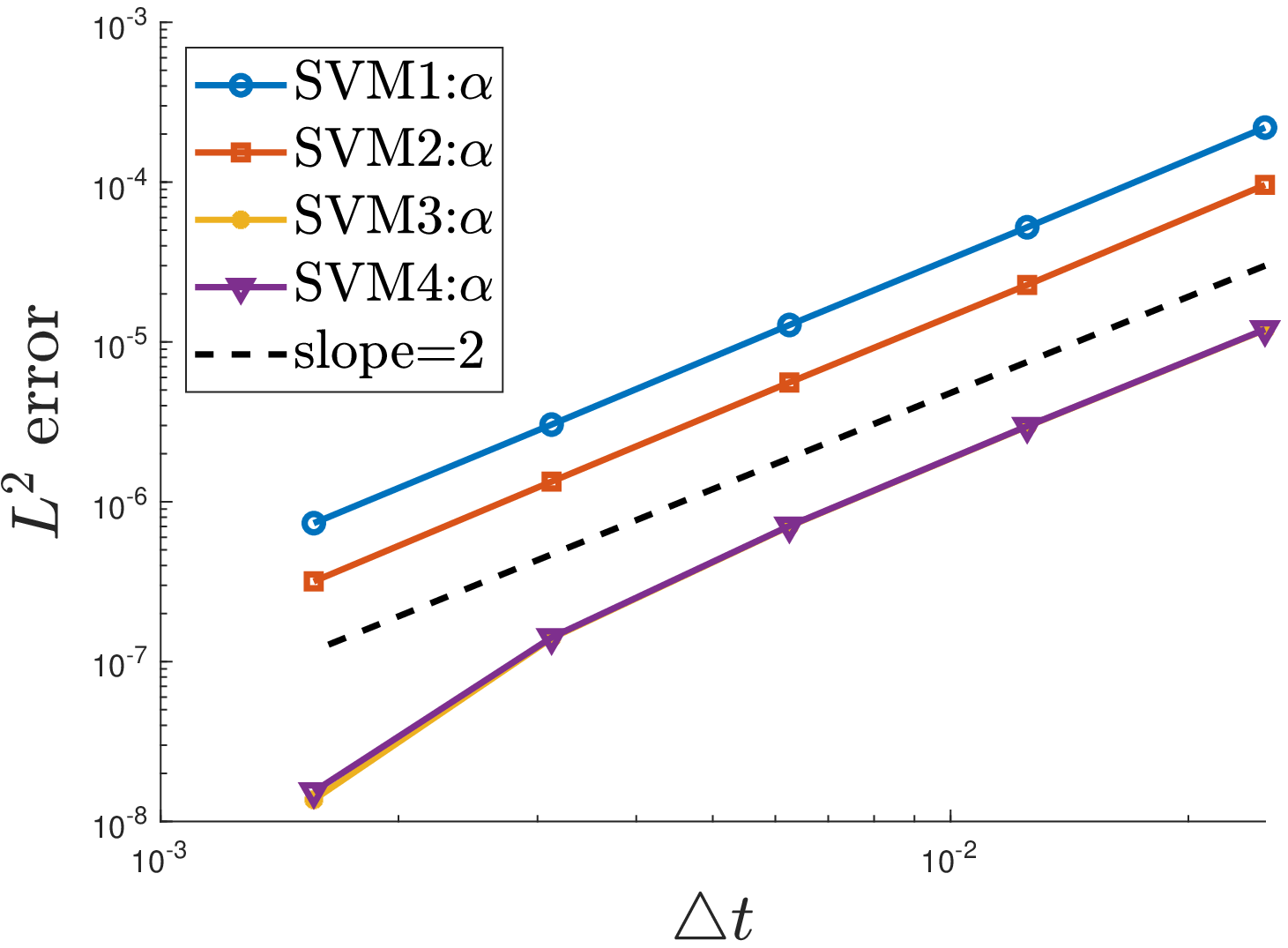}}
\subfigure[$N=128$]{\includegraphics[width=0.3\textwidth]{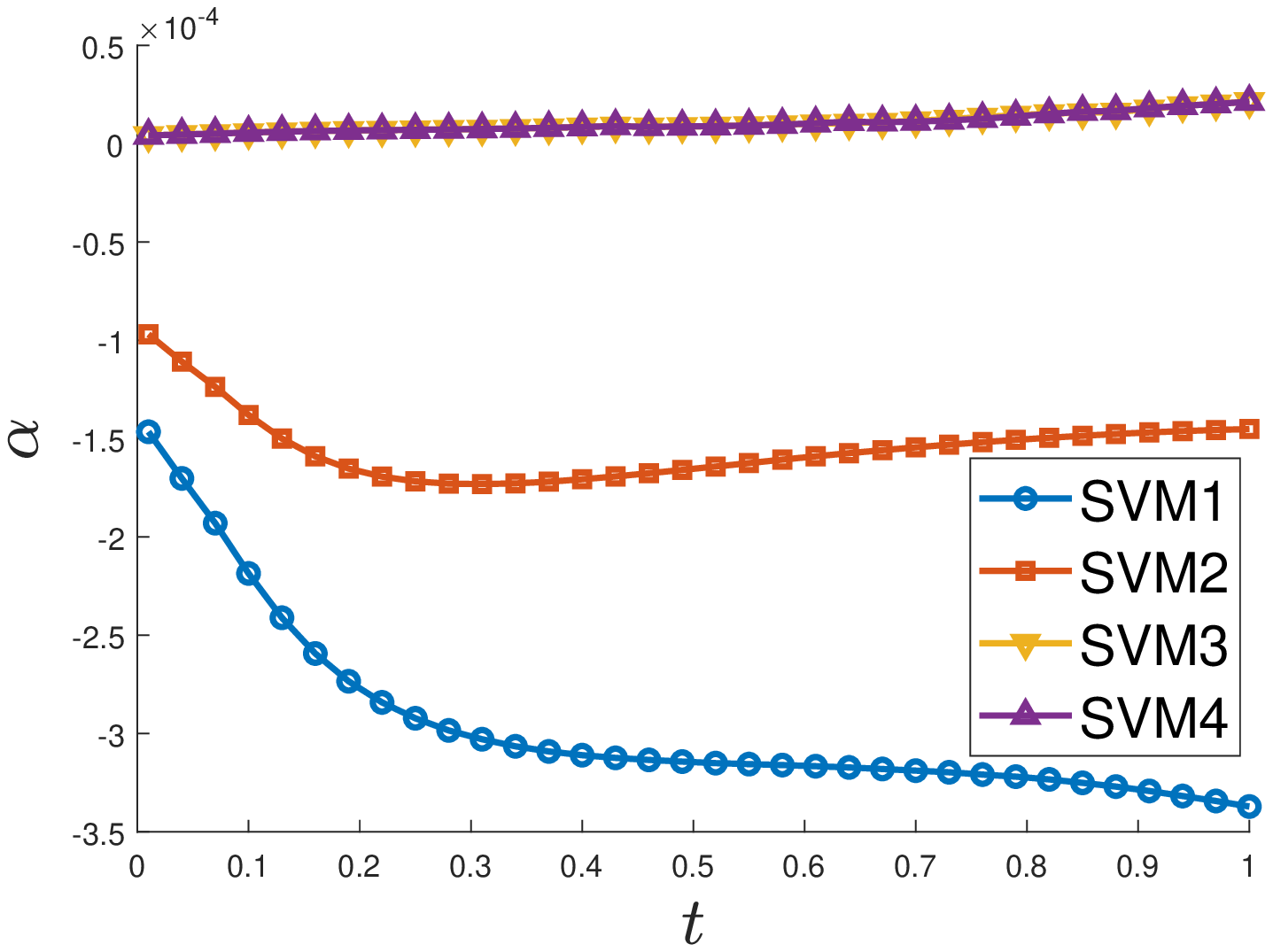}}
\subfigure[$N=256$]{\includegraphics[width=0.3\textwidth]{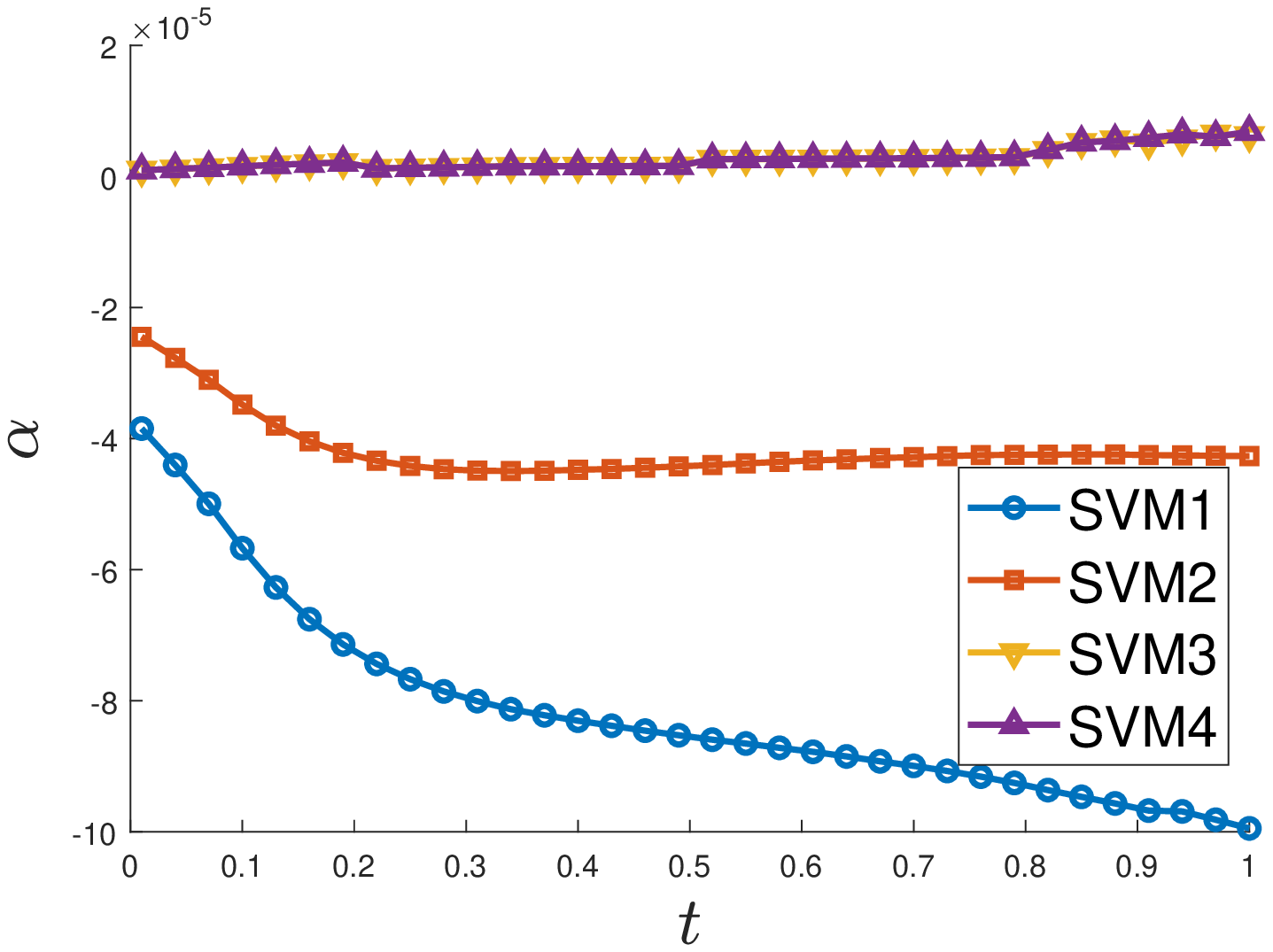}}
\\ \vspace{-0.3cm}
\caption{Supplementary variable $\alpha$ obtained in \textbf{SVM} schemes are within the expected order ($\alpha\sim\mathcal O(\triangle t^2)$). In comparison, $\alpha$ in \textbf{SVM1} and \textbf{SVM2} algorithms is larger than that in \textbf{SVM3} and \textbf{SVM4} algorithms. \label{constant_CH_accuracy_Alpha} }
\end{figure}

Taking $T=1,~\triangle t=10^{-4},~h=1/128$ and $h=1/256$ respectively, we use Example \ref{constant_CH mesh refinement} to examine the computational efficiency of the schemes. The results in CPU time(s) and convergence rate of supplementary variable $\alpha$ using each scheme are summarized in Figures \ref{constant_CH_accuracy_CPU} and \ref{constant_CH_accuracy_Alpha}, respectively. Firstly, $\alpha(t)$ scales with $C(h)\Delta t^2$ as expected. Secondly, $C(h)$ decays with $h$. Thirdly, SVM schemes compute faster than the EQ scheme on our workstation.

We use \textbf{SVM2} algorithm with $\triangle t=10^{-5}, h=1/128$ to conduct the following numerical investigations.

\subsection{Dynamics of the copolymer solution}
Notice that the volume fractions are dimensionless and we can easily nondimensionalize the model using the characteristic time and length scale. Thus, we assume the model we use in the numerical investigations are dimensionless as well. For details on the non-dimensionalization, please refer to \cite{Choksi_2005_Derive_copolymer_homopolymer}. We present 3 examples to show spot and lamellar structures in 2D copolymer solutions, which are present in copolymer melts as well. In this subsection, we fix the mobility matrix as
\begin{equation}\label{constant_CH_mobility_diagonal}
M=4\times 10^{-3}\times{\rm diag}[1,1,1] .
\end{equation}
All numerical experiments conducted by our algorithms respect energy dissipation rate and volume conservation although we do not show these in all examples.

\begin{figure}[H]
\centering
\subfigure[free energy]{\includegraphics[width=0.35\textwidth]{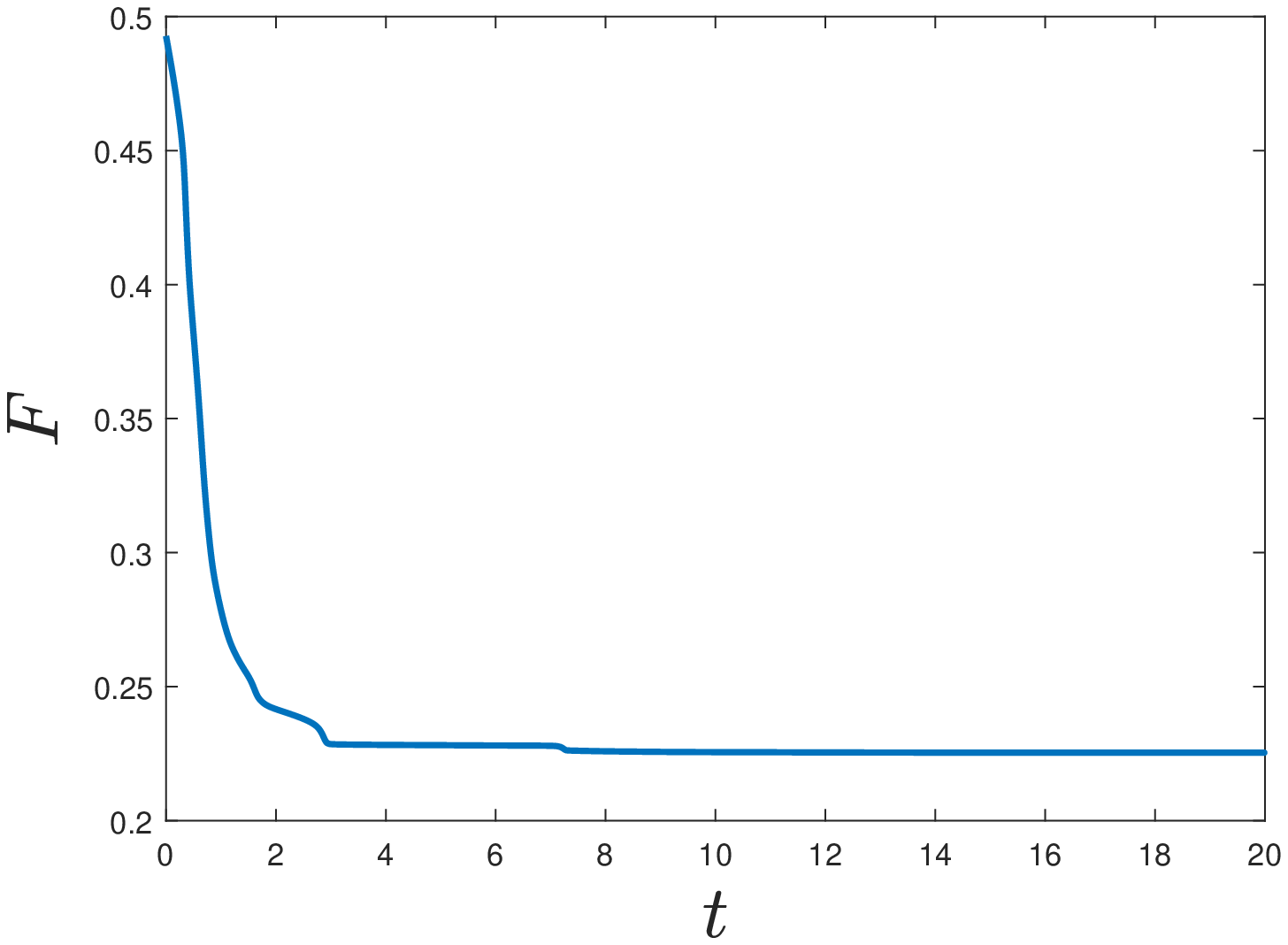}}
\qquad
\subfigure[average volume fractions]
{\includegraphics[width=0.35\textwidth]{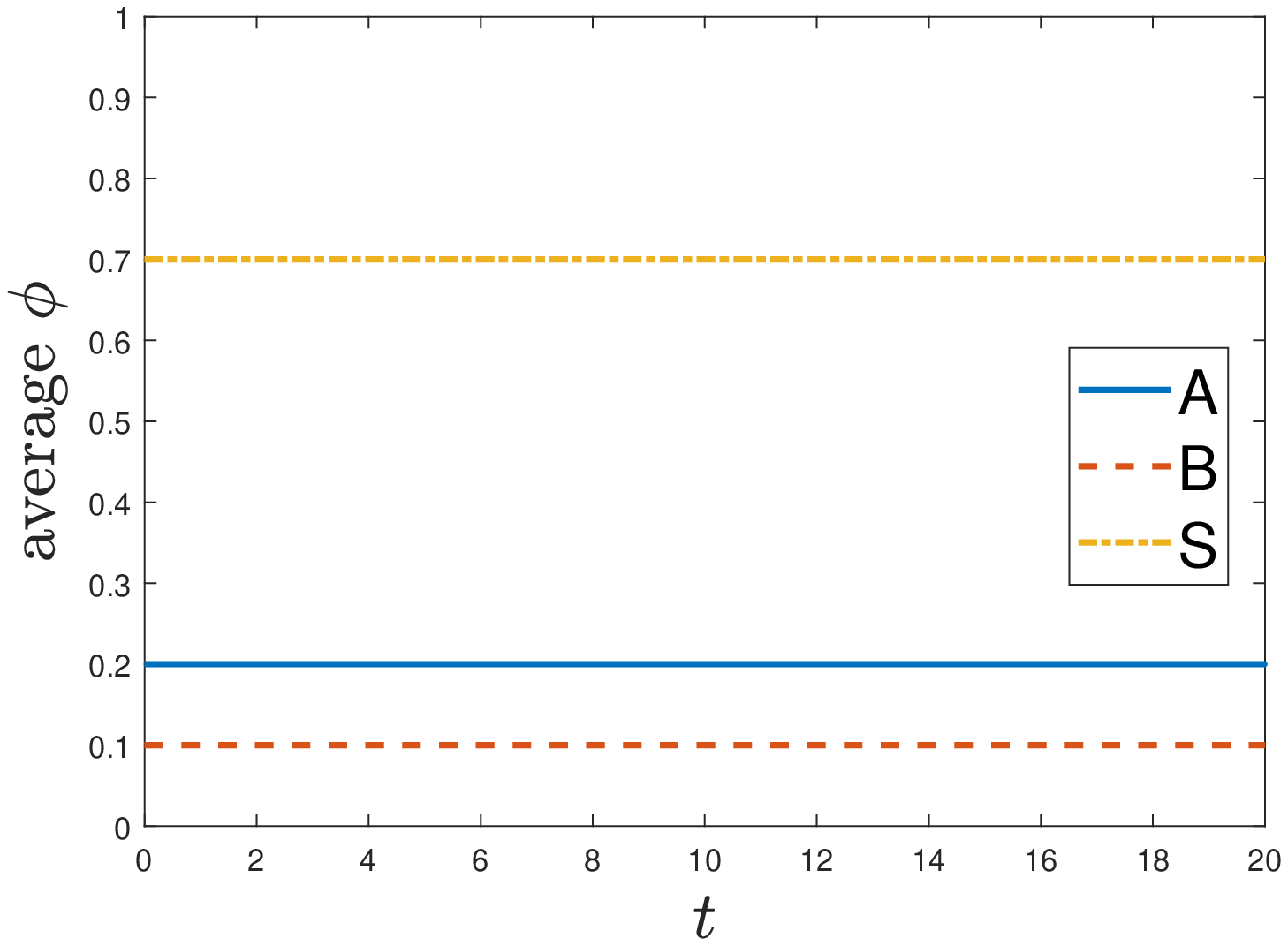}}
\\ \vspace{-0.3cm}
\caption{Energy decays and volume of each species conserves in $t\in[0,20]$, respectively. \label{constant_CH_spot_energy}}
\end{figure}

\begin{figure}[H]
\centering
\subfigure[$\phi_A-\phi_B$]{\includegraphics[width=0.7\textwidth]{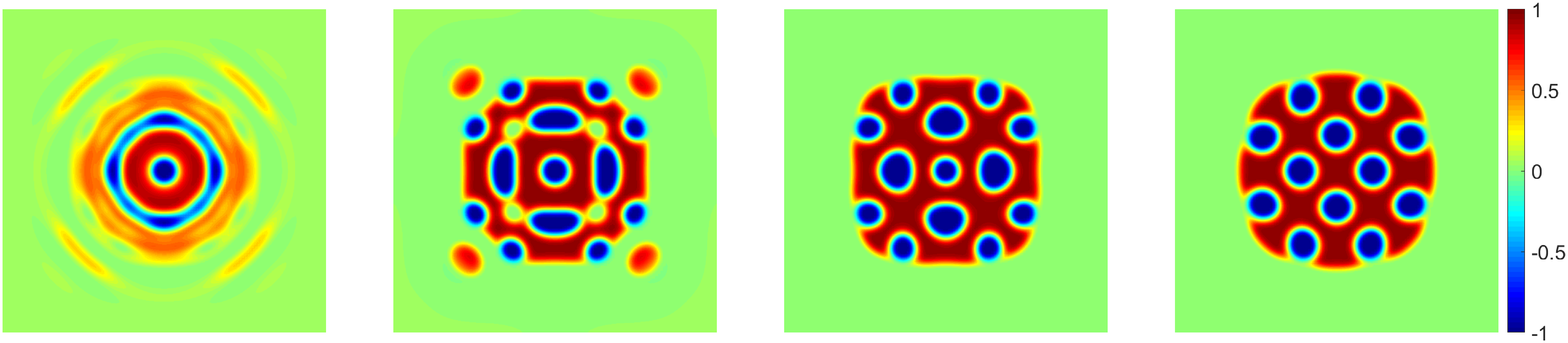}}
\\ \vspace{-0.3cm}
\subfigure[$\phi_S$]{\includegraphics[width=0.7\textwidth]{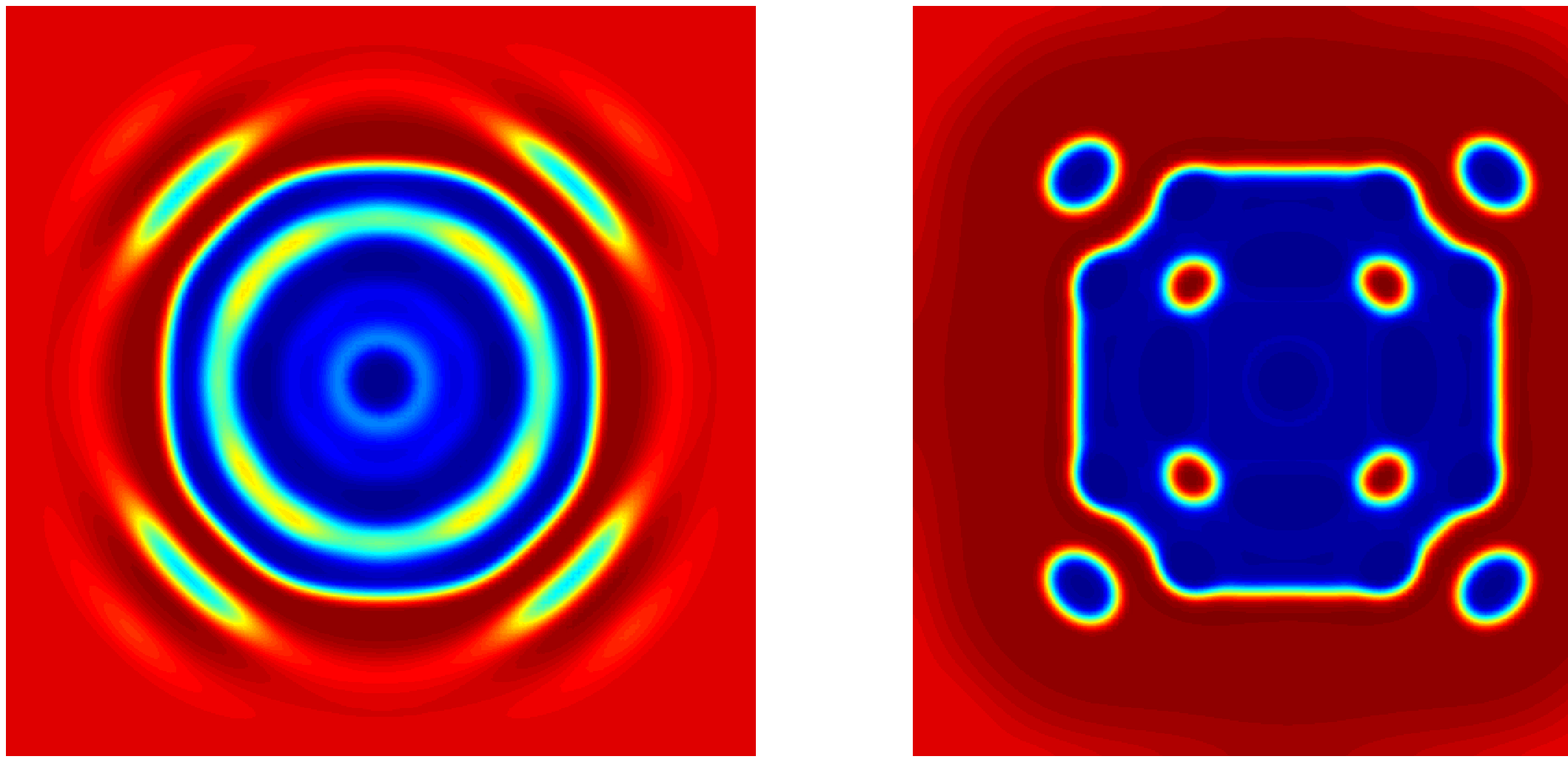}}
\\ \vspace{-0.3cm}
\caption{Microphase and macrophase separation at $t=0.5,~1,~3,~20$. At the steady state, we observe spot structures in microphase separation that $B$ forms into spots surrounded by $A$. In the 2nd row, the volume fraction of the solvent is shown. \label{constant_CH_spot_phase}}
\end{figure}

\begin{example}[\textbf{Spots}\label{constant_CH_spot}]
We set $N_A=2, ~N_B=N_S=1, ~\chi_{AA}=\chi_{BB}=\chi_{SS}=0, ~\chi_{AB}=6, ~\chi_{AS}=4, ~\chi_{BS}=8, ~\varepsilon=0.01, ~\gamma=10^3$ and $\phi^0_A=1/15+2/15(1-\cos2\pi x)(1-\cos2\pi y), ~\phi^0_B=0.5\phi^0_A, ~\phi^0_S=1-\phi^0_A-\phi^0_B$. The free energy and volume conservation are shown in Figure \ref{constant_CH_spot_energy}. Some phase portraits are plotted in Figure \ref{constant_CH_spot_phase}.
\end{example}

\begin{figure}[H]
\centering
\subfigure[$\phi_A-\phi_B$]{\includegraphics[width=0.7\textwidth]{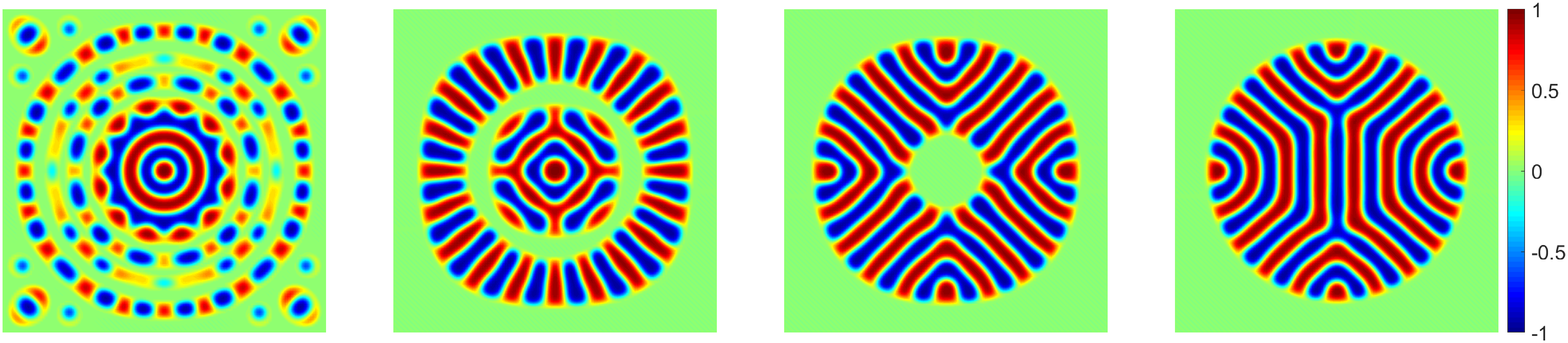}}
\\ \vspace{-0.3cm}
\subfigure[$\phi_S$]{\includegraphics[width=0.7\textwidth]{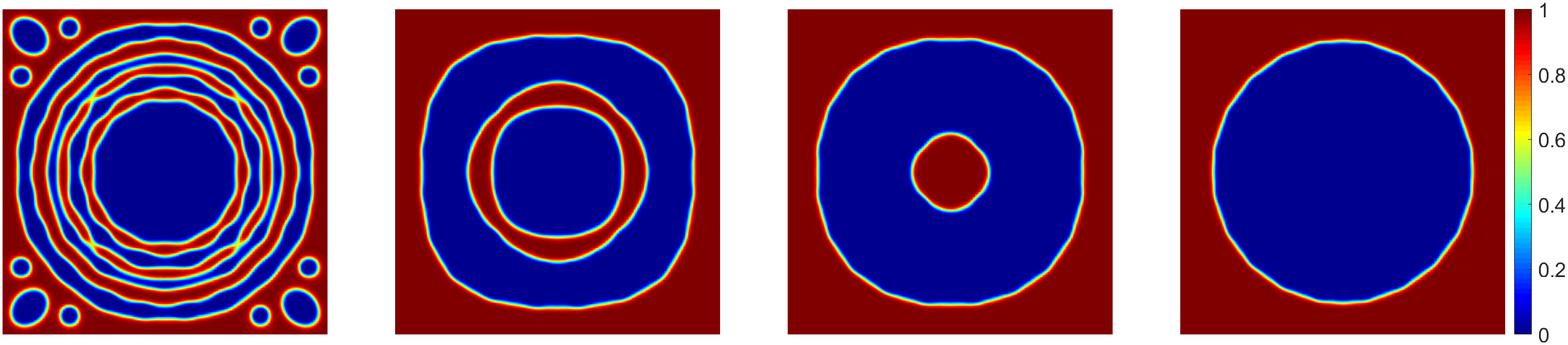}}
\\ \vspace{-0.3cm}
\caption{Microphase and macrophase separation at $t=1,~20,~40,~80$. At the steady state, we observe lamellar structures in microphase separation that both $A$ and $B$ forms into lamellae.  The solvent volume fraction is shown in the 2nd row. \label{constant_CH_lamella_phase}}
\end{figure}

\begin{example}[\textbf{Lamellae}]
We set $N_A=N_B=N_S=1, ~\chi_{AA}=\chi_{BB}=\chi_{SS}=0, ~\chi_{AB}=\chi_{AS}=6, ~\chi_{BS}=8, ~\varepsilon=0.01, ~\gamma=10^4$ and $\phi^0_A=\phi^0_B=3/16+1/16(1-\cos2\pi x)(1-\cos2\pi y), ~\phi^0_S=1-\phi^0_A-\phi^0_B$. Some phase portraits are plotted in Figure \ref{constant_CH_lamella_phase}.
\end{example}

\begin{figure}[H]
\centering
\subfigure[$\phi_A-\phi_B$]{\includegraphics[width=0.7\textwidth]{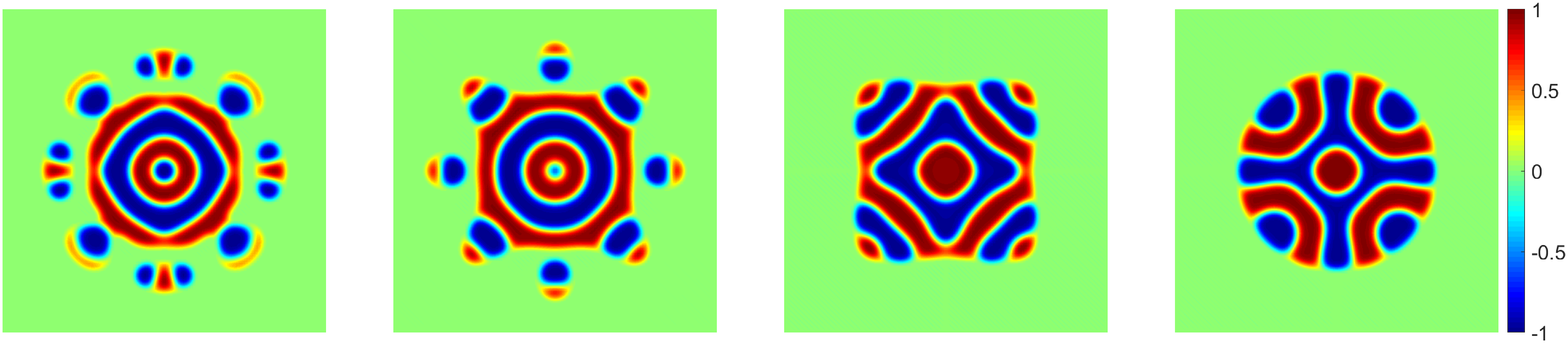}}
\\ \vspace{-0.3cm}
\subfigure[$\phi_S$]{\includegraphics[width=0.7\textwidth]{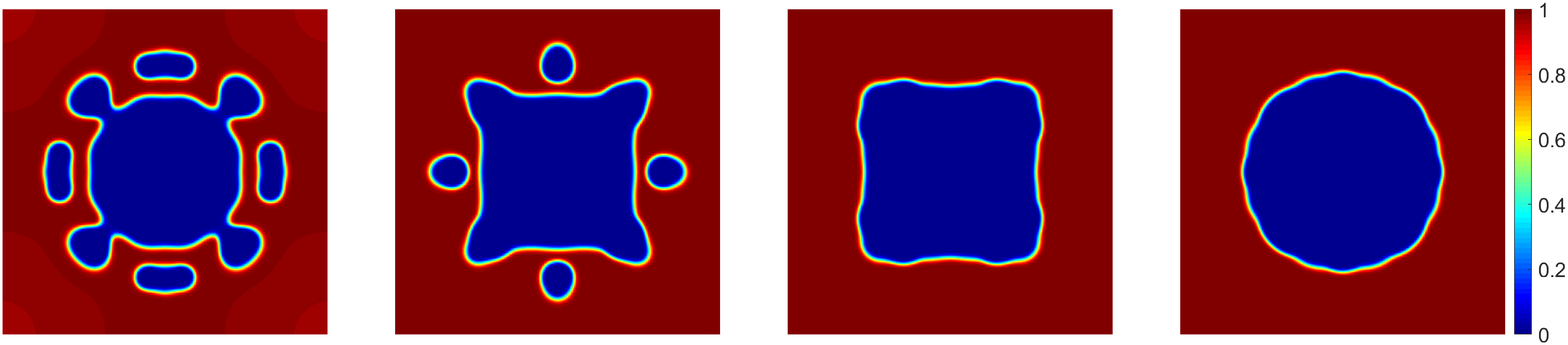}}
\\ \vspace{-0.3cm}
\caption{Microphase and macrophase separation at $t=1,~2,~5,~20$. At the steady state, we observe both lamellar and spot structures in microphase separation that there is an $A$ spot in the core surrounded by $B$ lamellae and there are 4 $B$ spots separated by $A$ lamellae scattered on the edge of the circular copolymer area. \label{constant_CH_LamellaSpot_phase}}
\end{figure}

\begin{example}[\textbf{Lamellae+spots}]
We set $N_A=N_B=N_S=1, ~\chi_{AA}=\chi_{BB}=\chi_{SS}=0, ~\chi_{AB}=\chi_{AS}=6, ~\chi_{BS}=8, ~\varepsilon=0.01, ~\gamma=10^3$ and $\phi^0_A=\phi^0_B=0.05+0.1(1-\cos2\pi x)(1-\cos2\pi y), ~\phi^0_S=1-\phi^0_A-\phi^0_B$. Some phase portraits are plotted in Figure \ref{constant_CH_LamellaSpot_phase}.
\end{example}

\subsection{Cross-coupling effect in mobility on dynamics}
In the past, nearly all numerical simulations on dynamics of copolymers focused on a mobility matrix without the cross-coupling effect, i.e., the mobility coefficient matrix is a diagonal matrix. Likewise, the corresponding friction matrix is diagonal as well. Given the fact that there exists friction between different species of polymers and long range interaction in the copolymer solution, we would like to know how the chemical potential derived from one species affects the dynamics of the other so as to impact on the overall dynamics in the mixture. This in fact investigates how forces generated from chemical potentials of different species impact on dynamics of the other species via frictions.

Notice that $M$ is the inverse of the friction matrix. So qualitatively, the smaller the molecule is, the larger the corresponding $M_{ii}$ should be. In theory, mobility matrices dictate dynamical paths of the system while relaxing back to the steady state. We all know that the steady state is determined by the free energy while the path towards the steady state in the phase space is dictated by the mobility. In this subsection, we choose mobility matrices with diagonally dominant elements and the cross-coupling effect is represented by the off-diagonal entries.

\begin{figure}[H]
\centering
\subfigure[$M1$]{\includegraphics[width=0.7\textwidth]{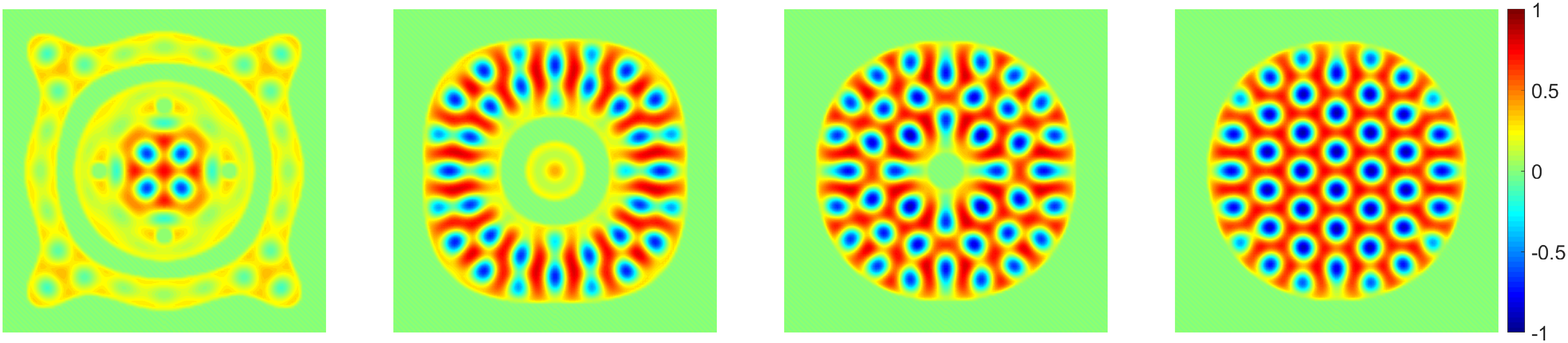}}
\\ \vspace{-0.3cm}
\subfigure[$M2$]{\includegraphics[width=0.7\textwidth]{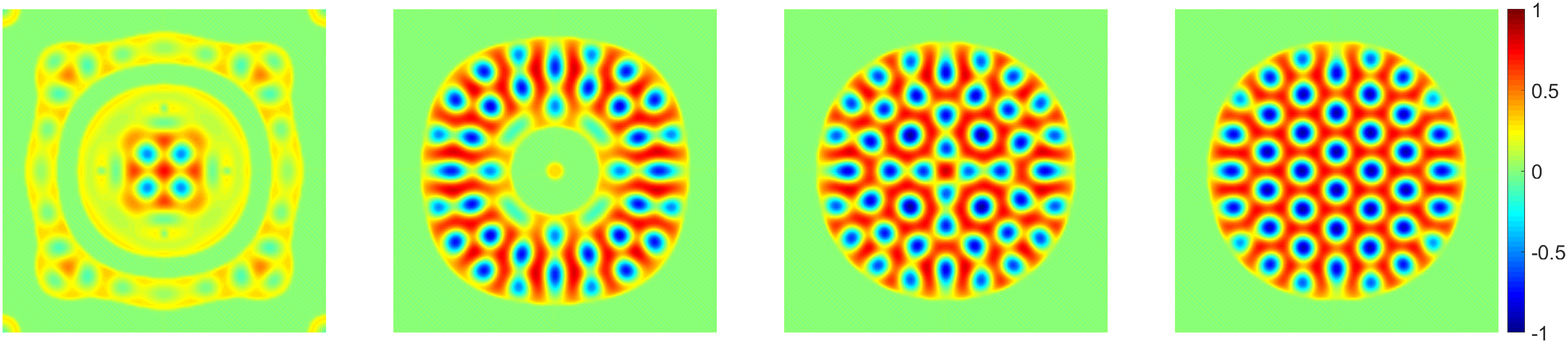}}
\\ \vspace{-0.3cm}
\subfigure[$M3$]{\includegraphics[width=0.7\textwidth]{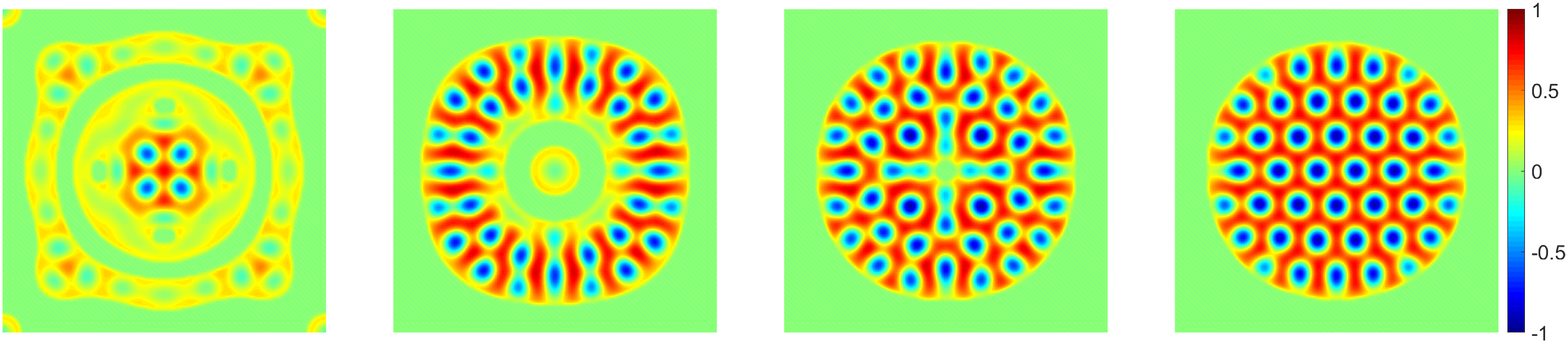}}
\\ \vspace{-0.3cm}
\subfigure[$M4$]{\includegraphics[width=0.7\textwidth]{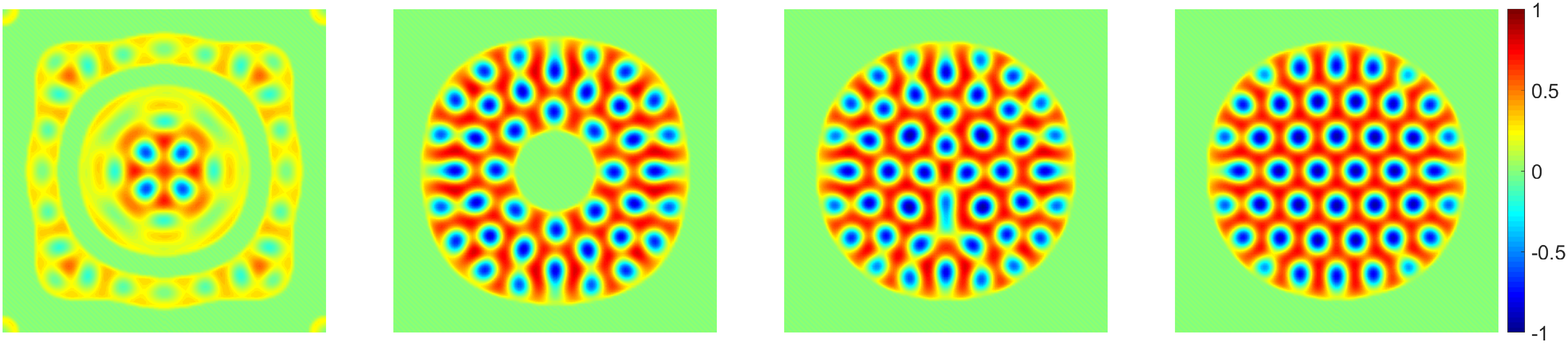}}
\\ \vspace{-0.3cm}
\subfigure[$M5$]{\includegraphics[width=0.7\textwidth]{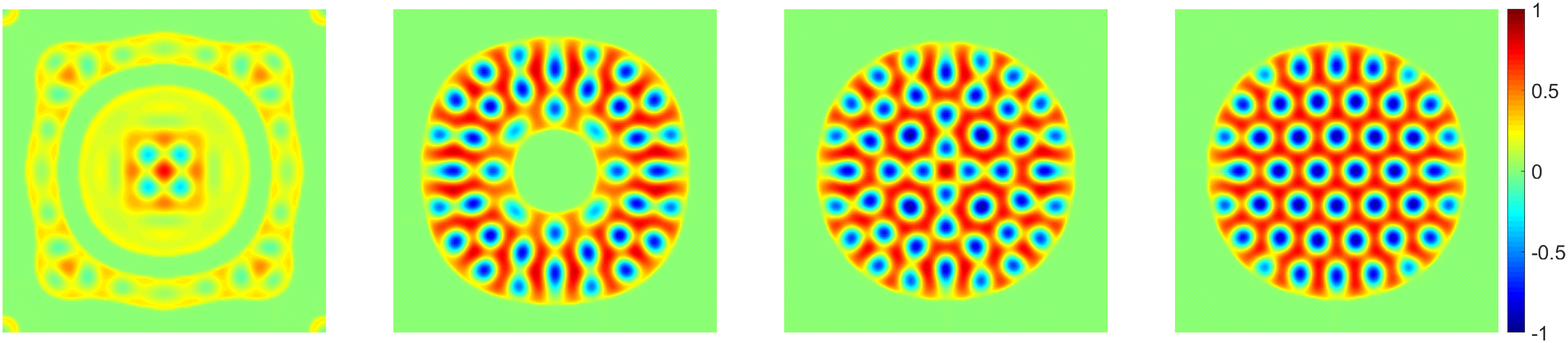}}
\\ \vspace{-0.3cm}
\subfigure[$M6$]{\includegraphics[width=0.7\textwidth]{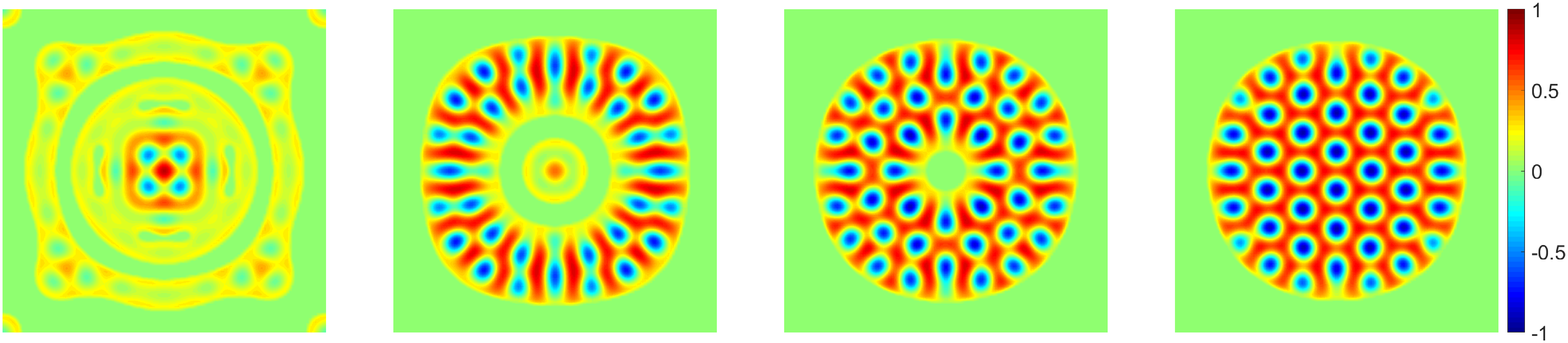}}
\\ \vspace{-0.3cm}
\caption{$\phi_A-\phi_B$ at $t=4,~25,~45,~120$ with respect to 6 different mobilites. From the core parts of the copolymer area, we observe different microscopic phase separations in transient dynamics decided by 6 different mobilities (It is helpful to observe macrophase separations simultaneously). At $t=120$, all 6 steady state structures are similar with (a)(b)(f) are 90 degree rotation of (c)(d)(e). Since there is no external field is present, the steady states are essentially the same modulo a rotational motion. \label{constant_CH_mobility_AB} }
\end{figure}

\begin{figure}[H]
\centering
\subfigure[$M1$]{\includegraphics[width=0.7\textwidth]{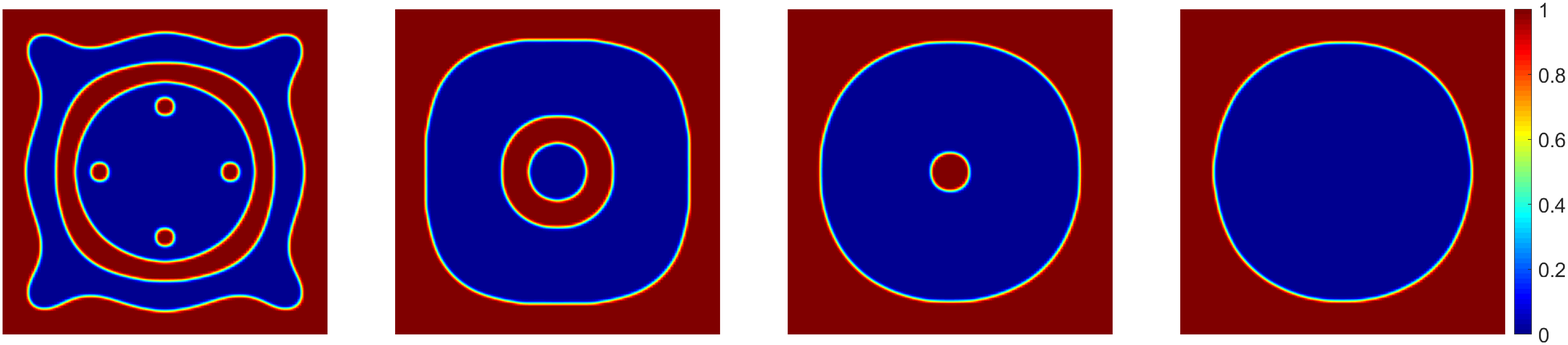}}
\\ \vspace{-0.3cm}
\subfigure[$M2$]{\includegraphics[width=0.7\textwidth]{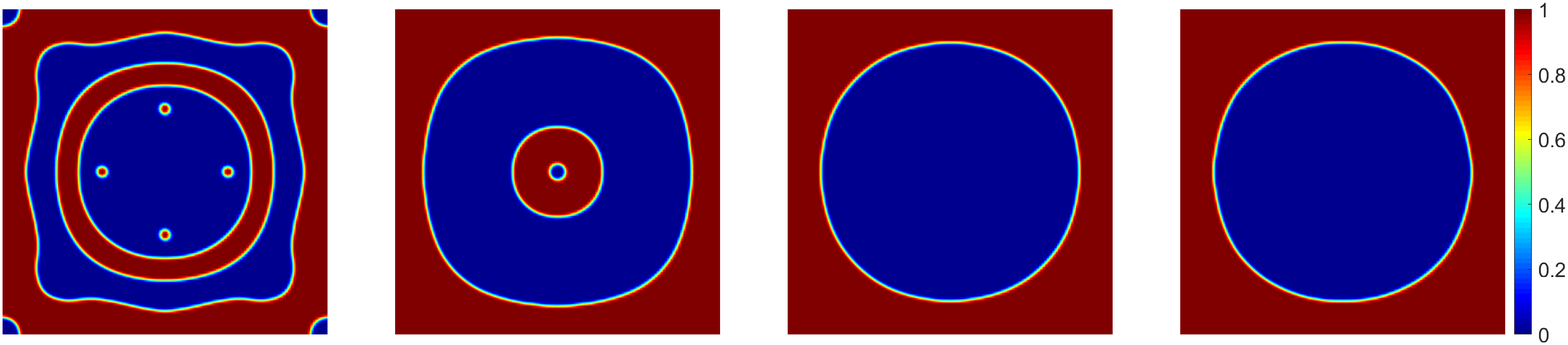}}
\\ \vspace{-0.3cm}
\subfigure[$M3$]{\includegraphics[width=0.7\textwidth]{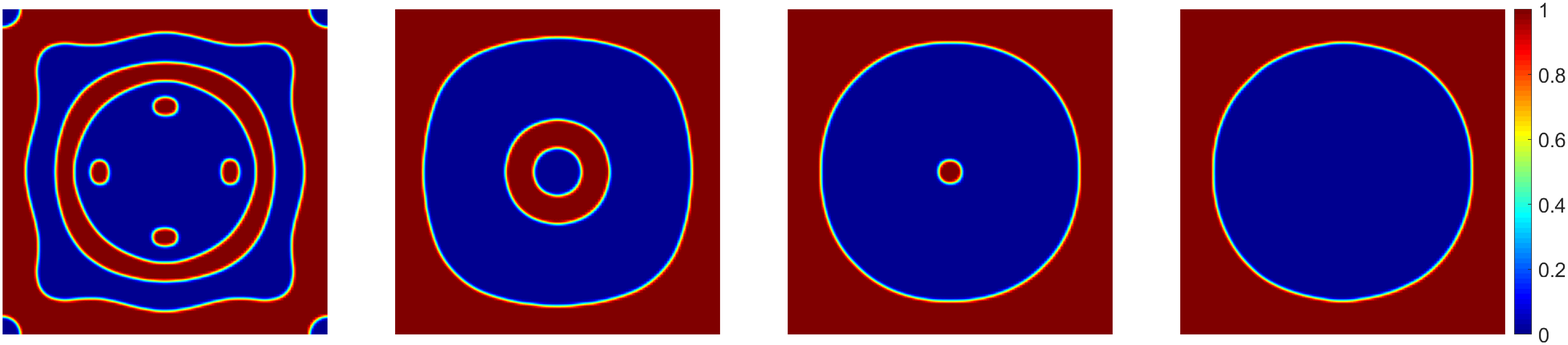}}
\\ \vspace{-0.3cm}
\subfigure[$M4$]{\includegraphics[width=0.7\textwidth]{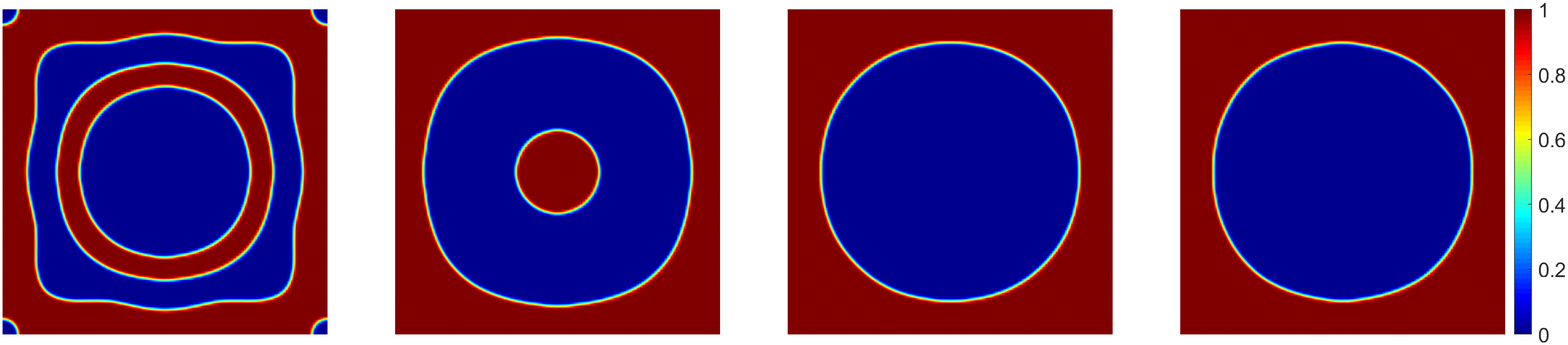}}
\\ \vspace{-0.3cm}
\subfigure[$M5$]{\includegraphics[width=0.7\textwidth]{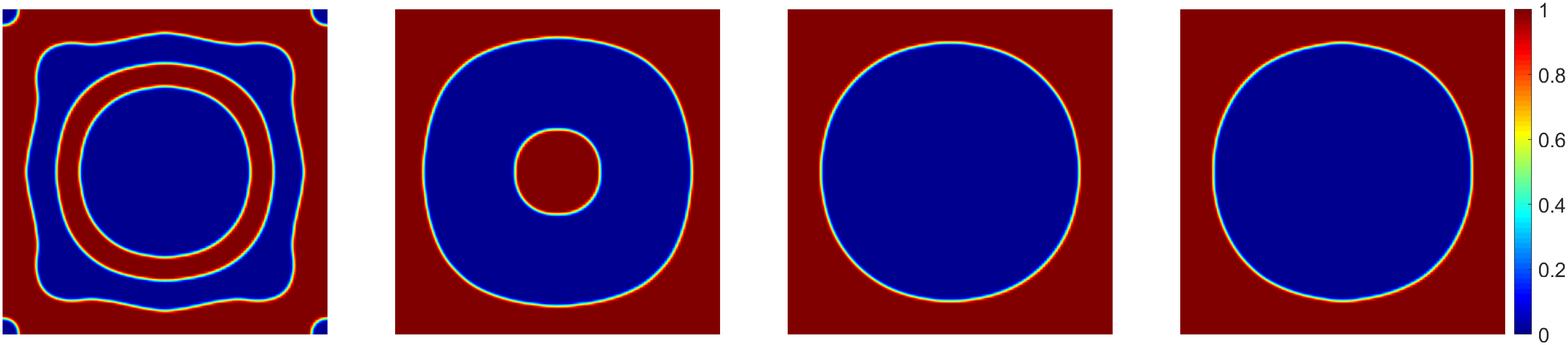}}
\\ \vspace{-0.3cm}
\subfigure[$M6$]{\includegraphics[width=0.7\textwidth]{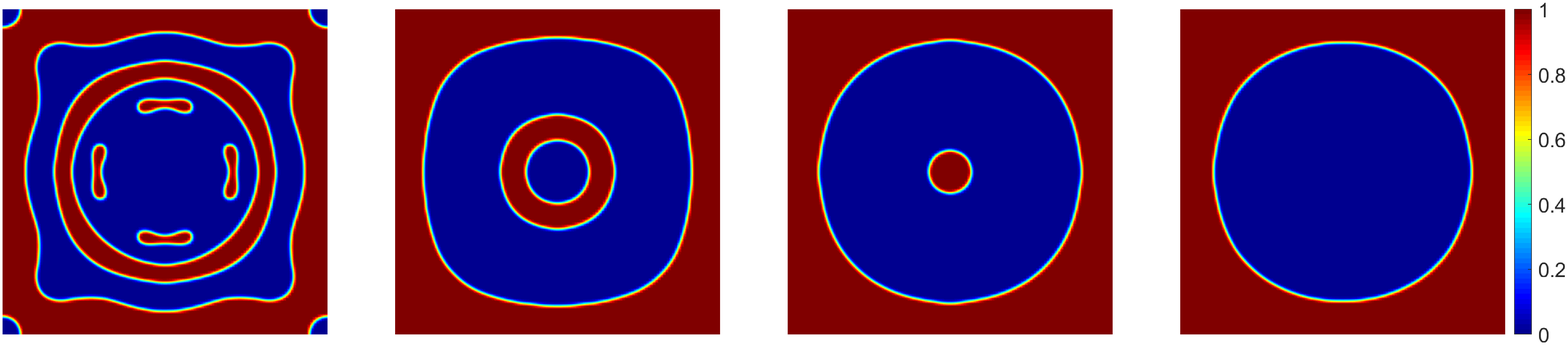}}
\\ \vspace{-0.3cm}
\caption{$\phi_S$ at $t=4,~25,~45,~120$ with respect to 6 different mobilites. Macroscopic phase separations are different in transient dynamics deictated by 6 different mobilities. At $t=120$, they are all the same. \label{constant_CH_mobility_S} }
\end{figure}

\begin{figure}[H]
\centering
\includegraphics[width=0.7\textwidth]{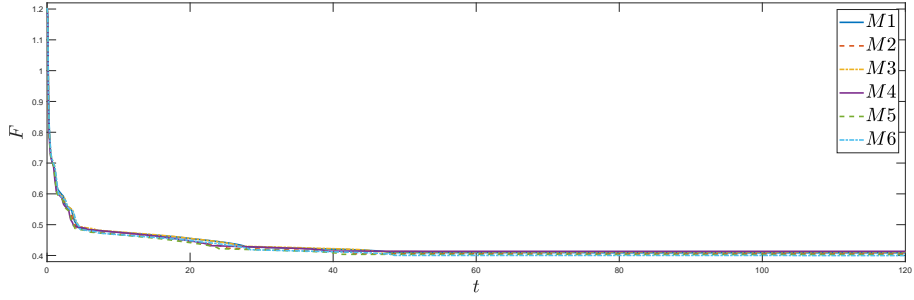}
\\ \vspace{-0.3cm}
\caption{Free energies in $t\in[0,120]$. There are slight differences among the results quantitatively. \label{constant_CH_mobility_energy} }
\end{figure}

\begin{example}[\textbf{Cross-coupling effect in mobility}\label{constant_CH_mobility}]
We set $N_A=3, ~N_B=2, ~N_S=1, ~\chi_{AA}=\chi_{BB}=\chi_{SS}=0, ~\chi_{AB}=4, ~\chi_{AS}=6, ~\chi_{BS}=8, ~\varepsilon=0.01, ~\gamma=10^4$ and $\phi^0_A=3/14+3/35(1-\cos2\pi x)(1-\cos2\pi y), ~\phi^0_B=2/3\phi^0_A, ~\phi^0_S=1-\phi^0_A-\phi^0_B$. We take 6 different mobility matrices respectively given by
\begin{equation}
\begin{gathered}
M1: 10^{-3} \!\times\!
\left[\begin{matrix} 4 & & \\ & 4 & \\ & & 4 \end{matrix}\right] ,
~
M2: 10^{-3} \!\times\!
\left[\begin{matrix} 3 & & \\ & 4 & \\ & & 5 \end{matrix}\right] ,
~
M3: 10^{-3} \!\times\!
\left[\begin{matrix} 3 & -0.5 & \\ -0.5 & 4 & \\  & & 5 \end{matrix}\right] ,
\\
M4: 10^{-3} \!\times\!
\left[\begin{matrix} 3 & -0.5 & -0.5 \\ -0.5 & 4 & -0.5 \\ -0.5 & -0.5 & 5 \end{matrix}\right] ,
~
M5: 10^{-3} \!\times\!
\left[\begin{matrix} 3 & 0.5 & \\ 0.5 & 4 & \\ & & 5 \end{matrix}\right] ,
~
M6: 10^{-3} \!\times\!
\left[\begin{matrix} 3 & 0.5 & 0.5 \\ 0.5 & 4 & 0.5 \\ 0.5 & 0.5 & 5 \end{matrix}\right] .
\end{gathered}
\end{equation}
Some phase portraits and free energies are plotted in Figures \ref{constant_CH_mobility_AB}, \ref{constant_CH_mobility_S} and \ref{constant_CH_mobility_energy}, respectively.
\end{example}

Figure \ref{constant_CH_mobility_AB} shows the microscopic phase separation between $A$ and $B$ polymers, while Figure \ref{constant_CH_mobility_S} shows the macroscopic phase separation between the copolymer and the solvent. At $t=120$, the system reaches the metastable/steady state, which is determined by the free energy. In this example, all mobilities eventually lead to the same  metastable/steady state as expected since mobility is instrumental to transient dynamics not to the final steady state. During transient dynamics, the free energy decays with time as shown in Figure \ref{constant_CH_mobility_energy}.

We reiterate that different mobilities result in different dynamics, while different mobilities do not necessarily lead to different metastable/steady state. In fact, if the steady state is globally stable and unique, it is independent of the mobility. In Example \ref{constant_ECH_mobility} and Example \ref{constant_MCH_mobility}, 3 different mobilities also lead to the same steady state. However, different mobilities lead to different metastable states as depicted in Figure \ref{constant_CH_mobility_AB}(b)(c), dictated by mobilities $M2,M3$. For discussion about metastable states and energy landscapes, we suggest readers refer to \cite{Choksi_2015_Metastable_states, Choksi_2011_2D_copolymer_melt, ZhangLei_2020_Energy_landscape, ZhangLei_2020_Find_saddle_points}.

\subsection{Dynamics of the copolymer solution coupled with the electric field}
When an electric field is applied to the copolymer solution, patterns normally emerge in long time. During the transient state from an irregular phase morphology to lamellar patterns, branches are formed momentarily and annihilated eventually by the electric field. This is a process of defect removal in material processing, where one refers to the branches as defects. In the following example, we would like to showcase the process at selected time slots.

\begin{figure}[H]
\centering
\subfigure[mobility 1]{\includegraphics[width=0.7\textwidth]{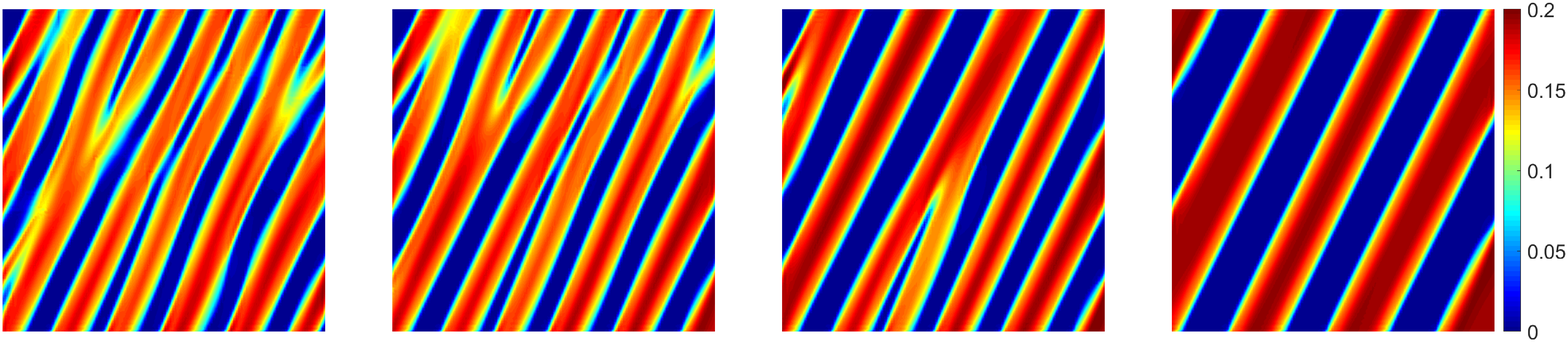}}
\\ \vspace{-0.3cm}
\subfigure[mobility 2]{\includegraphics[width=0.7\textwidth]{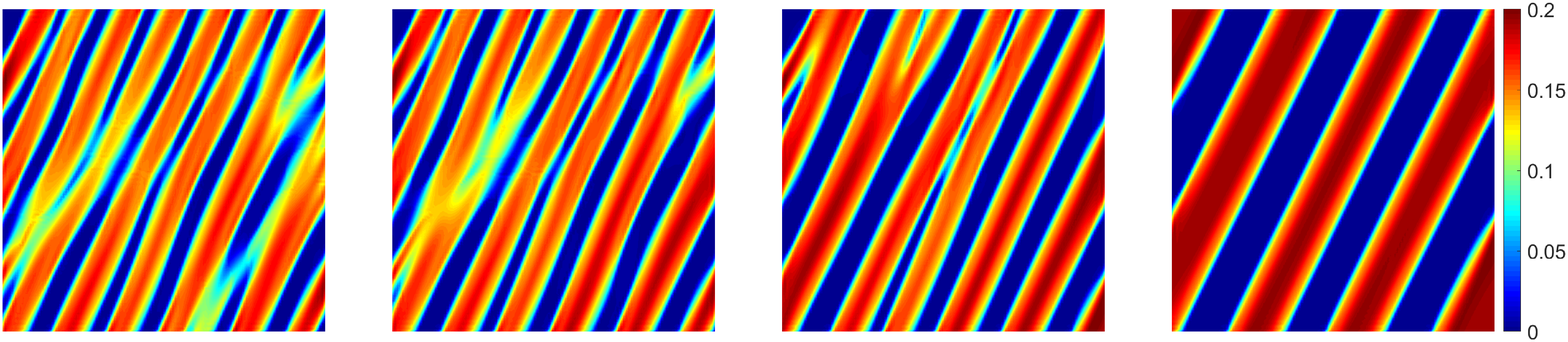}}
\\ \vspace{-0.3cm}
\subfigure[mobility 3]{\includegraphics[width=0.7\textwidth]{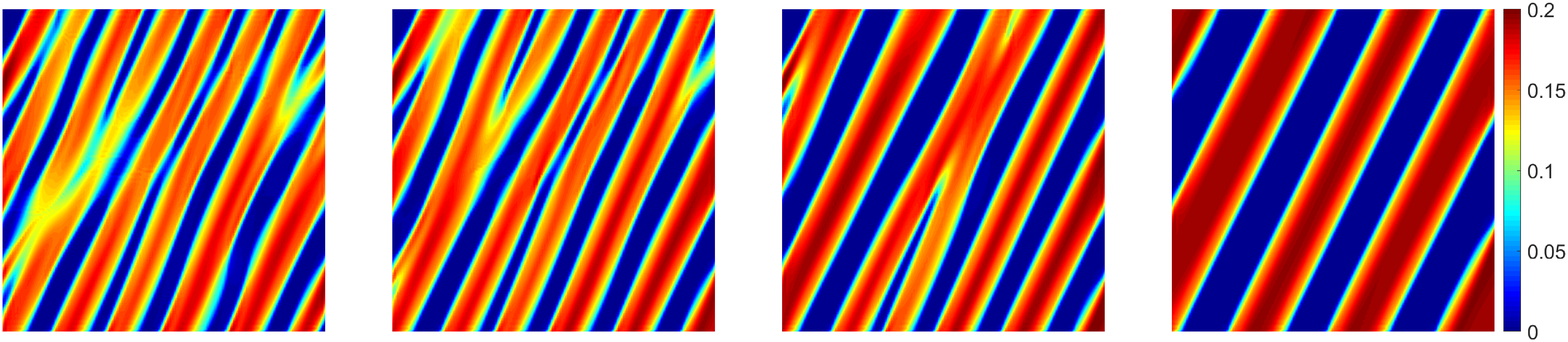}}
\\ \vspace{-0.3cm}
\caption{$\phi_A-\phi_B$ at $t=12,~15,~22,~200$ in dynamics with respect to 3 different mobilities. Defects in the microscopic phase separation are observed in transient. The electric field irons out the defect branches. The system generates a lamellar pattern in the steady state. \label{constant_ECH_mobility_AB} }
\end{figure}

\begin{figure}[H]
\centering
\subfigure[mobility 1]{\includegraphics[width=0.7\textwidth]{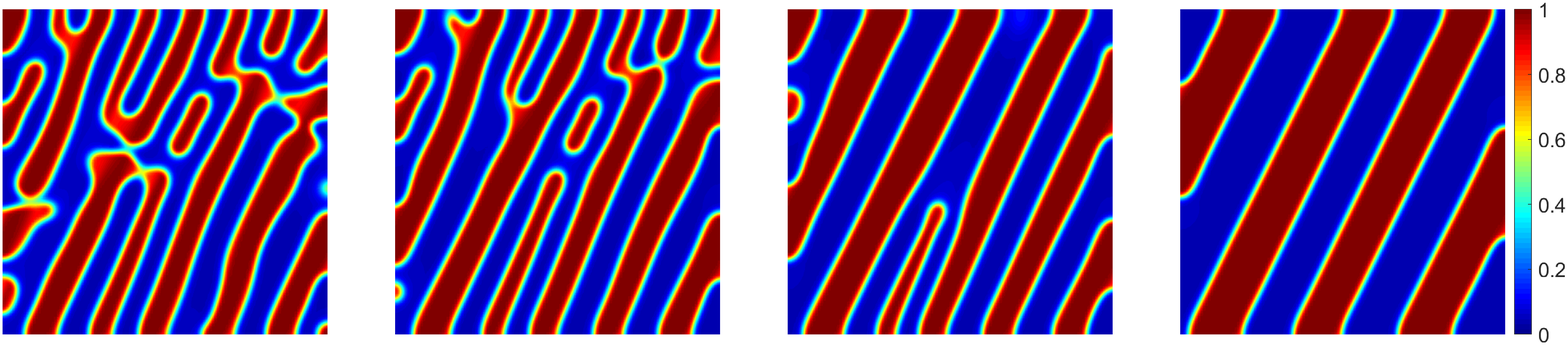}}
\\ \vspace{-0.3cm}
\subfigure[mobility 2]{\includegraphics[width=0.7\textwidth]{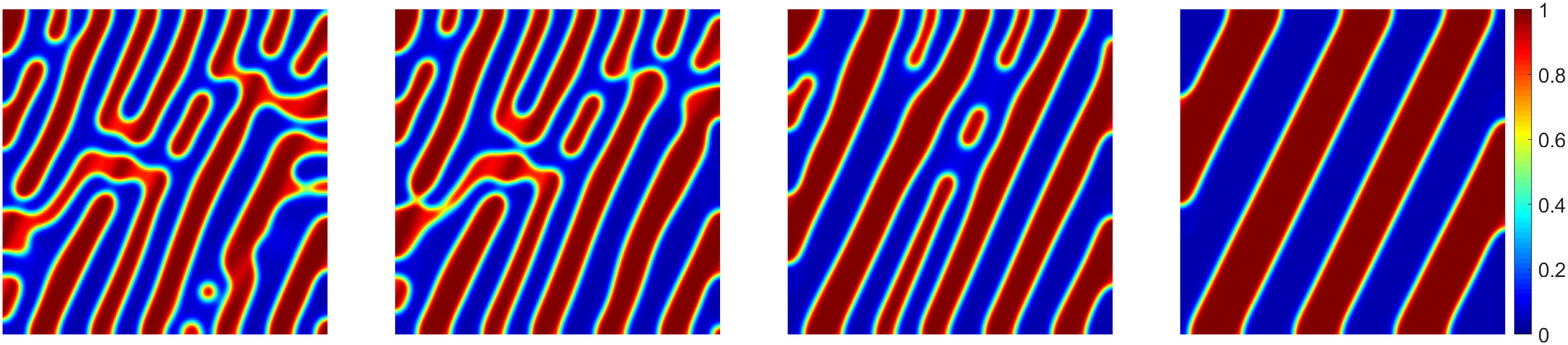}}
\\ \vspace{-0.3cm}
\subfigure[mobility 3]{\includegraphics[width=0.7\textwidth]{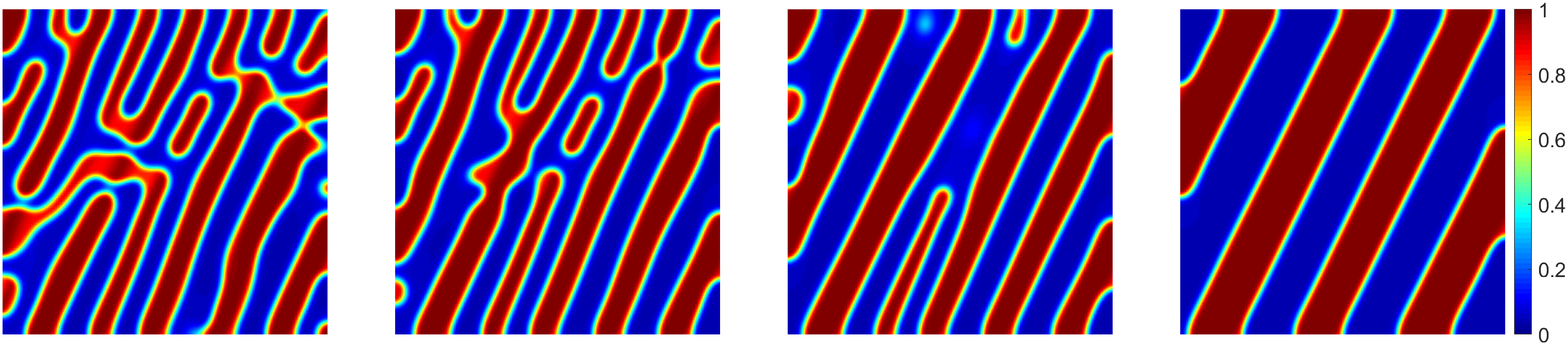}}
\\ \vspace{-0.3cm}
\caption{$\phi_S$ at $t=12,~15,~22,~200$ in dynamics with respect to 3 different mobilities. Defects in the macroscopic phase separation are observed in transient. Final lamellar patterned phase portraits are obtained through the Ostwald ripening. \label{constant_ECH_mobility_S} }
\end{figure}

\begin{figure}[H]
\centering
\includegraphics[width=0.7\textwidth]{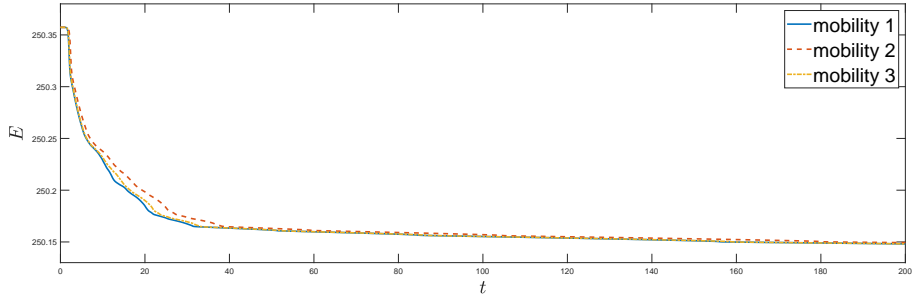}
\\ \vspace{-0.3cm}
\caption{Energies in $t\in[0,200]$ with respect to 3 different mobilities. \label{constant_ECH_mobility_energy} }
\end{figure}

\begin{example}[\textbf{Pattern formation driven by the electric field}\label{constant_ECH_mobility}]
We set $N_A=15, ~N_B=10, ~N_S=1, ~\chi_{AA}=\chi_{BB}=\chi_{SS}=0, ~\chi_{AB}=1, ~\chi_{AS}=2, ~\chi_{BS}=4, ~\varepsilon=0.01, ~\gamma=10^5, ~\epsilon_0=\epsilon_1=1, ~\bm E_0=[10,20]^\top$ and \begin{equation}\label{constant_CH_mobility_matrices}
{\rm mobility}~1 \!:\! 10^{-4} \!\times\!
\left[\begin{matrix} 4 & & \\ & 6 & \\ & & 20 \end{matrix}\right] \!,
~
{\rm mobility}~2 \!:\! 10^{-4} \!\times\!
\left[\begin{matrix} 4 & -1 & \\ -1 & 6 & \\ & & 20 \end{matrix}\right] \!,
~
{\rm mobility}~3 \!:\! 10^{-4} \!\times\!
\left[\begin{matrix} 4 & -1 & -1 \\ -1 & 6 & -1 \\ -1 & -1 & 20 \end{matrix}\right] \!.
\end{equation}
The initial conditions are
\begin{align}\label{constant_CH_rand_initial}
\phi^0_A = 0.3 + 0.001\times{\rm rand}(x,y) ,
\quad
\phi^0_B = 0.2 + 0.001\times{\rm rand}(x,y) ,
\quad
\phi^0_S = 1 - \phi^0_A - \phi^0_B ,
\end{align}
where ${\rm rand}(x,y)\in(-1,1)$ denotes a random number generator that follows the uniform distribution. Some phase portraits and the total energy are plotted in Figures \ref{constant_ECH_mobility_AB}, \ref{constant_ECH_mobility_S} and \ref{constant_ECH_mobility_energy}. Notice that the imposed electric field is slanted in this example.
\end{example}
Under the electric field, lamellar stripes emerge and so do defects. Defects move along the direction of the electric field and disappear at the boundary of the computational domain. The number of defects decreases to $0$ eventually and the morphology of the final meso-structures maintain a lamellar pattern with stripes along the direction of the electric field. The free energy curves decay with time as shown in Figure \ref{constant_ECH_mobility_energy}.

\subsection{Dynamics of the copolymer solution coupled with the magnetic field}
Like the electric field, the magnetic field aligns the copolymer in the solution and eliminate transient defects as well.

\begin{figure}[H]
\centering
\subfigure[mobility 1]{\includegraphics[width=0.7\textwidth]{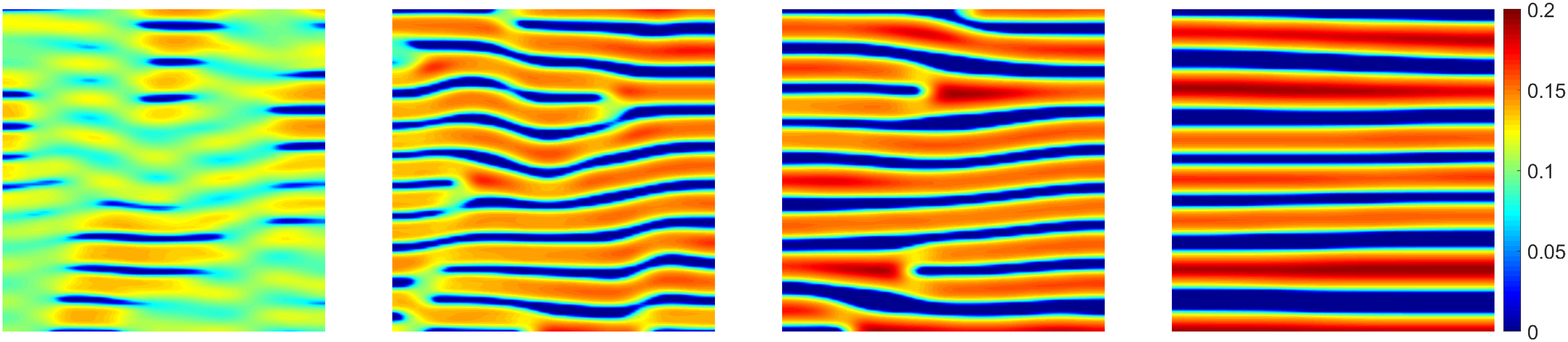}}
\\ \vspace{-0.3cm}
\subfigure[mobility 2]{\includegraphics[width=0.7\textwidth]{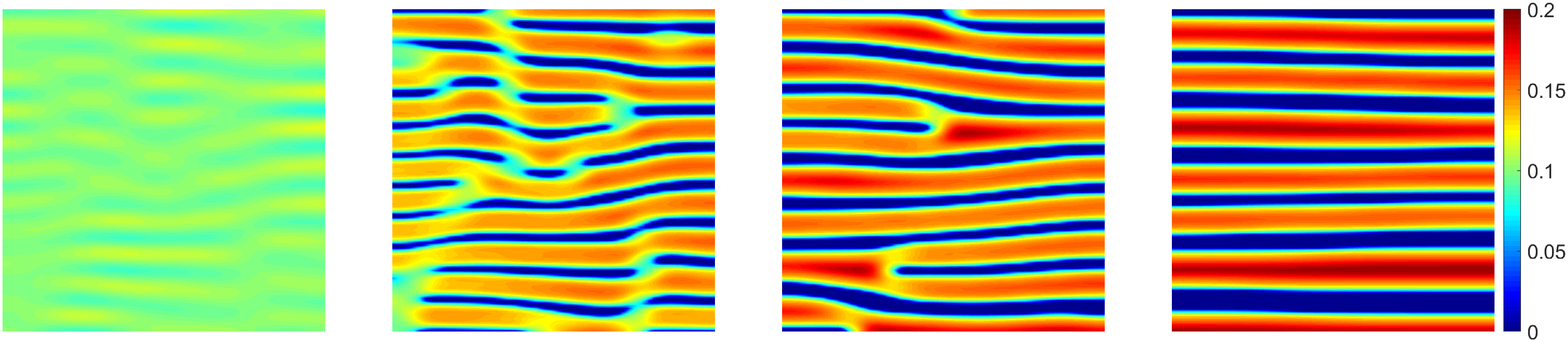}}
\\ \vspace{-0.3cm}
\subfigure[mobility 3]{\includegraphics[width=0.7\textwidth]{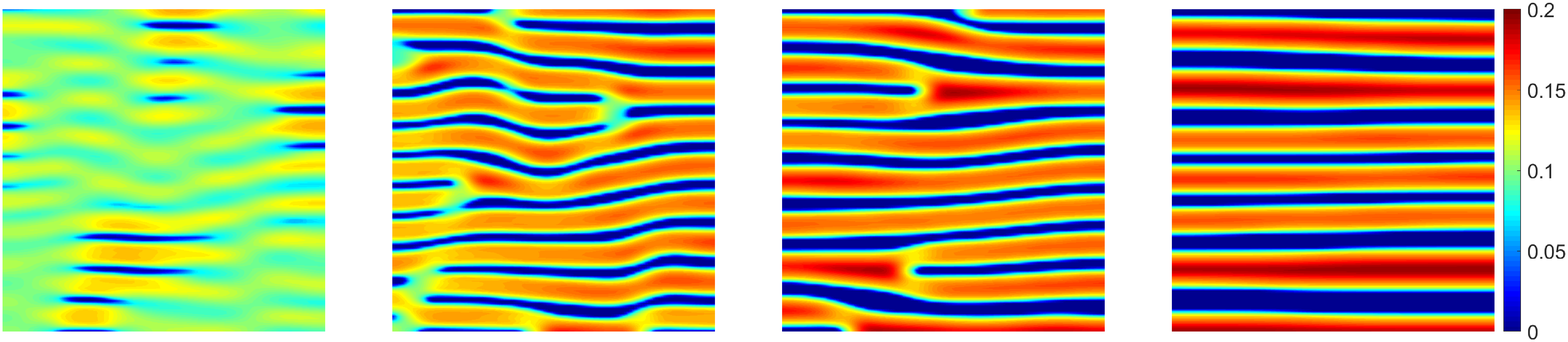}}
\\ \vspace{-0.3cm}
\caption{$\phi_A-\phi_B$ at $t=2.5,~3.5,~10,~100$. Microscopic phase separation is shown for 3 different mobilities. \label{constant_MCH_mobility_AB} }
\end{figure}

\begin{figure}[H]
\centering
\subfigure[mobility 1]{\includegraphics[width=0.7\textwidth]{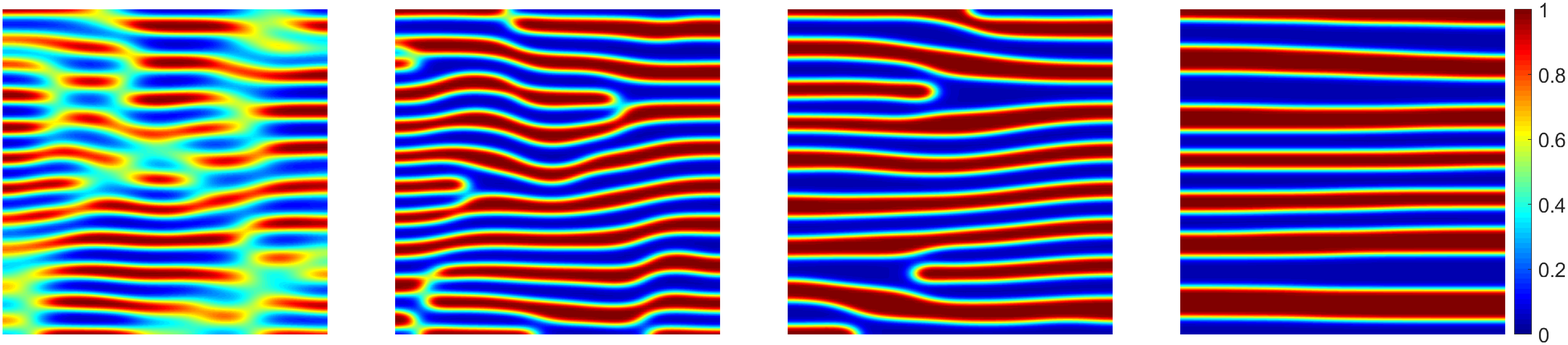}}
\\ \vspace{-0.3cm}
\subfigure[mobility 2]{\includegraphics[width=0.7\textwidth]{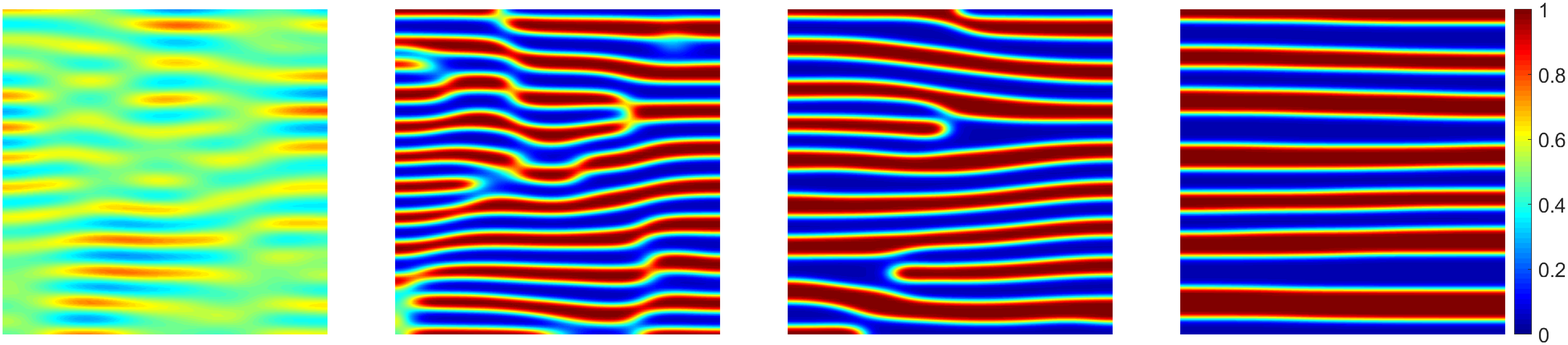}}
\\ \vspace{-0.3cm}
\subfigure[mobility 3]{\includegraphics[width=0.7\textwidth]{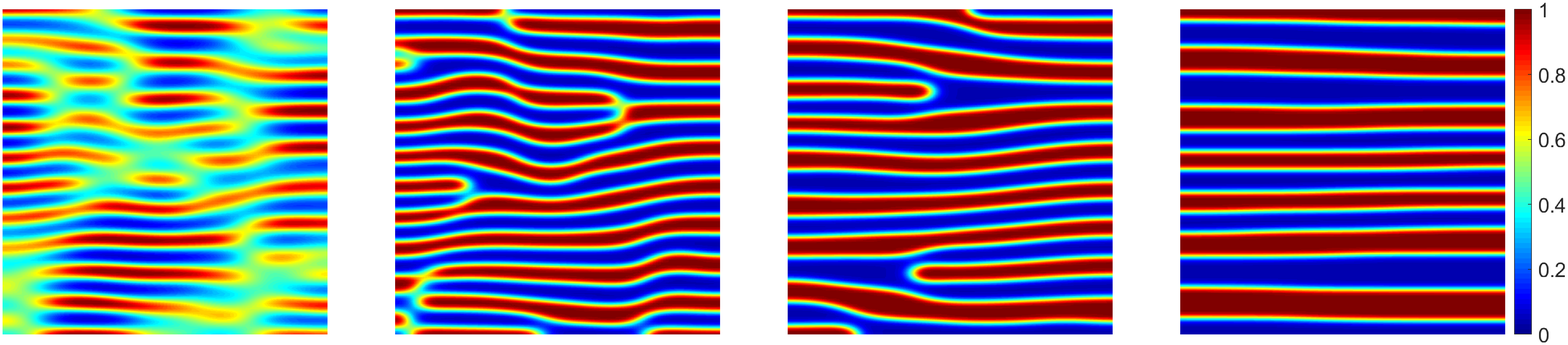}}
\\ \vspace{-0.3cm}
\caption{$\phi_S$ at $t=2.5,~3.5,~10,~100$. Macroscopic phase separation is shown for 3 different mobilities. \label{constant_MCH_mobility_S} }
\end{figure}

\begin{figure}[H]
\centering
\includegraphics[width=0.7\textwidth]{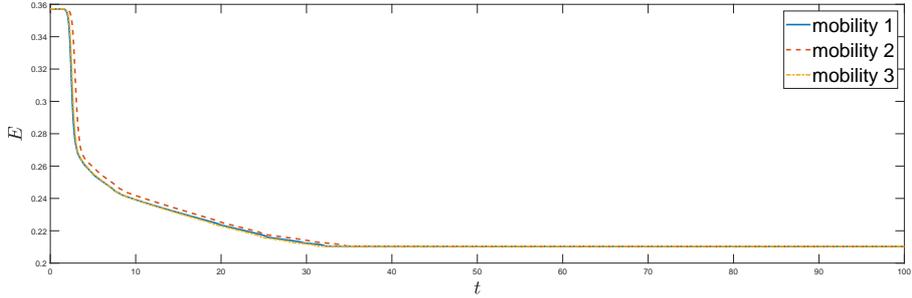}
\\ \vspace{-0.3cm}
\caption{Energies in $t\in[0,100]$ with respect to 3 different mobilities. \label{constant_MCH_mobility_energy} }
\end{figure}

\begin{example}[\textbf{Pattern formation driven by the magnetic field}\label{constant_MCH_mobility}]
We set $N_A=15, ~N_B=10, ~N_S=1,  ~\chi_{AA}=\chi_{BB}=\chi_{SS}=0, ~\chi_{AB}=1, ~\chi_{AS}=2, ~\chi_{BS}=4, ~\varepsilon=0.01, ~\gamma=10^5$ and $\gamma_m=10^{-3}, ~\bm B_0=[1,0]^\top$. We use mobility matrices in \eqref{constant_CH_mobility_matrices} and initial conditions in \eqref{constant_CH_rand_initial}. Some phase portraits and the total energy are plotted in Figures \ref{constant_MCH_mobility_AB}, \ref{constant_MCH_mobility_S} and \ref{constant_MCH_mobility_energy}. Notice the magnetic field is aligned in the $x$-direction in this simulation.
\end{example}
The effect of the magnetic field is similar to that of the electric field in that it induces lamellar stripes in the final state. During transient dynamics, both lamellar structures and defects emerge. Defects move along the direction of the magnetic field and disappear at the boundary of the computational domain. The number of defects eventually decreases to $0$ and the morphology of the lamellar structure shows stripes in the direction of the imposed magnetic field. The energy curves decay to a steady state as shown in Figure \ref{constant_MCH_mobility_energy}.

\subsection{Hysteresis effect related to the electric field}
In this example, we use the electric field coupled model to simulate the hysteresis effect due to the increase and decrease of the strength of an electric field. First, taking the steady state denoted as $S_1$ in Example \ref{constant_CH_spot} as the initial condition, we add an electric field $\bm E_0(t)$ to drive state $S_1$ into a new state $S_2$. Then, we shut down the electric field and we expect structure $S_2$ relaxes back to $S_1$.

\begin{figure}[H]
\centering
\subfigure[$\phi_A-\phi_B$]
{\includegraphics[width=0.7\textwidth]{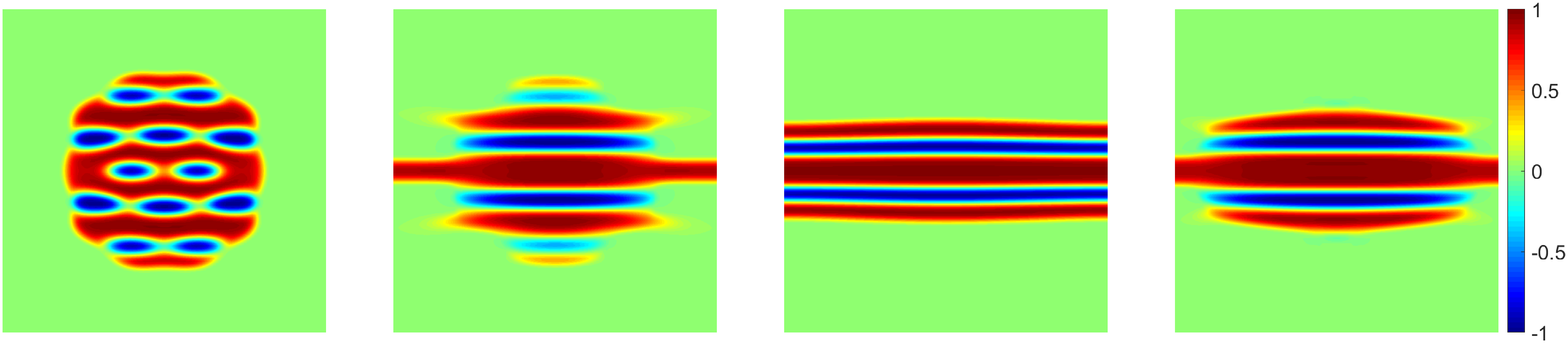}}
\\ \vspace{-0.3cm}
\subfigure[$\phi_S$]{\includegraphics[width=0.7\textwidth]{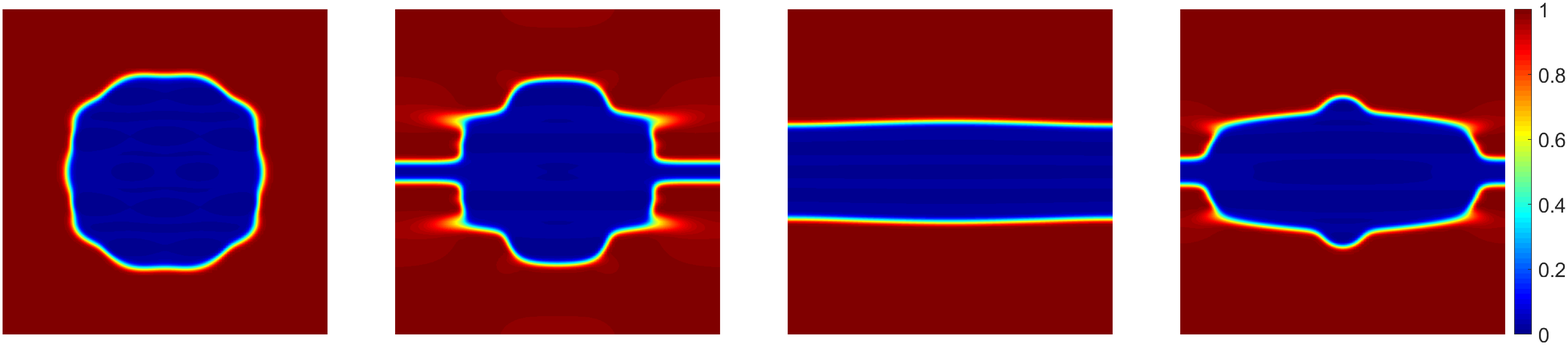}}
\\ \vspace{-0.3cm}
\caption{Microphase and macrophase separation at $t=1.3, ~5, ~15, ~20$ after the electric field is imposed. The electric field stretches the given spot structures into lamellar stripes.   \label{constant_CH_hysteresis_destroy}}
\end{figure}

\begin{figure}[H]
\centering
\subfigure[$\phi_A-\phi_B$]
{\includegraphics[width=0.7\textwidth]{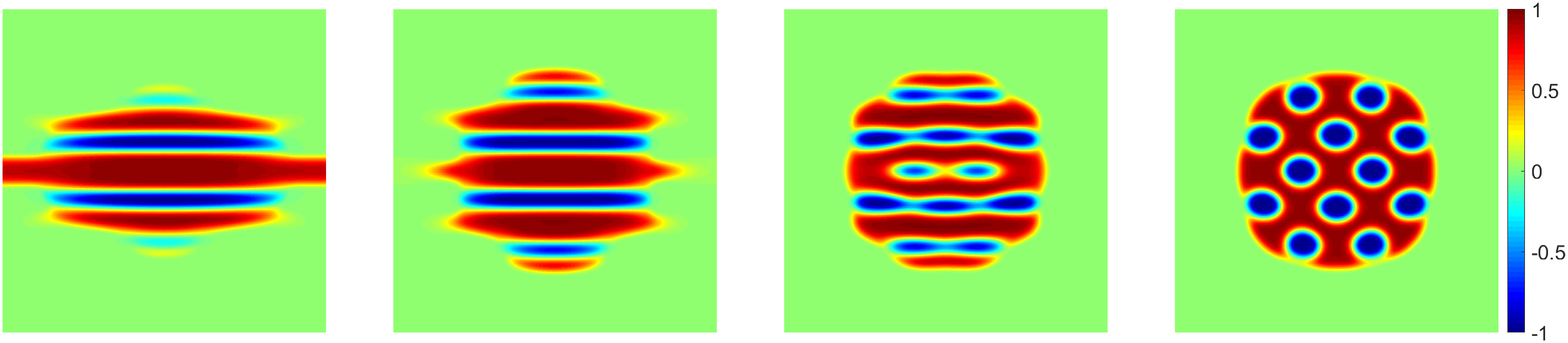}}
\\ \vspace{-0.3cm}
\subfigure[$\phi_S$]{\includegraphics[width=0.7\textwidth]{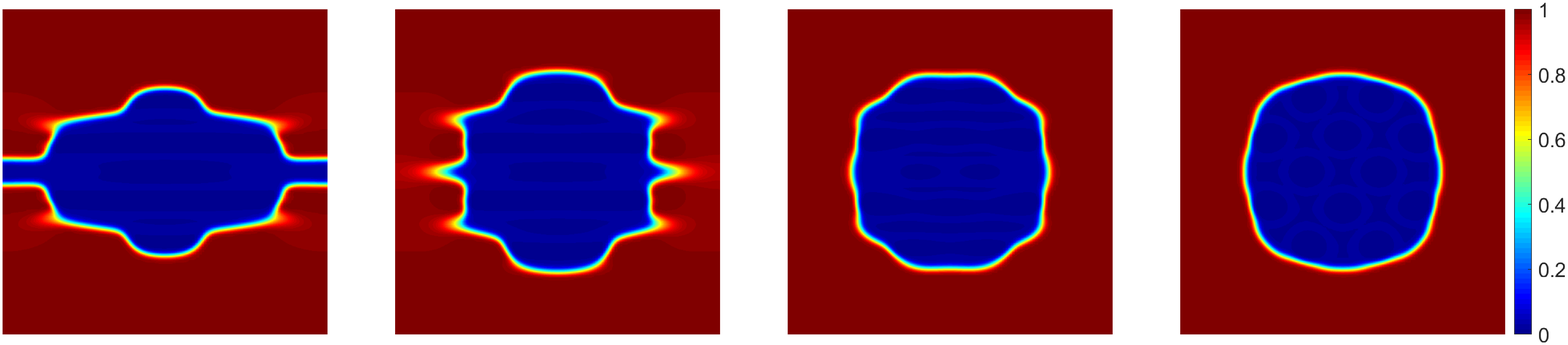}}
\\ \vspace{-0.3cm}
\caption{Microphase and macrophase separation at $t=22, ~25, ~35, ~70$ after the imposed electric field is reversed. Without the electric field, the lamellar stripe structures gradually relax back to spot structures. However, the time it takes for $S_2$ to relax back to $S_1$ is much longer than it evolves from $S_1$ to $S_2$. \label{constant_CH_hysteresis_recover}}
\end{figure}

\begin{figure}[H]
\centering
\subfigure[$-20\leq t\leq70$]{
\includegraphics[width=0.8\textwidth]{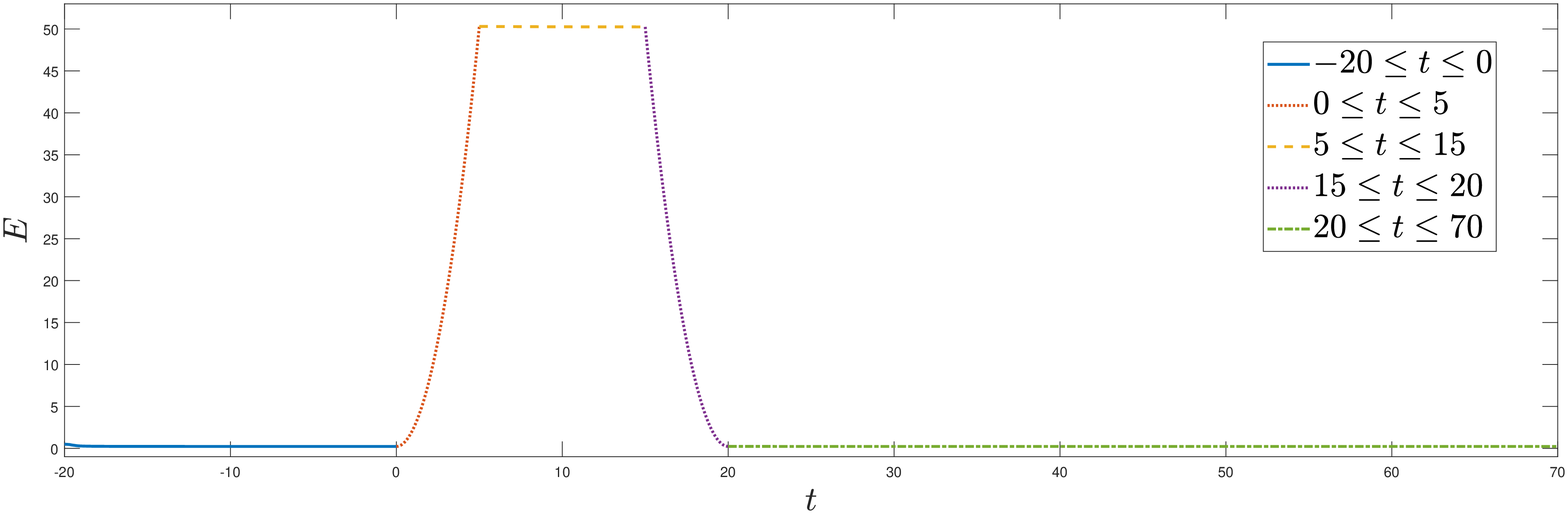}}
\\ \vspace{-0.3cm}
\subfigure[$-20\leq t\leq 0$]{
\includegraphics[width=0.3\textwidth]{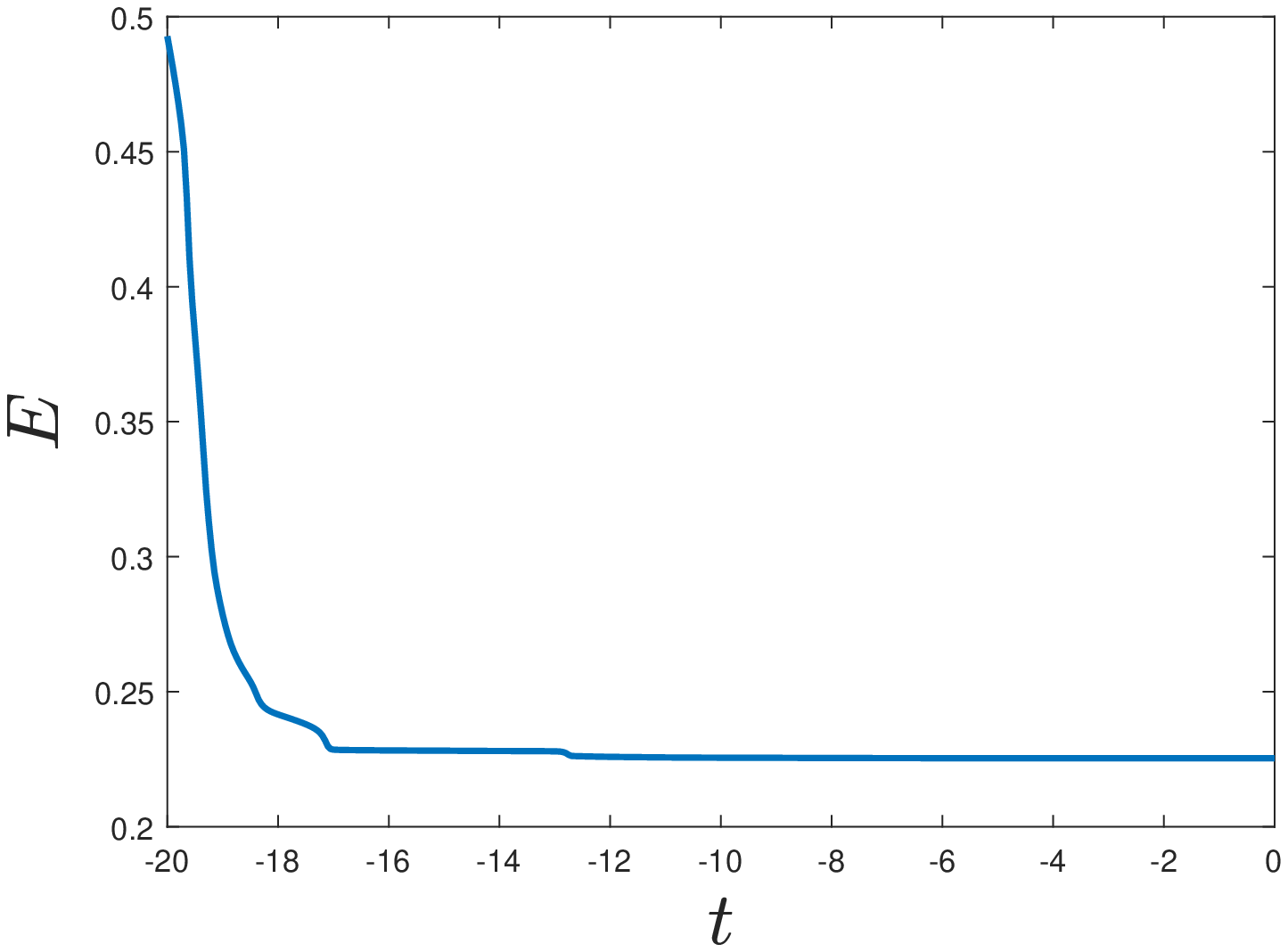}}
\subfigure[$5\leq t\leq 15$]{
\includegraphics[width=0.3\textwidth]{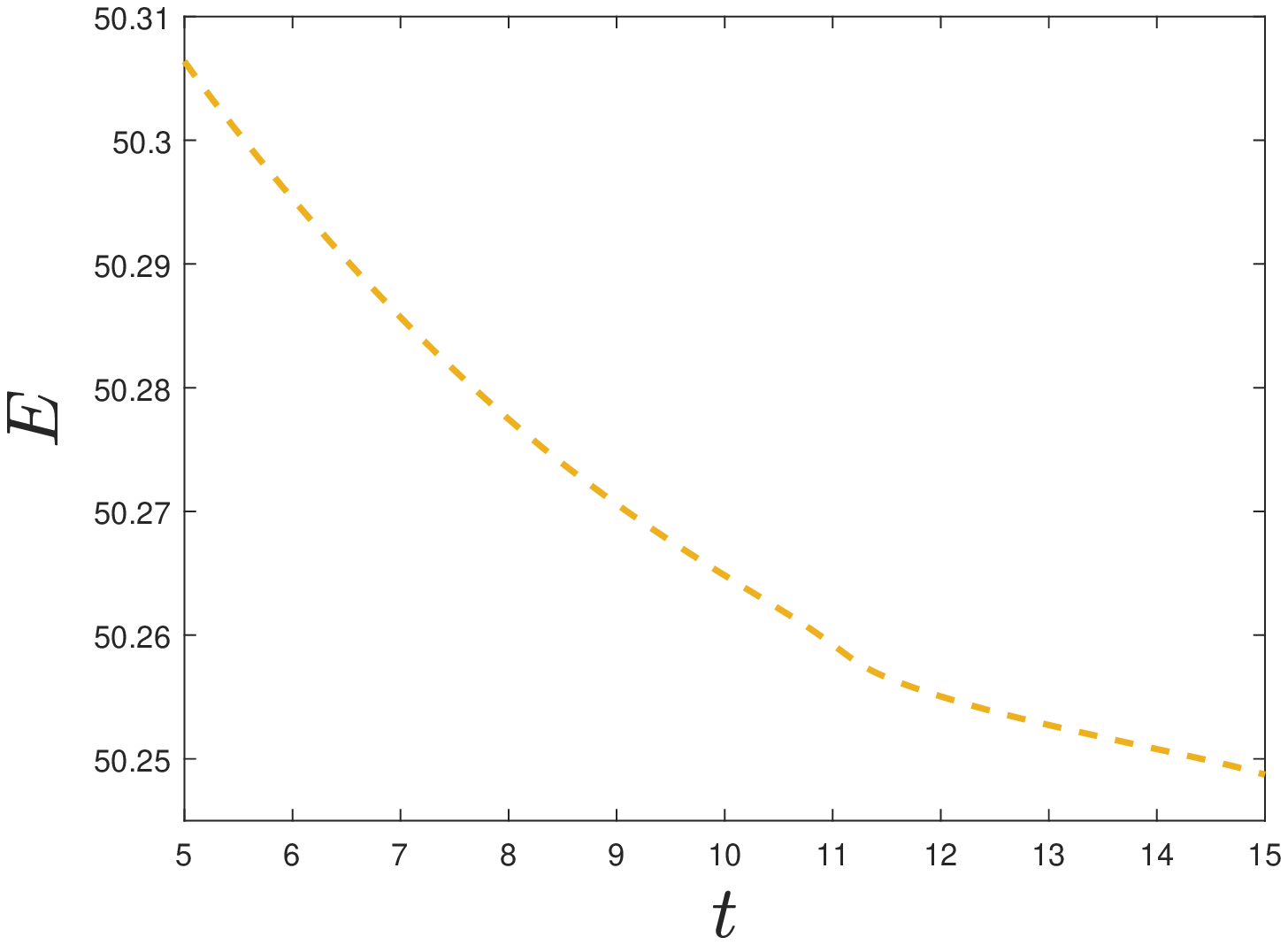}}
\subfigure[$20\leq t\leq 70$]{
\includegraphics[width=0.3\textwidth]{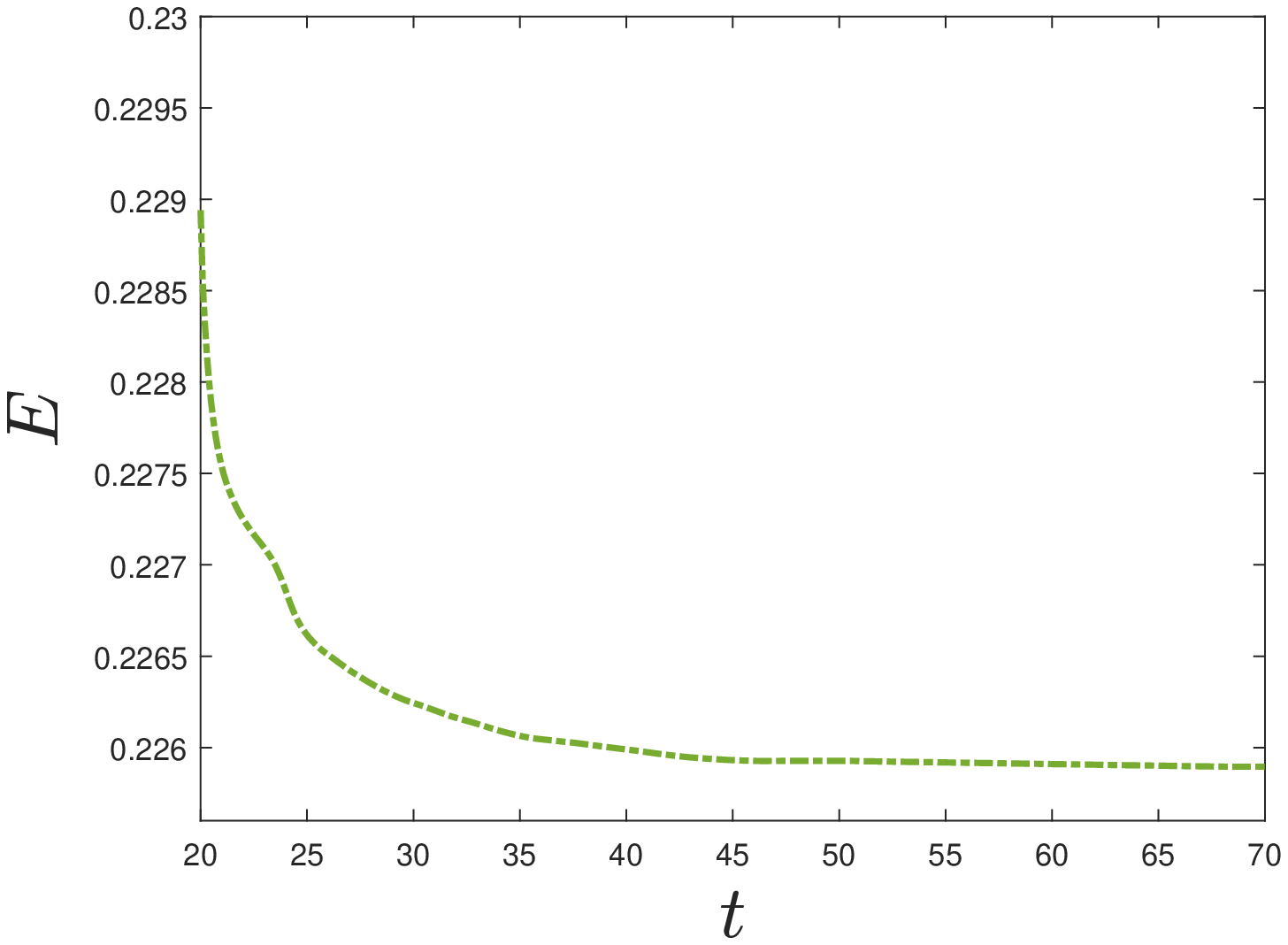}}
\\ \vspace{-0.3cm}
\caption{Energy in $t\in[-20,70]$. We use ``$t\in[-20,0]$'' to substitute expression ``$t\in[0,20]$'' from Example \ref{constant_CH_spot} to display the energy evolution. The result shows that the energy is not symmetric with respect to $t\in[0,20]$.  \label{constant_CH_hysteresis_energy}}
\end{figure}

\begin{figure}[H]
\centering
\subfigure[$\|\bm E_0-\nabla\Phi\|$]{
\includegraphics[width=0.7\textwidth]{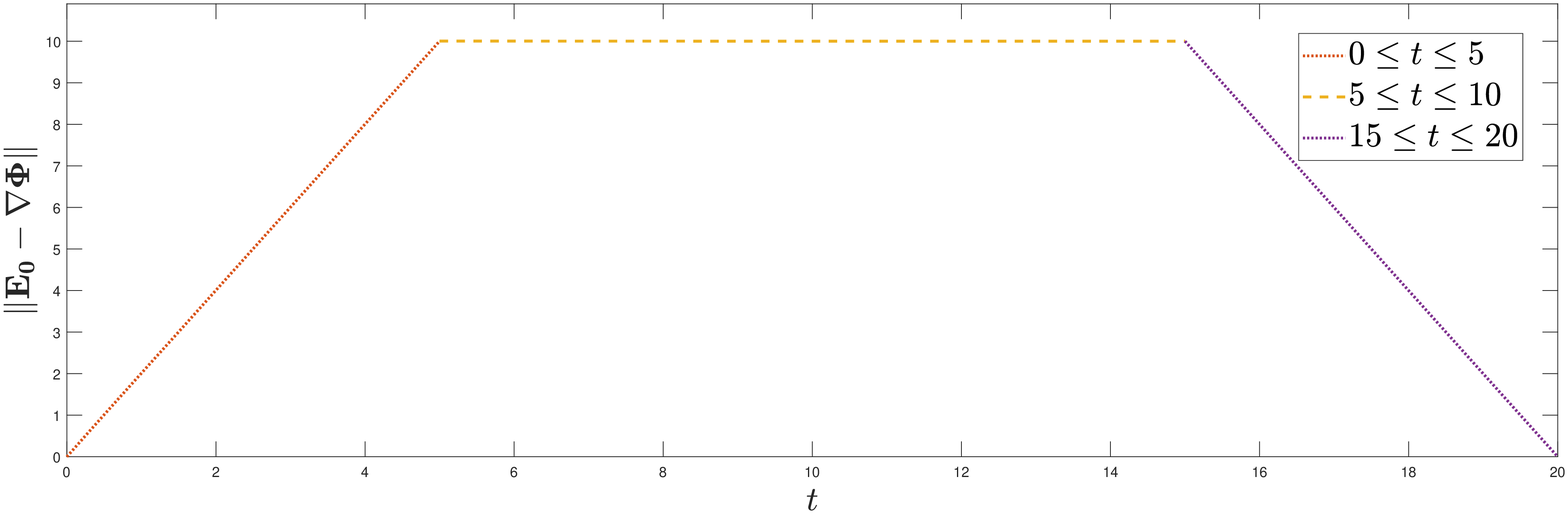}}
\\ \vspace{-0.3cm}
\subfigure[$\|\nabla\Phi\|$]{
\includegraphics[width=0.7\textwidth]{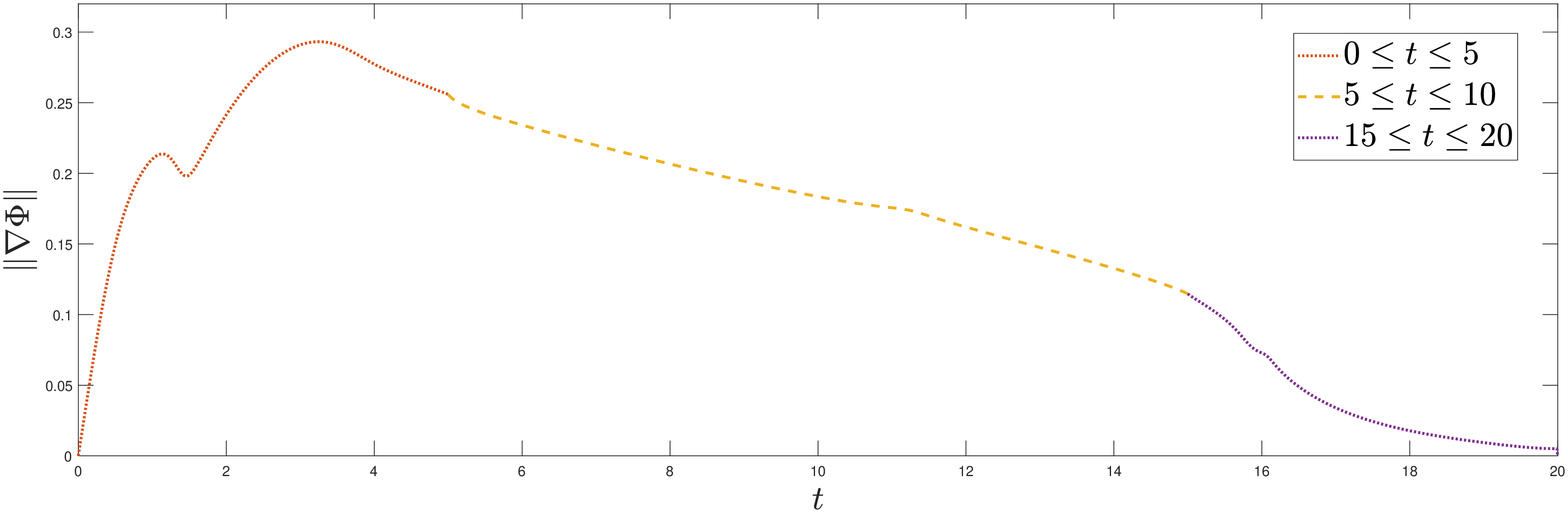}}
\\ \vspace{-0.3cm}
\caption{Magnitudes (discrete $L^2$ norm) of the electric field $\bm E=\bm E_0-\nabla\Phi$ and induced electric field $\nabla\Phi$ in $t\in[0,20]$. The magnitude of induced electric field $\|\nabla\Phi\|$ is small compared to $\|\bm E_0-\nabla\Phi\|$, but not negligible. \label{constant_CH_hysteresis_induced} }
\end{figure}

\begin{example}[\textbf{Hysteresis effect}]
We set $N_A=2, ~N_B=N_S=1, ~\chi_{AA}=\chi_{BB}=\chi_{SS}=0, ~\chi_{AB}=6, ~\chi_{AS}=4, ~\chi_{BS}=8, ~\varepsilon=0.01, ~\gamma=10^3, ~\epsilon_0=1, ~\epsilon_1=0.6$. We use the mobility matrix in \eqref{constant_CH_mobility_diagonal} and take the numerical solution at $t=20$ from Example \ref{constant_CH_spot} as the initial condition of this example. The imposed electric field is given by $\bm E_0(t)=[E_1(t),0]^\top$ where
\begin{equation}
E_1(t) = \left\{
\begin{aligned}
&2t, && 0\leq t\leq 5,
\\
&10, && 5\leq t\leq 15,
\\
&40-2t, && 15\leq t\leq 20,
\\
&0, && t\geq 20.
\end{aligned}
\right.
\end{equation}
Some phase portraits at different time frames and the total energy are plotted in Figures \ref{constant_CH_hysteresis_destroy}, \ref{constant_CH_hysteresis_recover} and \ref{constant_CH_hysteresis_energy}. The magnitudes of the electric field and induced electric field are presented in Figure \ref{constant_CH_hysteresis_induced}.
\end{example}


Notice that the induced electric potential $\Phi$ satisfies \eqref{Maxwell_equation}, i.e.
\begin{equation}
\nabla\cdot\Big[
\big(\epsilon_0+\epsilon_1(\phi_A-\phi_B-\overline\phi_A+\overline\phi_B)\big)(\bm E_0-\nabla\Phi)
\Big]=0 .
\end{equation}
So $\Phi$ implies some mesoscopic structure properties of the system. We observe that curve $\|\bm E_0-\nabla\Phi\|$ in Figure \ref{constant_CH_hysteresis_induced}(a) appears to be symmetric about $t=10$ at large scale, while curve $\|\nabla\Phi\|$ in Figure \ref{constant_CH_hysteresis_induced}(b) does not in fine scale. The asymmetry in the relaxation dynamics of the morphology of transient states and the induced electric field confirm the existence of hysteresis effect with respect to the increase and decrease of the electric field.

\section{Conclusion}
We have derived thermodynamically consistent models for diblock copolymer solutions coupled with the electric and magnetic field, respectively, using the generalized Onsager principle. We then devised a series of structure-preserving numerical algorithms for the models to preserve their energy dissipation rates using the EQ method and the supplementary variable method (SVM), respectively. Like EQ, SVM provides a paradigm for one to develop structure-preserving numerical algorithms that preserve the thermodynamical consistency of the models after discretization. Numerical experiments using the SVM scheme are performed to confirm the convergence rates and reveal the impact of the mobility, the electric and magnetic fields on controlling dynamics of the copolymer solution.

\section*{Acknowledgements}
This research is partially supported by the National Natural Science Foundation of China (award 11971051 and NSAF-U1930402). Qi Wang's research is partially supported by a DOE grant (DE-SC0020272), National Science Foundation grants (award DMS-1815921 and OIA-1655740) and a GEAR award from SC EPSCoR/IDeA Program.

\bibliographystyle{plain}
\bibliography{ref}

\begin{thebibliography}{10}

\bibitem{Badia_2011_EQ}
Santiago Badia, Francisco Guill\'en-Gonz\'alez, and Juan~Vicente
  Guti\'errez-Santacreu.
\newblock Finite element approximation of nematic liquid crystal flows using a
  saddle-point structure.
\newblock {\em J. Comput. Phys.}, 230(41):1686–1706, 2011.

\bibitem{Krausch_2003_Electric_copolymer_solution}
A.~B\"oker, H.~Elbs, H.~H\"ansel, A.~Knoll, S.~Ludwigs, H.~Zettl, A.~V.
  Zvelindovsky, G.~J.~A. Sevink, V.~Urban, V.~Abetz, A.~H.~E. M\"uller, and
  G.~Krausch.
\newblock Electric field induced alignment of concentrated block copolymer
  solutions.
\newblock {\em Macromolecules}, 36(21):8078--8087, 2003.

\bibitem{YangXiaofeng_2021_Triblock_copolymer}
Chuanjun Chen, Xi~Li, Jun Zhang, and Xiaofeng Yang.
\newblock {Efficient linear, decoupled, and unconditionally stable scheme for a
  ternary Cahn-Hilliard type Nakazawa-Ohta phase-field model for tri-block
  copolymers}.
\newblock {\em Appl. Math. Comput.}, 388(1):125463, 2021.

\bibitem{YangXiaofeng_2020_Hydrodynamic_copolymer_melt}
Chuanjun Chen, Jun Zhang, and Xiaofeng Yang.
\newblock {Efficient numerical scheme for a new hydrodynamically-coupled
  conserved Allen–Cahn type Ohta–Kawaski phase-field model for diblock
  copolymer melt}.
\newblock {\em Comput. Phys. Commun.}, 256:107418, 2020.

\bibitem{ShenJie_2020_Lagrange_SAV}
Qing Cheng, Chun Liu, and Jie Shen.
\newblock {A new Lagrange multiplier approach for gradient flows}.
\newblock {\em Comput. Methods Appl. Mech. Engrg.}, 367:113070, 2020.

\bibitem{ShenJie_2017_EQ_copolymer}
Qing Cheng, Xiaofeng Yang, and Jie Shen.
\newblock Efficient and accurate numerical schemes for a hydro-dynamically
  coupled phase field diblock copolymer model.
\newblock {\em J. Comput. Phys.}, 341:44--60, 2017.

\bibitem{Choksi_2011_2D_copolymer_melt}
Rustum Choksi, Mirjana Maras, and J.~F. Williams.
\newblock {2D phase diagram for minimizers of a Cahn–Hilliard functional with
  long-range interactions}.
\newblock {\em Siam J. Appl. Dyn. Syst.}, 10(4):1344–1362, 2011.

\bibitem{Choksi_2009_3D_copolymer_melt}
Rustum Choksi, Mark~A. Peletier, and J.~F. Williams.
\newblock {On the phase diagram for microphase separation of diblock
  copolymers: an approach via a nonlocal Cahn–Hilliard functional}.
\newblock {\em Siam J. Appl. Math.}, 69(6):1712–1738, 2009.

\bibitem{Choksi_2003_Derive_copolymer_melt}
Rustum Choksi and Xiaofeng Ren.
\newblock On the derivation of a density functional theory for microphase
  separation of diblock copolymers.
\newblock {\em J. Stat. Phys.}, 113(1-2):151–176, 2003.

\bibitem{Choksi_2005_Derive_copolymer_homopolymer}
Rustum Choksi and Xiaofeng Ren.
\newblock Diblock copolymer/homopolymer blends: derivation of a density
  functional theory.
\newblock {\em Physica D}, 203(1-2):100--119, 2005.

\bibitem{Copetti_1992_CH_log_energy}
Maria In\^es~Martins Copetti and Charlie~M. Elliott.
\newblock {Numerical analysis of the Cahn-Hilliard equation with a logarithmic
  free energy}.
\newblock {\em Numer. Math.}, 63(1):39--65, 1992.

\bibitem{Faghihi_2018_Magnetic_islands}
Niloufar Faghihi, Simiso Mkhonta, Ken~R. Elder, and Martin Grant.
\newblock Magnetic islands modelled by a phase-field-crystal approach.
\newblock {\em Eur. Phys. J. B}, 91:55, 2018.

\bibitem{GongYuezheng_2020_SVM}
Yuezheng Gong, Qi~Hong, and Qi~Wang.
\newblock Supplementary variable method for developing structure-preserving
  numerical approximations to thermodynamically consistent partial differential
  equations.
\newblock {\em arXiv:2006.04348}, 2020.

\bibitem{Gonzalez_2013_EQ}
Francisco Guill\'en-Gonz\'alez and Giordano Tierra.
\newblock {On linear schemes for a Cahn-Hilliard diffuse interface model}.
\newblock {\em J. Comput. Phys.}, 234:140--171, 2013.

\bibitem{Gonzalez_2014_EQ_CH}
Francisco Guill\'en-Gonz\'alez and Giordano Tierra.
\newblock {Second order schemes and time-step adaptivity for Allen–Cahn and
  Cahn–Hilliard models}.
\newblock {\em Comput. Math. Appl.}, 68(8):821--846, 2013.

\bibitem{Hairer_GNI}
Ernst Hairer, Christian Lubich, and Gerhard Wanner.
\newblock {\em Geometric numerical integration: structure-preserving algorithms
  for ordinary differential equations}.
\newblock Springer, 2006.

\bibitem{Hamley_Copolymer_physics}
Ian~W. Hamley.
\newblock {\em The physics of block copolymers}.
\newblock Oxfor University Press, 1998.

\bibitem{Hamley_Copolymer_science}
Ian~W. Hamley.
\newblock {\em Developments in block copolymer science and technology}.
\newblock John Wiley, 2004.

\bibitem{Hamley_Copolymer_solution}
Ian~W. Hamley.
\newblock {\em Block copolymers in solution: fundamentals and applications}.
\newblock John Wiley, 2005.

\bibitem{HongQi_2020_SVM}
Qi~Hong, Jun Li, and Qi~Wang.
\newblock Supplementary variable method for structure-preserving approximations
  to partial differential equations with deduced equations.
\newblock {\em Appl. Math. Lett.}, 110:106576, 2020.

\bibitem{Landau_Electrodynamics}
Lev~Davidovich Landau and Evgenii~Mikha\'{\i}lovich Lifshitz.
\newblock {\em Electrodynamics of continuous media}.
\newblock Elsevier, 1984.

\bibitem{Makatsoris_2019_3D_diblock_copolymer}
Dung~Q. Ly and Charalampos Makatsoris.
\newblock {Effects of the homopolymer molecular weight on a diblock copolymer
  in a 3D spherical confinement}.
\newblock {\em BMC Chem.}, 13(1):24, 2019.

\bibitem{Ohta_1986_Copolymer_melts}
Takao Ohta and Kyozi Kawasaki.
\newblock Equilibrium morphology of block copolymer melts.
\newblock {\em Macromolecules}, 19(10):2621--2632, 1986.

\bibitem{Ohta_1998_Copolymer_homopolymer}
Takao Ohta and Makiko Nonomura.
\newblock Elastic property of bilayer membrane in copolymer-homopolymer
  mixtures.
\newblock {\em Eur. Phys. J. B}, 2(1):57--68, 1998.

\bibitem{Onsager_1931_RRIP1}
Lars Onsager.
\newblock {Reciprocal relations in irreversible processes. I.}
\newblock {\em Phys. Rev.}, 37:405--426, 1931.

\bibitem{Onsager_1931_RRIP2}
Lars Onsager.
\newblock {Reciprocal relations in irreversible processes. II.}
\newblock {\em Phys. Rev.}, 38:2265--2279, 1931.

\bibitem{OnsagerMachlup_1953_Irreversible}
Lars Onsager and Stefan Machlup.
\newblock Fluctuations and irreversible processes.
\newblock {\em Phys. Rev.}, 91:1505--1512, 1953.

\bibitem{Glasner_2016_Electric_copolymer_melt}
Saulo Orizaga and Karl Glasner.
\newblock Instability and reorientation of block copolymer microstructure by
  imposed electric fields.
\newblock {\em Phys. Rev. E}, 93(5):052504, 2016.

\bibitem{Boker_2017_Electric_copolymer}
Christian~W. Pestera, Clemens Liedel, Markus Ruppel, and Alexander B\"oker.
\newblock Block copolymers in electric fields.
\newblock {\em Prog. Polym. Sci.}, 64:182--214, 2017.

\bibitem{Zvelindovsky_2009_Electric_copolymer_lamellae}
Marco Pinna, Ludwig Schreier, and Andrei~V. Zvelindovsky.
\newblock Mechanisms of electric-field-induced alignment of block copolymer
  lamellae.
\newblock {\em Soft Matter}, 5:970--973, 2009.

\bibitem{Seymour_2015_Phase_field_magnet_electric}
Matthew Seymour, F.~Sanches, Ken Elde, and Nikolas Provatas.
\newblock Phase-field crystal approach for modeling the role of microstructure
  in multiferroic composite materials.
\newblock {\em Phys. Rev. E}, 92:184109, 2015.

\bibitem{ShenJie_2018_SAV}
Jie Shen, Jie Xu, and Jiang Yang.
\newblock {The scalar auxiliary variable (SAV) approach for gradient flows}.
\newblock {\em J. Comput. Phys.}, 353:407--416, 2018.

\bibitem{ShenJie_2019_SAV}
Jie Shen, Jie Xu, and Jiang Yang.
\newblock A new class of efficient and robust energy stable schemes for
  gradient flows.
\newblock {\em Siam Rev.}, 61(3):474--506, 2019.

\bibitem{Choksi_2015_Metastable_states}
David Shirokoff, Rustum Choksi, and Jean-Christophe Nave.
\newblock Sufficient conditions for global minimality of metastable states in a
  class of non-convex functionals: a simple approach via quadratic lower
  bounds.
\newblock {\em J. Nonlinear Sci.}, 25(3):539–582, 2015.

\bibitem{Williams_2017_2D_Ohta-Kawasaki}
Jan~Bouwe van~den Berg and J.~F. Williams.
\newblock {Validation of the bifurcation diagram in the 2D Ohta–Kawasaki
  problem}.
\newblock {\em Nonlinearty}, 30(4):1584–1638, 2017.

\bibitem{Williams_2019_3D_Ohta-Kawasaki}
Jan~Bouwe van~den Berg and J.~F. Williams.
\newblock {Rigorously computing symmetric stationary states of the
  Ohta-Kawasaki problem in three dimensions}.
\newblock {\em Siam J. Math. Aanl.}, 51(1):131--158, 2019.

\bibitem{WangQi_2021_GOP_application}
Qi~Wang.
\newblock Generalized onsager principle and its applications.
\newblock In Xiangyang Liu, editor, {\em Frontiers and progress of current soft
  matter research}, chapter~3. Springer Singapore, 2021.

\bibitem{LiBaohui_2015_Double_hydrophilic_copolymer}
Jiaping Wu, Zheng Wang, Yuhua Yin, Run Jiang, Baohui Li, and An~chang Shi.
\newblock A simulation study of phase behavior of double-hydrophilic block
  copolymers in aqueous solutions.
\newblock {\em Macromolecules}, 48(24):8897–8906, 2014.

\bibitem{WuXiangfa_2008_Copolymer_with_electric}
Xiangfa Wu and Yuris~A. Dzenis.
\newblock Phase-field modeling of the formation of lamellar nanostructures in
  diblock copolymer thin films under inplanar electric fields.
\newblock {\em Phys. Rev. E}, 77(511):031807, 2008.

\bibitem{Russel_2004_Electric_copolymer}
Ting Xu, A.~V. Zvelindovsky, G.~J.~A. Sevink, Oleg Gang, Ben Ocko, Yuqing Zhu,
  Samuel~P. Gido, and Thomas~P. Russell.
\newblock Electric field induced sphere-to-cylinder transition in diblock
  copolymer thin films.
\newblock {\em Macromolecules}, 37(18):6980--6984, 2004.

\bibitem{Russell_2005_Electric_copolymer}
Ting Xu, A.~V. Zvelindovsky, G.~J.~A. Sevink, K.~S. Lyakhova, H.~Jinnai, and
  Thomas~P. Russell.
\newblock Electric field alignment of asymmetric diblock copolymer thin films.
\newblock {\em Macromolecules}, 38(26):10788--10798, 2005.

\bibitem{Yamada_2014_Copolymer_solvent}
Kohtaro Yamada, Emiko Yasuno~Youhei Kawabata, and Tohru Okuzono~Tadashi Kato.
\newblock Mesoscopic simulation of phase behaviors and structures in an
  amphiphile-solvent system.
\newblock {\em Phys. Rev. E}, 89(6):062310, 2014.

\bibitem{YangXiaofeng_2019_CH_log_energy}
Xiaofeng Yang and Jia Zhao.
\newblock {On linear and unconditionally energy stable algorithms for variable
  mobility Cahn-Hilliard equation with logarithmic Flory-Huggins potential}.
\newblock {\em Comm. Comp. Phys.}, 25:703--728, 2019.

\bibitem{WangQi_2016_GOP}
Xiaogang Yang, Jun Li, M.~Gregory Forest, and Qi~Wang.
\newblock {Hydrodynamic theories for flows of active liquid crystals and the
  Generalized Onsager principle}.
\newblock {\em Entropy}, 18(6):202, 2016.

\bibitem{DongSuchuan_2020_gPAV}
Zhiguo Yang and Suchuan Dong.
\newblock A roadmap for discretely energy-stable schemes for dissipative
  systems based on a generalized auxiliary variable with guaranteed positivity.
\newblock {\em J. Comput. Phys.}, 404:109121, 2020.

\bibitem{ZhangLei_2020_Energy_landscape}
Jianyuan Yin, Yiwei Wang, Jeff Z.~Y. Chen, Pingwen Zhang, and Lei Zhang.
\newblock Construction of a pathway map on a complicated energy landscape.
\newblock {\em Phys. Rev. Lett.}, 124(9):090601, 2020.

\bibitem{ZhangLei_2020_Find_saddle_points}
Bing Yu and Lei Zhang.
\newblock Global optimization-based dimer method for finding saddle points.
\newblock {\em Discrete Cont. Dyn.-B}, 26(1):741--753, 2021.

\bibitem{ChenChuanjun_2020_Copolymer_melt}
Jun Zhang, Chuanjun Chen, and Xiaofeng Yang.
\newblock {Efficient and energy stable method for the Cahn-Hilliard phase-field
  model for diblock copolymers}.
\newblock {\em Appl. Numer. Math.}, 151:263--281, 2020.

\bibitem{ChenChuanjun_2020_AC_copolymer}
Jun Zhang, Chuanjun Chen, Xiaofeng Yang, and Kejia Pan.
\newblock {Efficient numerical scheme for a penalized Allen–Cahn type
  Ohta–Kawasaki phase-field model for diblock copolymers}.
\newblock {\em J. Comput. Appl. Math.}, 378:112905, 2020.

\bibitem{YangXiaofeng_2020_Magnetic_copolymer}
Jun Zhang and Xiaofeng Yang.
\newblock {A new magnetic-coupled Cahn–Hilliard phase-field model for diblock
  copolymers and its numerical approximations}.
\newblock {\em Appl. Math. Lett.}, 107:106412, 2020.

\bibitem{WangQi_2018_EQ}
Jia Zhao, Xiaofeng Yang, Yuezheg Gong, Xueping Zhao, Xiaogang Yang, Jun Li, and
  Qi~Wang.
\newblock A general strategy for numerical approximations of non-equilibrium
  models–part i: thermodynamical systems.
\newblock {\em Int. J. Numer. Anal. Mod.}, 15(6):884--918, 2018.

\end{thebibliography}
\addcontentsline{toc}{section}{Bibliography}
\end{document}